\documentstyle[12pt]{article}
\def\m@th{\mathsurround=0pt }
\def\eqalign#1{\null\,\vcenter{\openup\jot \m@th
   \ialign{\strut\hfil$\displaystyle{##}$&$
      \displaystyle{{}##}$\hfil \crcr#1\crcr}}\,}
\newtheorem{theorem}{Theorem}[section]
\newtheorem{proposition}[theorem]{Proposition}
\newtheorem{corollary}[theorem]{Corollary}
\newtheorem{lemma}[theorem]{Lemma}

\def\ad{{\rm ad }\ }
\def\gra{{\rm gr }\ }
\def\adr{{\rm ad}_r}

\font\goth=eufm10

\def\gg{\hbox{\goth g}}

\def\gn{\hbox{\goth n}}

\def\gh{\hbox{\goth h}}

\def\gs{\hbox{\goth s}}

\def\gl{\hbox{\goth l}}
\def\go{\hbox{\goth o}}
\def\gr{\hbox{\goth r}}

\begin{document} 

\begin{center}
{\large{\bf  Quantum  Zonal Spherical  
Functions and Macdonald Polynomials } }
\vskip 5mm Gail Letzter\footnote{supported by NSA grant no. 
MDA904-01-1-0033.

AMS subject classification 17B37}

Mathematics Department

Virginia Polytechnic Institute

and State University

Blacksburg, VA 24061

letzter@math.vt.edu

\end{center}
\medskip
\begin{abstract}
A unified theory of quantum symmetric pairs is applied to $q$-special functions.  Previous work
characterized 
certain left coideal subalgebras in the 
quantized enveloping algebra and established an appropriate framework for 
quantum zonal spherical functions.  Here 
a distinguished family of such functions, invariant under the Weyl 
group associated to the restricted roots, is shown to be a family of 
Macdonald polynomials, as conjectured by Koornwinder and Macdonald. 
Our results place earlier work for  Lie algebras  of 
classical type  in a  general context and extend to the 
exceptional cases.
\end{abstract}

\centerline{INTRODUCTION}
\bigskip
A beautiful classical result shows that the zonal spherical functions associated to 
 real compact symmetric spaces can be realized as Jacobi polynomials.
An analogous result was proved for zonal spherical functions of  p-adic symmetric spaces
([M1]).
In the late 1980's, I.G. Macdonald [M2] introduced his  two-parameter family of orthogonal 
polynomials which  provided a  unified context for the  polynomial families used in
these parallel   theories. 
 With the discovery of quantum groups, also in the 
1980's,  both Koornwinder [K2, Section 6.4] and Macdonald [M2] 
asked whether one could obtain a  similar 
description of  quantum zonal spherical functions. This conjecture 
was investigated extensively in the 1990's  by Noumi, 
Sugitani, Dijkhuisen,  Stokman and others when the underlying Lie algebra is 
of classical type. 
Quantum zonal spherical functions are realized as particular $q$ special functions 
in [N], [NS], [S], [DN], [NDS], and [DS]
using  case-by-case 
computations.
The main result of
this paper  is an answer to  Macdonald and Koornwinder's question for all quantum symmetric pairs
with reduced restricted root systems, thus generalizing this earlier 
work to include Lie algebras of exceptional type.   Our methods 
are new and do not involve extensive case work.  Instead, we  rely on the 
unified theory  of quantum symmetric pairs developed in   [L2], [L3], and 
[L4].  As a consequence, we provide   simple, uniform, formulas for  
the parameters which appear in the Macdonald polynomials corresponding 
to quantum zonal spherical functions.

 The first problem faced in the quantum case was the actual definition of quantum
symmetric spaces. Classically, a symmetric pair of Lie algebras is a  pair
$\gg,\gg^{\theta}$ where $\gg$ is a complex Lie algebra and ${\theta}$ is a Lie algebra 
involution. The initial breakthrough  was made by Koornwinder [K1], who constructed
quantum
$2$ spheres  and computed their zonal spherical functions by inventing quantum analogs of
the symmetric pair $\gs\gl \ 2,\gs\go\ 2$.  Noumi ([N]) extended Koornwinder's approach to
two
 families of symmetric pairs for $\gg$  of type $A_n$,  using two-sided coideal 
analogs of $\gg^{\theta}$  inside  the 
quantized enveloping algebra $U_q(\gg)$. He showed that the 
corresponding zonal spherical functions were, indeed, Macdonald polynomials. 
In [NS], Noumi and Sugitani introduced one-sided coideal analogs of 
$U(\gg^{\theta})$ inside of $U_q(\gg)$ for other cases when $\gg$ is of 
classical type; analysis of the corresponding zonal spherical functions can be 
found in [NS], [S], [DN], [NDS], and [DS]. In his  comprehensive survey of the early history of quantum symmetric
spaces,  Dijkhuisen ([D]) conjectures that 
 coideals are the  correct objects to use in order to develop a
 general theory of quantum  symmetric pairs, the corresponding symmetric spaces,
and their zonal  spherical functions.
 In [L2] and [L3],  a universal method was developed for constructing left coideal 
subalgebras of $U_q(\gg)$, which are quantum analogs of $U(\gg^{\theta})$. These analogs
are further characterized as the appropriately unique maximal left coideal subalgebras of $U_q(\gg)$
which specialize  to $U(\gg^{\theta})$ as $q$ goes to $1$.
A complete list of the generators and relations of the coideal subalgebras
for all possible symmetric pairs $\gg, \gg^{\theta}$ is  given in [L4].  Furthermore, 
although  the
symmetric pairs used in [N] and [NS] are formally defined differently, it turns out
that these examples are subsumed by this new theory     (see [L2, Section 6] and the last paragraph of [L3]).

Using the new quantum analogs of $U(\gg^{\theta})$ inside of $U_q(\gg)$, 
one can define and study quantum symmetric spaces and their zonal 
spherical functions (see [L3, Section 7] and [L4].) 
In particular, fix an irreducible 
symmetric pair $\gg,\gg^{\theta}$.  Let ${\cal B}_{\theta}$ denote the orbit of the  analogs of
$U(\gg^{\theta})$ in $U_q(\gg)$ under the group ${\bf H}$ of Hopf algebra
automorphisms of $U_q(\gg)$ fixing the Cartan part. For each pair $B,B'$ in ${\cal B}$, one can
associate the space ${}_{B'}{\cal H}_{B}$ of left $B'$ and right $B$ invariants inside the
quantized function algebra $R_q[G]$ corresponding to $\gg$.  Using the
interpretation of elements of 
$R_q[G]$  as functions on $U_q(\gg)$,   restriction  to the Cartan part of $U_q(\gg)$
yields an algebra homomorphism $\Upsilon$ from 
${}_{B'}{\cal H}_{B}$ into the character group ring associated to the restricted root
system ([L4, Theorem 4.2]). Zonal spherical functions
 are joint eigenvectors in ${}_{B'}{\cal H}_{B}$ with respect to the
action of the center of $U_q(\gg)$.
A zonal spherical family associated to ${\cal B}_{\theta}$  is the image of a 
specially chosen basis of ${}_{B'}{\cal H}_{B}$ consisting
 of zonal spherical functions
for some $B,B'$ in ${\cal B}_{\theta}$. Let 
$W_{\Theta}$ denote the Weyl group associated to  the restricted 
roots. In [L4,  Theorem 6.5], it is shown that each ${\bf H}\times {\bf H}$ orbit 
of ${\cal B}_{\theta}\times {\cal B}_{\theta}$ contains a virtually unique $W_{\Theta}$
invariant zonal spherical family.
Moreover, the natural map from
${\cal B}_{\theta}\times {\cal B}_{\theta}$ to the set of zonal spherical functions
associated to
${\cal B}_{\theta}\times {\cal B}_{\theta}$ is 
${\bf H}\times {\bf H}$ equivariant ([L4, Theorem 6.3]).  In this context, the action of an
element of
${\bf H}\times {\bf H}$ on a zonal spherical family corresponds to the image of this
family under an automorphism of the restricted character group ring.   Thus in order to
compute the zonal spherical families associated to ${\cal B}_{\theta}$, it is only
necessary to analyze the unique $W_{\Theta}$ invariant  representative for each   ${\bf
H}\times {\bf H}$ orbit.

The set  ${\cal B}_{\theta}$ contains a distinguished ${\bf H}$ orbit
corresponding to the standard analogs of $U(\gg^{\theta})$ in $U_q(\gg)$. 
It follows from [L4, Theorem 6.5 and subsequent discussion] (see Theorem 1.1 below)
that there is a unique family of $W_{\Theta}$ invariant 
quantum zonal spherical functions associated to the
standard analogs of
$U(\gg^{\theta})$.  We further assume
here that  the  restricted root space associated to
$\gg,\gg^{\theta}$ is reduced.    Our goal in this work is to show  that this family of  
$W_{\Theta}$ invariant
quantum zonal spherical functions is precisely a family  of   Macdonald  polynomials
where both parameters are equal to  powers of 
$q$.  Moreover, Theorem 8.2 provides  a simple formula involving the multiplicity
of restricted roots which relates one parameter to another. 
(We should mention that for most pairs 
$\gg,\gg^{\theta}$, the standard analogs are the only possible analogs.  
However, under special circumstances, there is a one parameter family 
of analogs of $U(\gg^{\theta})$; the standard analogs appear when 
 the parameter is set  equal to zero. See [L3, Section 7, Variation 2] 
for more information.)

Our overall method of determining the zonal spherical functions is 
inspired by the work of [N], [NS], and [D].   Indeed, we show that the radial 
components of certain central elements correspond to the  difference operators, which 
then define  Macdonald polynomials.  However, the strategy we employ is 
quite  different.  We do not use 
  $L$ operators as in [N], [NS], and [D]  to  express the Casimir 
elements   because, in part, this cannot be done for the exceptional Lie 
algebras. Moreover, our argument  avoids an explicit   
expression of central elements inside
$U_q(\gg)$ for
$\gg$ larger  than ${\gs\gl\ 2}$,   thus  making a detour around difficult
computations and extensive case work. Instead, we rely on a variety of 
representation theory techniques and draw upon the description of the center and the locally finite
part of
$U_q(\gg)$ developed in [JL1] and [JL2].

It should be noted that many of the techniques of this paper extend  to the
case when $\Sigma$ is not reduced  and to the nonstandard analogs.  
In a
future paper, we  adapt the methods used here to identify the
zonal spherical functions of these other cases as $q$ hypergeometric
functions.      
Ultimately, this will generalize work of [NDS] and [DS] on
nonstandard analogs of type AIII and work of [NS], [S],[DN],[NDS], and [DS] 
when the restricted root system is not reduced.

Since this paper is long, we present a detailed description of its 
organization.
The first section
sets notation  and recalls background on the author's theory of  quantum symmetric pairs and their zonal spherical 
functions. 
Section 2 
presents four versions of an Iwasawa type tensor product decomposition (Theorem 2.2).
Let $B$ denote a quantum analog of $U(\gg^{\theta})$ inside $U_q(\gg)$ and write
$B_+$ for the augmentation ideal of $U_q(\gg)$. Further analysis of projections of
$U_q(\gg)$ modulo
$B_+U_q(\gg)$ and
$U_q(\gg)B_+$ are obtained in Theorem 2.3. These maps are used in Section 3 (Theorem 3.2)
to construct a function ${\cal X}$ from
$U_q(\gg)$ to the ring of endomorphisms of the restricted character group ring which
corresponds to the action of $U_q(\gg)$ on the quantum zonal spherical functions. Upon
restriction to the center of $U_q(\gg)$, ${\cal X}$ is exactly the function which computes
the  quantum radial components.  Radial components of central elements, and more
generally, radial components of elements in the centralizer of $B$ in $U_q(\gg)$, are
shown to be invariant under the action of the restricted Weyl group
$W_{\Theta}$ in Theorem 3.4 and Theorem 3.6.

  Section 4 is devoted  to analyzing the rank one radial components of 
central  elements.  We first prove that in the rank one case, the center of $U_q(\gg)$ has
a Casimir-like central element which looks remarkably like the Casimir element in the
center of
$U_q(\gs\gl\ 2)$ (Theorem 4.5). The radial component of this special central element is
then computed using a straightforward $U_q(\gs\gl\ 2)$ calculation (Theorem 4.7). 

In Section 5,  a filtration, similar to the ad invariant filtration of
$U_q(\gg)$ studied in [JL2] and [J1, Chapter 7], is introduced. The filtration used here
is carefully chosen so that the algebra $B$
lies in degree zero. In particular, the resulting graded algebra, $\gra U_q(\gg)$,
contains   
$B$ as a subalgebra.   A theory of graded zonal spherical functions    is  presented. It is shown that up to a shift of weight,
the graded zonal spherical functions are all equal to each other (Lemma 5.7). Furthermore, the 
action of the graded image of a central element of $U_q(\gg)$ on a graded zonal spherical 
function agrees with  the action of the top degree of the 
corresponding radial component.  This allows us to write the top degree terms of  radial
components in a simple form using  the graded zonal spherical functions and
the Harish-Chandra projection (Theorem  5.8).  The resulting expression is reminiscent of
Harish-Chandra's formula and its generalizations for the  radial components of central
elements in
$U(\gg)$ with respect to the adjoint action of the corresponding Lie group on $\gg$ (see
[W, 7.A.2.9 and 7.A.3.7]).

In Section 6, the graded zonal spherical functions, and thus the 
top degree terms of radial components, are  determined  first in the 
rank one case using 
 Section 4 (see Lemma 6.6). This information is then glued 
 together using the 
 Weyl invariance of radial components in Theorem 6.7 and Corollary 
 6.8. The argument is delicate since the graded zonal spherical 
 functions are  elements of the  formal power series ring 
 corresponding to the restricted character group ring. 
 Unfortunately, the restricted Weyl group does not act on this ring.  
We overcome this obstacle  by finding the possible weights of  highest weight vectors
with respect to rank one subalgebras of $\gra U_q(\gg)$ inside 
Verma-like $\gra U_q(\gg)$ modules.  This allows one to express   
 the graded zonal spherical function as the product of all the rank one 
 formulas times a $W_{\Theta}$ invariant term.  Taking into account 
 the highest weight summand of the graded zonal spherical function 
 yields that this $W_{\Theta}$ invariant term must be $1$.  In the 
 end, we show that the formula for the graded zonal spherical function 
 is just the inverse of the element in the restricted character group 
 ring used to define Macdonald's inner product (Theorem 6.7).  It 
 should be noted that the 
 proof here is remarkably similar to an argument used in a completely 
 different setting, the factorization of the affine PRV determinant ([JL3]).
 
It is sometimes necessary to pass to the slightly larger  simply connected quantized
enveloping algebra $\check U$ (which is just a small extension of
$U_q(\gg)$). By construction, the Cartan subalgebra of $\check U$ is just the group algebra
of the torus corresponding to the weight lattice. In Section 7, one finds an element in
the centralizer $C_{\check U}(B)$ of $B$ in $\check U$ whose top degree term with respect
to the Harish-Chandra projection corresponds to a minuscule or pseudominuscule restricted
weight.   In most cases this element is in the center of $\check U$ and is easily
determined using   the description of the image of
central elements under the  Harish-Chandra map ([JL1]).  For a few
exceptional irreducible symmetric pairs, the central elements turn out to be too ``large" to
correspond to minuscule or pseudominuscule restricted weights.  This failure is 
related to the classical fact that the map from the center of 
of $U(\gg)$ to the set of invariant differential operators on the 
symmetric space corresponding to the pair $\gg,\gg^{\theta}$ is
not always surjective ([H]). An
analysis of  the locally finite part of $\check U$ is used in order     to locate the suitable
element in
$C_{\check U}(B)$ in the problematic cases.  The radial  components of these 
``small''   elements in $C_{\check U}(B)$ are then determined (Theorem 
7.7) by taking the
sum of the terms in   
the $W_{\Theta}$ orbit of the top degree part (found in Section 6) plus a possible zero
degree term.   

Section  8 recalls basic facts about Macdonald polynomials.  The radial components
studied in Section 7 are identified with difference  operators associated to minuscule and
pseudominuscule weights.  This  in turn establishes our main result that zonal spherical
functions are particular  Macdonald polynomials. We conclude the paper with an appendix 
listing all irreducible symmetric pairs $\gg,\gg^{\theta}$ with 
reduced root system and the
values of the parameters in the Macdonald polynomials corresponding 
to the quantum zonal spherical functions. 

\section{Background and Notation}
Let ${\bf C}$ denote the complex numbers, ${\bf Q}$ denote the rational numbers,
${\bf Z}$ denote the integers, and ${\bf N}$ denote the nonnegative integers.
If $G$ is a multiplicative monoid and ${\bf F}$ is a field, then we write ${\bf F}[G]$ for  the 
corresponding monoid algebra
 over ${\bf F}$. (This is the obvious generalization of ``group 
 algebra''.) Unfortunately many   monoids come to us additively. In the special case of the additive monoid
${\bf Q}$, we invent the symbol $q$ and temporarily identify ${\bf Q}$ with the multiplicative monoid $\{q^{r}|r\in {\bf
Q}\}$ where
$q^rq^s=q^{r+s}$.  Let ${\cal C}$ be the algebraic closure of the field of fractions for ${\bf C}[{\bf Q}]$.
  We write ${\cal C}^{\times}$ for the 
nonzero elements of ${\cal C}$.

Given a root system
$\Phi$, let $\Phi^+$ denote the positive roots,
$Q(\Phi)$ denote the root lattice,
$P(\Phi)$ denote the weight lattice, $Q^+(\Phi)$ denote the ${\bf N}$ span of the
elements in $\Phi^+$, and  $P^+(\Phi)$ denote the set of
dominant integral weights.  Write $2\Phi$ for the root system 
$\{2\alpha|\alpha\in \Phi\}$ given the same inner product as 
$\Phi$.  
If $H$ is an {\em additive} submonoid of ${\bf C}\Phi$, we invent the formal variable $z$,
so that ${\cal C}[H]$ consists of the ${\cal C}$ linear combinations of the basis elements $z^{\lambda}$
for $\lambda\in H$.   

Let
$\gg$ denote a semisimple Lie algebra over the complex numbers
${\bf C}$ with a chosen triangular decomposition $\gg=\gn^-\oplus \gh\oplus\gn^+$.  Let
$\Delta$ be the set of roots  for
$\gg$ and let $\pi=\{\alpha_1,\dots,\alpha_n\}$ be the set of (positive) simple roots in
$\Delta$ corresponding to the root vectors in $\gn^+$.  Write $(\ ,\ )$ for the Cartan inner product
on $\gh^*$ with respect to the root system $\Delta$.  Let $W$ denote
the Weyl group associated to the root system $\Delta$. Set $\rho$ equal to   the half sum of the
positive roots in $\Delta$.  

Let $\theta$ be a maximally split involution with respect to the fixed Cartan subalgebra
$\gh$ of $\gg$ and  triangular decomposition  (see [L3, (7.1), (7.2), and
(7.3)].)  We assume throughout the paper that $\gg,\gg^{\theta}$ is an irreducible pair
in the sense of [A] (see also [L4, Section 7]).
 Recall
that
a complex semisimple Lie algebra $\gg'$ with maximally split involution $\theta'$  can
be written as a direct sum of semisimple Lie subalgebras
$\gr_i$ such that each 
$\gr_i,\gr_i^{\theta'}$ is an irreducible pair. Using such a  direct sum decomposition, the
results of this paper easily extend to arbitrary symmetric pairs.

The involution  $\theta$ on $\gg$ induces an involution on $\gh^*$ which we refer to as
$\Theta$.  Furthermore, $\Theta$ restricts to an involution on $\Delta$.
Set $\pi_{\Theta}=\{\alpha_i\in \pi|\Theta(\alpha_i)=\alpha_i\}$.
Recall [L3, Section 7, (7.5)] that there is a
permutation
$p$ on
$1,\dots, n$ such that
$$-\Theta(\alpha_i)\in \alpha_{p(i)}+Q^+(\pi_{\Theta}).$$  Set $\pi^*=\{\alpha_i\in
\pi\setminus\pi_{\Theta}|i\leq p(i)\}$.  

Given $\alpha\in \gh^*$, set $\tilde\alpha= (\alpha-\Theta(\alpha))/2$.  The subset
$$\Sigma=\{\tilde\alpha|\alpha\in \Delta{\rm \  and \ }\Theta(\alpha)\neq \alpha\}$$
of $\gh^*$ is the
restricted root system associated to the pair $\gg,\gg^{\theta}$.   Note that the set
of (positive) simple 
restricted roots   is just   $\{\tilde\alpha_i|\alpha_i\in \pi^*\}$, while 
$\Sigma^+=\{\tilde\alpha|\alpha\in \Delta^+$ and $\Theta(\alpha)\neq \alpha\}$.
Let $W_{\Theta}$ denote the Weyl group associated to $\Sigma$.

We make the following assumption  throughout this paper:

\medskip
\centerline{$\Sigma$  is  a   reduced  root system.}

\medskip
\noindent
A complete list of the possible irreducible pairs $\gg,\gg^{\theta}$ with 
$\Sigma$ reduced using Araki's classification [A] can be found in the appendix 
of this paper.

Let $U=U_q(\gg)$ be the quantized enveloping algebra generated by $x_1,\dots, x_n,$
$y_1,\dots, y_n,$ $t_1^{\pm 1},\dots, t_n^{\pm 1}$ over
${\cal C}$. Recall that $U$ is a Hopf algebra with coproduct, 
 counit, and antipode.  (See [J1, 3.2.9] or [L3, Section 1,
(1.4)-(1.10)] for relations and Hopf algebra structure.)  Let
$U_+$ denote the  augmentation ideal of $U$.  Given a subalgebra $S$ of $U$ and a 
$U$ module $V$, we write 
$S_+$ for $S\cap U_+$ and set $V^S=\{v\in V|sv=0$ for all $s\in S_+\}$.   

Let $T$ denote the group generated by the
$t_i$ for $1\leq i\leq n$. Set $U^0$ equal to the group algebra ${\cal
C}[T]$.  Recall that there is an isomorphism
$\tau:Q(\pi)\rightarrow T$ such that $\tau(\alpha_i)=t_i$ for each $i$.
Given a $T$ module $N$  and a weight $\gamma\in \gh^*$, the $\gamma$
weight subspace of $N$  is just the set $$\{m\in 
N|\tau(\beta)m=q^{(\beta,\gamma)}m {\rm \ for\ all\ } \tau(\beta)\in T\}.$$
Note that $U$ is a module over itself via the (left) adjoint action denote 
by ad (see for example [J1, 1.3.1]).  If $N$ is a subspace of $U$, then 
$N_{\gamma}$ is the $\gamma$ weight space of $N$ as an $\ad T$ module.

Given a vector subspace  $F$  and a subset $S$ of an algebra   over ${\cal C}$, we set $FS$ 
equal to the ${\cal C}$ span of the set $\{as|a\in F, s\in S\}$. Similarly, write 
$SF$ for the vector space over ${\cal C}$ spanned by the set $\{sa|a\in F, s\in S\}$.
  Now suppose that $S$ is a submonoid of $T$ and $F$ is
both a subalgebra and
$\ad S$ submodule of
$U$. Note that   the subalgebra of $U$ 
generated by $F$ and $S$ is equal to the vector space $FS$.  

Set
$T_{\Theta}=\{\tau(\beta)| \beta\in Q(\pi)$ and $\Theta(\beta)=\beta\}$. Let
 ${\cal M}$ be the  subalgebra of $U_q(\gg)$ generated by $x_i,y_i,t_i^{\pm 1}$ for
$\alpha_i\in \pi_{\Theta}$.   By [L3, Theorem 7.1], we can
lift
$\theta$ to a
${\bf C}$ algebra automorphism $\tilde \theta$ of $U$ which sends $q$ to $q^{-1}.$ 
For each $\alpha_i\in\pi\setminus\pi_{\Theta}$, set 
$$B_i=y_it_i+\tilde\theta(y_i)t_i.\leqno(1.1)$$  Let
$B_{\theta}$ denote the subalgebra of $U$ generated by the 
$B_i,\alpha_i\in\pi\setminus\pi_{\Theta}$, ${\cal M}$, and $T_{\Theta}$.   
By [L3, Theorem 7.2 and the discussion following the proof of Theorem 
7.4], $B_{\theta}$ is a left coideal subalgebra which specializes to 
$U(\gg^{\theta})$ as $q$ goes to $1$. 

Let ${\bf H}$ denote the group of  Hopf algebra automorphisms of $U$ which 
fix elements of $T$.  Note that ${\bf H}$ acts on the set of left coideal
subalgebras of $U$.  Set
${\cal B}$ equal  to the orbit of $B_{\theta}$ under the action of ${\bf H}$. Of course, ${\cal 
B}$ depends on the choice of pair $\gg,\gg^{\theta}$, but this will be 
understood from context. For most 
irreducible pairs $\gg,\gg^{\theta}$, the orbit  ${\cal B}$ equals the set ${\cal 
B}_{\theta}$ defined in [L4, Section 2] and is just the orbit under 
the action of ${\bf H}$ of the quantum analogs of $U(\gg^{\theta})$ 
inside of $U_q(\gg)$.
 There are, however, a  few cases of irreducible pairs for which
the set
${\cal B}_{\theta}$  is strictly larger than ${\cal B}$.
Since we are assuming that
$\Sigma$ is  reduced, this occurs when $U(\gg^{\theta})$ has  
nonstandard  analogs --- also referred  to as analogs of Variation 2
(see [L3, Section 7,  Variation 2]).   We  only  consider the
standard analogs in this paper. 
 A complete list of the   quantum analogs of
$U(\gg^{\theta})$ in $U$ associated to all possible irreducible symmetric pairs can be
found in [L4, Section 7]. (In the
notation of [L4, Sections 2 and 7],
 the
irreducible pair
$\gg,\gg^{\theta}$ admits nonstandard analogs when 
${\cal S}$ is nonempty and consists of one special root $\alpha_i\in \pi^*$. 
When this happens, the orbits of ${\cal B}_{\theta}$ are parametrized by one variable,
$s_i$, and ${\cal B}$ is the
${\bf H}$ orbit in ${\cal B}_{\theta}$ associated to $s_i=0$.)

Given $\lambda\in P^+(\pi)$, let $L(\lambda)$ denote the finite 
dimensional simple $U$ module with highest weight $\lambda$. Write 
$L(\lambda)^*$ for the dual of $L(\lambda)$ given its natural right 
module structure.   Recall the quantum Peter-Weyl theorem ([J1, 9.1.1 and 1.4.13], see also [L3,
(3.1)]): the quantized function
algebra $R_q[G]$  is  isomorphic as a
$U$ bimodule to a direct sum of the 
$L(\lambda)\otimes L(\lambda)^*$ as $\lambda$ varies over the dominant integral weights
$P^+(\pi)$.
 As explained in [L4, Section 4], elements of $L(\lambda)\otimes 
 L(\lambda)^*$, and thus of $R_q[G]$, can be thought of as functions on $U$.   Restriction to the torus $T$
yields an algebra homomorphism, denoted by $\Upsilon$, from $R_q[G]$ into ${\cal C}[P(\pi)]$. 
   
Let $G$ denote the connected, simply connected Lie group $G$ with Lie algebra $\gg$ and
let $K$ be the compact Lie group corresponding to the Lie algebra $\gg^{\theta}$.  In the classical
case, zonal spherical functions are $K$ invariant functions on the symmetric space $G/K$
which are also eigenfunctions for the action of the center of $U(\gg)$.
The   quantum
symmetric space $R_q[G/K]_{B'}$, or more precisely,  quantum analog    of the
ring of regular functions on
$G/K$  associated to
$B'\in {\cal B}$, is the algebra of  left
$B'$ invariants inside the quantized function algebra (see [L3, Section 7 and Theorem 7.8]).
Thus quantum zonal spherical functions at the pair $(B,B')$ in ${\cal B}\times {\cal B}$
are right $B$ invariant elements of $R_q[G/K]_{B'}$ which are eigenfunctions for the action of
the center of
$U$ on
$R_q[G]$. In particular,  quantum zonal spherical functions live inside the space
${}_{B'}{\cal H}_{B}$
 of left $B'$ and  right $B$ invariants of $R_q[G]$. Moreover, the eigenspaces in
$R_q[G]$ for the action of the center of
$U$ are just the subspaces $L(\lambda)\otimes L(\lambda)^*$ for $\lambda$ dominant
integral. 
Hence the  quantum zonal spherical functions at $\lambda$ associated to the pair 
$(B,B')$ are the nonzero elements in the space 
${\ }_{B'}{\cal H}_{B}(\lambda)$ defined by
$${\ }_{B'}{\cal H}_{B}(\lambda)=\{f\in L(\lambda)\otimes L(\lambda)^*|\
B_+'f=fB_+=0\}.$$  By [L4, Theorem 3.4 and Section 4], $_{B'}{\cal H}_{B}(\lambda)$ is one
dimensional if $\lambda\in P^+(2\Sigma)$ and zero otherwise. In 
particular, ${}_{B'}{\cal H}_B$ is a direct sum of the subspaces ${\ 
}_{B'}{\cal H}_{B}(\lambda)$ as $\lambda$ runs over the elements in 
$P^+(2\Sigma)$ ([L4, (4.2)]). By
[L4, Theorem 4.2], the map   
$\Upsilon$ from $R_q[G]$ into ${\cal 
C}[P(\pi)]$ restricts  to an injective algebra  homomorphism from ${}_{B'}{\cal H}_B$
  into ${\cal
C}[P(2\Sigma)]$.  Furthermore, 
[L4, Lemma 4.1] ensures that 
$\Upsilon({}_{B'}{\cal H}_B)$
contains a distinguished basis $\{\varphi^{\lambda}_{B,B'}|\lambda\in 
P^+(2\Sigma)\}$ such that $\varphi^{\lambda}_{B,B'}\in \Upsilon(_{B'}{\cal H}_{B}(\lambda))$ and $$\varphi^{\lambda}_{B,B'}\in 
z^{\lambda}+\sum_{\beta<\lambda}{\cal C}z^{\beta}\leqno{(1.2)}$$ for all 
$\lambda\in P^+(2\Sigma)$.

We recall the notion of zonal spherical families introduced in [L4, 
Section 6].
In particular, a function $\lambda\mapsto \psi_{\lambda}$ from $P^+(2\Sigma)$ 
to ${\cal C}[P(2\Sigma)]$ is called a zonal spherical family 
associated to ${\cal B}$ if there exists $B$ and $B'$ in ${\cal B}$ 
such that $\psi_{\lambda}=\varphi^{\lambda}_{B,B'}$ for all 
$\lambda\in P^+(2\Sigma)$.  With a slight abuse of notation, we generally identify the zonal spherical family above with its
image,  $\{\psi_{\lambda}|\lambda\in 
P^+(2\Sigma)\}$. We have the following
modification of [L4, Theorem 6.5].

\begin{theorem}
Let $\lambda\in P^+(2\Sigma)$.  There exists a unique $W_{\Theta}$ invariant zonal spherical
family $\{\varphi_{\lambda}|\lambda\in P^+(2\Sigma)\}$ associated to ${\cal B}$.
\end{theorem}

\noindent
{\bf Proof:}
In the case when $U(\gg^{\theta})$ does not admit a nonstandard analog, then this is
a consequence of [L4, Theorem 4.2 and Theorem 6.5] as explained in the discussion 
following the proof of the [L4, Theorem 6.5].  (This case corresponds to (i) of 
[L4].)  Now assume that $U(\gg^{\theta})$ does admit nonstandard 
analogs and in particular, ${\cal S}=\{\alpha_i\}$. Let $B\in 
{\cal B}$.  Set 
$$N={\bf Z}\tilde\alpha_i+\sum_{\alpha_j\in \pi^*\setminus{\cal S}}{\bf 
Z}2\tilde\alpha_j.$$
Consider a spherical vector $\xi_{\lambda}\in L(\lambda)$ with respect 
to $B$. By [L4, 
Theorem 3.6 and its proof], $\xi_{\lambda}$ is a sum of weight vectors 
of weight $\lambda-\beta$ where $\beta\in N\cap Q^+(\Sigma)$.   Now $B$ is 
the image under an element in $\bf H$ of the analog in ${\cal 
B}_{\theta}$ associated to $s_i=0$. Hence the proof of [L4, Theorem 
3.6] actually shows that $\xi_{\lambda}$ is a sum of weight vectors of 
weight $\lambda-\beta$ where $\beta$ is in the smaller set 
$Q^+(2\Sigma)$. Thus arguing as in [L4, Lemma 4.1], 
the space ${ }_{B'}{\cal H}_{B}(\lambda)$ is a subspace of 
$z^{\lambda}{\cal C}[Q(2\Sigma)]$ for all pairs $B$ and $B'$ in ${\cal 
B}$.

 Recall 
that ${\cal B}$  is 
the single orbit of $B_{\theta}$ under the action of ${\bf H}$. Thus 
${\cal B}\times {\cal B}$ is a single ${\bf H}\times {\bf H}$ orbit  contained in 
${\cal B}_{\theta}\times {\cal 
B}_{\theta}$.
Note that elements of ${\rm Hom}(N,{\cal C}^{\times})$ act on 
${\cal C}[N]$ as linear transformations where 
$g(z^{\beta})=g(\beta)z^{\beta}$ for all ${\beta}\in N$ and $g\in   
 {\rm Hom}(N, 
{\cal C}^{\times})$. Set ${\cal Z}_{\lambda}$ equal to the set of $W_{\Theta}$
invariant zonal spherical functions at $\lambda$  associated to
${\cal B}$ with top degree term equal to $z^{\lambda}$. By [L4, Theorem
6.5]  there exists a $W_{\Theta}$ invariant zonal spherical family  
 $\{\varphi_{\lambda}|\lambda\in P^+(2\Sigma)\}$  associated to 
 ${\cal B}$
 such that the  $${\cal
Z}_{\lambda}=\{z^{\lambda}g(z^{-\lambda}\varphi_{\lambda})|\ g\in {\rm 
 Hom}(N,{\cal C}^{\times}) {\rm \ and\ }g{\rm\ acts\ trivially\ on\
}\ Q(2\Sigma)\}.$$
 However, by the previous paragraph, $z^{-\lambda}\varphi_{\lambda}\in 
 {\cal C}[Q(2\Sigma)]$.   Hence ${\cal Z}_{\lambda}$ contains exactly one 
 element
 $\varphi_{\lambda}$.  Therefore $\{\varphi_{\lambda}|\lambda\in P^+(2\Sigma)\}$
 is the unique $W_{\Theta}$ invariant 
 zonal spherical family associated to ${\cal B}$.  $\Box$
 
\medskip
Using Theorem 1.1, we write $\{\varphi_{\lambda}|\lambda\in 
P^+(2\Sigma)\}$ for the unique $W_{\Theta}$ invariant zonal spherical family 
associated to ${\cal B}$.   Of course this family depends on the 
choice of $\gg,\gg^{\theta}$, but this will be understood from context.
As an immediate consequence of Theorem 1.1,  if 
$\Upsilon({}_{B'}{\cal H}_B)$ is $W_{\Theta}$ invariant  then 
$\varphi_{B,B'}^{\lambda}=\varphi_{\lambda}$ for all $\lambda\in 
P^+(2\Sigma)$. By [L4, Corollary 5.4], we can choose $B'_{\theta}\in 
{\cal B}$ such that $\Upsilon({}_{B_{\theta}'}{\cal H}_{B_{\theta}})$ is $W_{\Theta}$ invariant.
We drop the subscript $\theta$ and abbreviate $B_{\theta}$ as $B$ 
and $B_{\theta}'$ as $B'$ after Section 2 is completed.

\section{Decompositions and Related Projections}

 In this section, we present   tensor product  decompositions
 of $U$ with respect to a subalgebra $B\in {\cal
B}_{\theta}$ similar to the quantum Iwasawa decomposition of [L1] 
and [L2].
 This, in turn, is used to analyze various projections of
  elements in  $U$ modulo $B_+U$ and $UB_+$.

 Let $T'$ be the subgroup of $T$
generated by $\{t_i|\alpha_i\in \pi^*\}$ of $T$. 
  Note that
$T_{\Theta}\times T'=T$ and so the multiplication map defines a
vector space isomorphism
$$U^0\cong {\cal
C}[T_{\Theta}]\otimes {\cal C}[T'].\leqno{(2.1)}$$

  Let $G^-$ be the subalgebra of $U$ generated by
$y_it_i, 1\leq i\leq n$. Set ${\cal M}^-={\cal M} 
\cap
G^-$ 
 and ${\cal M}^+={\cal M}\cap U^+$.
Let $N^-$ be the subalgebra
  of $G^-$ generated by the $(\ad {\cal
M}^-)$ module $(\ad {\cal M}^-){\cal C}[y_it_i |\alpha_i\notin
\pi_{\Theta}]$. Similarly,  let $N^+$ be the subalgebra of $U^+$
generated by the $(\ad {\cal M}^+)$ module $(\ad {\cal M}^+){\cal
C}[x_i |\alpha_i\notin \pi_{\Theta}]$.   Note that both
$N^-$ and $N^+$ can be written as a direct sum
of weight spaces.

Given a subset $S$ of $U$ and a weight $\beta\in Q(\pi)$, we write
$S_{\beta,r}$ for the restricted weight space of $S$
corresponding to $\beta$.   In particular,
$$S_{\beta,r}=\sum_{\{\beta'|\tilde\beta'=\tilde\beta\}}S_{\beta'}.$$

Recall the standard partial ordering on $\gh^*$: 
For all distinct pairs of elements $\alpha$ and $\beta$ in $\gh^*$, 
$\alpha\leq\beta$ provided that
$\beta-\alpha\in Q^+(\pi)$.  Now suppose $\alpha$ and $\beta$ are in 
$\gh^*$ and $\tilde\alpha<\tilde\beta$.   It follows that 
$\tilde\beta-\tilde\alpha$ is an element of $ Q^+(\pi)\cap (\sum_{\alpha\in \pi^*} {\bf 
C}\tilde\alpha)$. This latter set is contained in $Q^+(\Sigma)$.
 In particular, the partial 
ordering on $\gh^*$ restricts to the standard partial ordering on the 
restricted weights.  

Set $T'_{\geq}$ equal to the multiplicative monoid generated by the 
$t_i^2$ for $\alpha_i\in \pi^*$.  Note that
${\cal C}[T'_{\geq}]$ is just the polynomial ring 
${\cal C} [t_i^{2}|\alpha_i\in \pi^*]$.

\begin{lemma}
For each $B\in {\cal B}$, all $\beta, \gamma\in Q^+(\pi)$,
and $Y\in U^+_{\gamma}G^-_{-\beta}$, we have
$$Y\in
N^+_{\beta+\gamma,r}B+ \sum_{\tilde\beta'<\tilde\beta+\tilde\gamma}
N^+_{\beta',r} T'_{\geq}B\leqno{(2.2)}$$
and
$$Y\in
BN^-_{-\beta-\gamma,r}+\sum_{\tilde\beta'<\tilde\beta+\tilde\gamma} B T'_{\geq}N^-_{-\beta',r}.\leqno{(2.3)}$$

\end{lemma}

\noindent {\bf Proof:}  Let $B\in {\cal B}$. Note that
any Hopf algebra automorphism which fixes $T$ restricts to the identity on $T'_{\geq}$ and  an
automorphism of   $N^-$ and
$N^+$.  Hence, without loss of generality, we may  assume
that $B=B_{\theta}$.

Choose $\beta\in Q^+(\pi)$. By 
construction, $N^+$ is an $\ad {\cal M}^+$ module. Note that $$(\ad
x_i)a=x_ia-t_iat_i^{-1}x_i\leqno{(2.4)} $$ for all $i$ and for all $a\in U$.  Hence 
$${\cal M}^+N_{\beta}^+\subset
N^+_{\beta,r}{\cal M}^+.$$
Now if $\alpha_i\notin\pi_{\Theta}$, it follows that $x_i\in N^+$.
Hence
$$U_{\gamma}^+N_{\beta}^+\subset
N^+_{\beta+\gamma,r}{\cal M}^+$$ for all $\gamma\in Q^+(\pi)$.
 Thus it is
sufficient to prove (2.2) for $Y\in G^-_{-\beta}$.

 It follows from the defining relations of $U$ that
$$G^-_{- \beta}x_i \subset x_i G^-_{-\beta}+G^-_{-\beta+\alpha_i}+
G^-_{-\beta+\alpha_i}t_i^{2}$$ for each $1\leq i\leq n$.  
If $\alpha_i\in \pi_{\Theta}$
then  $ t_i^{2}\in
T_{\Theta} $.  If $\alpha_i\in \pi^*$ then $t_i^{2}\in T'_{\geq}$. Finally, if $\alpha_i\notin
\pi_{\Theta}\cup\pi^*$, then $t_i^{2}=t_{p(i)}^{2}
(t_i^{2}t_{p(i)}^{-2})$ . In particular, $t_i^{2}$ is an element
of $T'_{\geq}T_{\Theta}$ for each $i, 1\leq i\leq n$. Hence
$$ G^-_{-\beta}x_i \subset x_i G^-_{-\beta}+G^-_{-\beta+\alpha_i}
T'_{\geq}T_{\Theta}
\leqno{(2.5)}$$
for all $\alpha_i\in \pi$.

Now consider a weight vector $Y$ of weight $-\beta$ in $G^-$.
Without loss of generality, we may assume that $Y$ is a monomial
in the $y_it_i$, say  $y_{i_1}t_{i_1}\cdots y_{i_m}t_{i_m}$. 
Note that if the restricted weight of $Y$ is zero, or if $m=0$, 
then $Y$ is an element of ${\cal M}^-$, and hence of $B$.  
Thus (2.2) holds in these cases. We
proceed by induction on both $m$ and the restricted weight of $Y$.
In particular, we
 assume that (2.2) holds for all monomials in the
$y_it_{i}$ of length strictly smaller than $m$ as well as for all
elements in $U^+_{\gamma}G^-_{-\lambda}$ with $\gamma\in Q^+(\pi)$
and $\tilde\lambda<\tilde\beta$.

Note that if $\alpha_{i_m}\in \pi_{\Theta}$ then
$y_{i_m}t_{i_m}\in  B_+$.  By the inductive hypothesis  (2.2) holds for
$y_{i_1}t_{i_1}\cdots y_{i_{m-1}}t_{i_{m-1}}$.  It follows that (2.2) holds whenever
$\alpha_{i_m}\in \pi_{\Theta}$. Thus, we may assume that $\alpha_{i_m}\notin
\pi_{\Theta}$.

Recall the definition of $B_i$ (1.1) and note that $B_i\in B_+$
(see for example [L4, (2.1) and (2.2)].  Hence 
$$y_{i_1}t_{i_1}\cdots y_{i_m}t_{i_m}+
y_{i_1}t_{i_1}\cdots
y_{i_{m-1}}t_{i_{m-1}}\tilde\theta(y_{i_m})t_{i_m}
\in y_{i_1}t_{i_1}\cdots
y_{i_{m-1}}t_{i_{m-1}}B_+.$$
Applying the inductive
hypothesis again to
$y_{i_1}t_{i_1}\cdots y_{i_{m-1}}t_{i_{m-1}}$ shows that 
$ Y$ is an element of 
$$ -y_{i_1}t_{i_1}\cdots
y_{i_{m-1}}t_{i_{m-1}}\tilde\theta(y_{i_m})t_{i_m}
+\sum_{\tilde\gamma\leq\tilde\beta-\tilde\alpha_{i_m}
}N^+_{\gamma,r}T'_{\geq}B $$
 By [L4, (2.1) and (2.2)],  $\tilde\theta(y_{i_m})t_{i_m}\in
{\cal M}^+x_{p(i_m)}{\cal M}^+T_{\Theta}$.  Note that $\tilde\beta=0$ for
all $\beta\in Q(\pi_{\Theta})$. Thus by (2.5), $Y$ is contained in
$$
\eqalign{&{\cal M}^+ x_{p(i_m)}{\cal
M}^+T_{\Theta}G^-_{-\beta +\alpha_{i_m},r}+{\cal
M}^+T'_{\geq}T_{\Theta}G^-_{-\beta
+2\alpha_{i_m},r}\cr&+\sum_{\tilde\gamma\leq\tilde\beta-\tilde\alpha_{i_m}
}N^+_{\gamma,r}T'_{\geq}B.\cr}$$

Note that both $\tilde\beta-\tilde\alpha_{i_m}$ and
$\tilde\beta-2\tilde\alpha_{i_m}$ are strictly smaller than
$\tilde\beta$. The result
(2.2) now follows by   induction on $\tilde\beta$.
It follows from [L4, Theorem 3.1] (see also [L4, (3.3) and the proof of Theorem 3.4])
that
$B$ contains elements
$x_i+Y_i$, for $\alpha_i\notin \pi_{\Theta}$,
where $Y_i\in G^-T_{\Theta}$ is a weight vector of weight $\Theta(\alpha_i)$. The proof
of (2.3) is similar to that of (2.2) using the elements $x_i+Y_i$ instead of $B_i$. 
$\Box$
\medskip

Using the above lemma, we obtain four tensor product decompositions
of $U$.

\begin{theorem} For all $B\in {\cal B}$, there
are isomorphisms of vector spaces via the multiplication map
\begin{enumerate}
	\item[(i)]$ N^+ \otimes {\cal C}[T']\otimes B \cong U $
\item[(ii)]$    B\otimes {\cal C}[T']\otimes N^+\cong U$
\item [(iii)]$ N^- \otimes {\cal C}[T']\otimes B \cong U$
\item[(iv)]$    B\otimes {\cal C}[T']\otimes N^-\cong U$
\end{enumerate}
\end{theorem}

\noindent {\bf Proof:} We prove the theorem for $B=B_{\theta}$.   The general result
follows from the fact that $N^+$, $N^-$, ${\cal C}[T']$ and $U$ are all preserved by 
Hopf algebra automorphisms in ${\bf H}$.  

Recall that
$U$ admits a triangular decomposition ([R]) or equivalently, an isomorphism of vector
spaces via the multiplication map:
$$U\cong G^-\otimes U^0\otimes U^+.\leqno{(2.6)}$$ 
By [Ke], (see also [L3, Section 6 and (6.2)]), 
$$U^+\cong {\cal M}^+\otimes N^+\leqno{(2.7)}$$ as vector spaces using the 
multiplication map. 
Combining (2.6) and (2.7) with  (2.1) yields the following
 vector space isomorphism 
$$U\cong  G^-\otimes   {\cal
M}^+\otimes {\cal C}[T_{\Theta}]\otimes{\cal C}[T'] \otimes
N^+\leqno{(2.8)}$$ induced by the multiplication map. 

Set $B_i=y_it_i$ for $\alpha_i\in \pi_{\Theta}$.
Give an $m$-tuple $J=(j_1,\dots, j_m)$, set $y_J=y_1t_1\cdots y_mt_m$ and 
$B_J=B_1\cdots B_m$.   Let ${\cal J}$ be a set of $m$-tuples, where 
$m$ varies, such that the set $\{y_J|J\in {\cal J}\}$ is a basis for 
$G^-$.  By  the proof of [L3, Theorem 7.4], we have
$$B=\bigoplus_{J\in {\cal J}}(B_J{\cal M}^+T_{\Theta}).\leqno{(2.9)}$$ Note
that  when $B_J$ is written as a direct sum of weight vectors, the lowest 
weight term is just $y_J$.  Hence  (2.9)
ensures the lowest weight term of an element of
$B$ is contained in $G^-{\cal M}^+T_{\Theta}$.  It
follows from (2.8) that $Bv\cap Bv'=0$ for any two linearly
independent elements of $T' N^+$.
 This forces the map induced by multiplication from
 $B\otimes {\cal
C}[T']\otimes N^+$ to $U$ to be injective.

Using the  triangular decomposition (2.6) and the relations  of $U$, we 
have that
$U=U^+G^-U^0$. By the previous lemma and (2.1),  $$U = 
U^+G^-U^0\subseteq
 BT'N^+ \subseteq U.$$  
Hence $U= BT'N^+$. Thus the multiplication map
induces a surjection from $B\otimes {\cal C}[T']\otimes
N^+$ onto $U$ which proves (ii).

By [L4, Theorem 3.1], there is a  ${\cal C}$ algebra anti-involution 
$\kappa$ of $U$ 
which fixes elements of $T$,
sends each $x_i$ to $c_iy_it_i$ for some nonzero scalar $c_i$, 
and restricts to a ${\cal C}
$ algebra antiautomorphism of 
$B$. It follows that $\kappa((\ad x_i)b)=c_i(\ad y_i)\kappa(b)$
 for all $b\in U$.  In particular,  $\kappa(N^+)=N^-$.   Thus  assertion 
 (iii) 
 follows from applying $\kappa$ to
(ii).

Let $\iota$ denote the ${\bf C}$ algebra antiautomorphism of $U$ defined by 
$\iota(x_i)=x_i$, $\iota(y_i)=y_i$, $\iota(t_i)=t_i$, and 
$\iota(q)=q^{-1}$. We use ${\it\Delta}$ to denote the coproduct of $U$. It is
straightforward to check using the Hopf  algebra relations of $U$ that $(\iota\otimes
\iota)\circ {\it\Delta}(a)= {\it\Delta}(\iota(a))$ for all $a\in U$. Hence $\iota(B)$ is
also a  left coideal subalgebra of $U$. We recall briefly the notion of specialization 
at $q=1$ (see [L3, Section 1]). Let $\hat U$ denote the ${\bf 
C}[q,q^{-1}]_{(q-1)}$ subalgebra of $U$ generated by $x_i,y_i$, 
$t_i^{\pm 1}$, and $(t_i-1)/(q-1)$ for $1\leq i\leq n$.  Recall 
that $\hat U/(q-1)\hat U$ is isomorphic to $U(\gg)$. 
 Note that $\iota(a)=a+(q-1)\hat U$  and hence the images of $\iota(a)$ and $a$
 are equal 
in $\hat U/(q-1)\hat U$ for all $a\in \hat U$.   It 
follows that both $B$ and  $\iota(B)$ specializes at $q=1$ to the same 
subalgebra, $U(\gg^{\theta})$, of $U(\gg)$.   Now the algebra $\iota(B)$ cannot be an analog of $U(\gg)$
of Variation 1 ([L3, Section 7]) 
since $\Sigma$ is reduced.  A check of the generators of $\iota(B)$ 
shows that $\iota(B)$ cannot be an analog of $U(\gg)$ of  Variation 2
([L3, Section 7]).   Hence, by [L3, Theorem 7.5], $\iota(B)\in 
{\cal B}$. 
Now $\iota$ restricts 
to an antiautomorphism of ${\cal C}[T']$. Furthermore, a 
straightforward computation yields that $\iota((\ad 
x_i)x_j)$ is a scalar multiple of $(\ad x_i)x_j$ for all $i$ and $j$.
It follows that $\iota(N^+)=N^+$.   Similarly $\iota(N^-)=N^-$.  
Hence applying $\iota$ to (ii) yields (i) and applying $\iota$ to 
(iii) yields (iv). 
$\Box$
\medskip

Let ${\cal A}$ denote the subgroup of $T$ generated by
$\tau(2\tilde\alpha)$ as $\alpha $ ranges over $ \pi^*$. Alternatively, we can view 
${\cal A}$ as the image under $\tau$ of the group 
$Q(2\Sigma)$. Let ${\cal A}_{\geq}$
be the semigroup generated by the
$\tau(2\tilde\alpha_i)$ for $\alpha_i\in \pi^*.$  Note that
for $\alpha_i\in \pi^*$,
$$t_i^2=\tau(2\tilde\alpha_i)\tau(\alpha_i+\Theta(\alpha_i))\in
{\cal A}_{\geq}T_{\Theta}.$$  Hence
$${\cal C}[T_{\geq}']\subseteq {\cal C}[{\cal A}_{\geq}]{\cal C}[T_{\Theta}]=
{\cal C}[{\cal A}_{\geq}]+{\cal C}[{\cal A}_{\geq}]{\cal 
C}[T_{\Theta}]_+.\leqno{(2.10)}$$

The following direct sum decompositions of  vector spaces follows
immediately from Theorem 2.2:
$$U= (UB _++U{\cal C}[T']_+)\oplus N^+ \leqno{(2.11)}$$
and
$$U= ( B_+U+{\cal C}[T']_+U)\oplus N^-\leqno{(2.12)} $$  for all $B$ in
${\cal B}$.   Given $B\in {\cal B}$, let
$P_{B}$ be the projection of $ U$ onto $N^+ $ using (2.11)
 and $R_{B}$ be the projection of $U$
onto $N^-$ using (2.12). 
In the next theorem, the projections $P_B$ and $R_B$ are used to
construct particular linear isomorphisms between the restricted
weight spaces of $N^+$ and $N^-$.

\medskip

\begin{theorem}For each $B\in {\cal B}$
and  each $\tilde\beta$, with $\beta\in Q^+(\pi)$, there exist
linear isomorphisms
$$P_{{\tilde\beta},B}:N^-_{-\beta,r}
\rightarrow N^+_{\beta,r}$$ and
$$R_{{\tilde\beta},B}:N^+_{\beta,r}
\rightarrow N^-_{-\beta,r}$$ such that
$$Y-P_{{\tilde\beta},B}(Y)\in N^+_{\beta,r}B_++\sum_{\tilde\gamma<\tilde\beta}N^+_{\gamma,r}{\cal
A}_{\geq}B\leqno(2.13)$$ and
$$X-R_{{\tilde\beta},B}(X)\in B_+N^-_{-\beta,r}+\sum_{\tilde\gamma<\tilde\beta} B{\cal A}_{\geq}N^-_{-\gamma,r}$$
for all $Y\in N^-_{-\beta,r}$ and $X\in
N^+_{\beta,r}.$
\end{theorem}

\noindent {\bf Proof:} Let $B\in {\cal B}$. Fix $\beta\in Q^+(\pi)$. By Lemma 2.1, 
$$P_B(\sum_{\tilde\gamma<\tilde\beta} N^-_{-\gamma,r}) \subseteq
\sum_{\tilde\gamma<\tilde\beta}N^+_{\gamma,r}.$$
Theorem 2.2(iii) ensures that $N^-\cap UB_+=0$. Thus  $P_B$ is
injective.  
Let $\gn^-_{\Theta}$ be the Lie algebra generated by the root 
vectors in $\gg$ corresponding to the set $\{-\gamma|\gamma \in 
{\Delta}^+\setminus Q(\pi_{\Theta})\}$.   By [L3, Section 6], we have 
the equality of formal characters: ${\rm 
ch\ }N^-={\rm ch\ }U(\gn^-_{\Theta})$. Similarly, 
${\rm ch}\ N^+={\rm ch}U(\gn^+_{\Theta})$ where 
$\gn^+_{\Theta}$ is the Lie algebra generated by the root 
vectors in $\gg$ corresponding to the set $\{\gamma|\gamma \in 
{\Delta}^+\setminus Q(\pi_{\Theta})\}$.    Hence both
$\sum_{\tilde\gamma<\tilde\beta}N^-_{-\gamma,r}$ and
$\sum_{\tilde\gamma<\tilde\beta}N^+_{\gamma,r}$ are
finite-dimensional and have the same dimension. It follows that $P_B$
restricted to the former subspace is a bijection. Thus
$$P_B(\sum_{\tilde\gamma<\tilde\beta}N^-_{-\gamma,r}) =
\sum_{\tilde\gamma<\tilde\beta}N^+_{\gamma,r}.\leqno(2.14)$$ 
Note that there is a projection of $N^+$ onto   $N^-_{\beta,r}$ 
with respect to the direct sum decomposition of
$N^+$ into restricted weight spaces.
Set $P_{{\tilde\beta},B}$ equal to the composition of $P_B$ with
this projection. By (2.14), $P_{\tilde\beta,B}$ is an isomorphism of
$N^-_{-\tilde\beta}$ onto $
N^+_{\tilde\beta}.$ Furthermore,   Lemma 2.1
ensures that $Y-P_{\tilde\beta, B}(Y)\in
N^+_{\beta,r}B_++\sum_{\tilde\gamma<\tilde\beta}N^+_{\gamma,r}{\cal 
A}_{\geq}B.$
Thus (2.13) follows. A similar argument constructs
 $R_{\tilde\beta,B}$. $\Box$

\section{Action of the Center on Spherical Functions.}

Set $Q_{\Sigma}=Q(\Sigma)$.  Note that $Q_{\Sigma}$ is a subset of 
$P(\Sigma)$ and hence ${\cal C}[Q_{\Sigma}]$ is a 
subring of ${\cal C}[P(\Sigma)]$.   Thus ${\cal C}[P({\Sigma})]$ is a right 
${\cal C}[Q_{\Sigma}]$ module where elements of ${\cal C}[Q_{\Sigma}]$ act as right
multiplication.  It follows that   we may embed
${\cal C}[Q_{\Sigma}]$ into the (right) 
endomorphism ring  ${\rm End}_r\ {\cal C}[P(\Sigma)]$ of ${\cal C}[P(\Sigma)]$.   

Note that ${\cal C}[P(\Sigma)]$ is also a right
${\cal C}[{\cal A}]$ module where $$z^{\lambda}*\tau(\mu)=q^{(\lambda,\mu)}z^{\lambda}$$   
 for all $z^{\lambda}\in {\cal C}[P(\Sigma)]$ and $\tau(\mu)\in {\cal A}$.  Since the
Cartan inner product restricts to a nondegenerate bilinear form on $P(\Sigma)\times
Q_{\Sigma}$, it follows that the action of ${\cal C}[{\cal A}]$ on ${\cal C}[P(\Sigma)]$
is faithful.   Hence ${\cal C}[{\cal A}]$ also embeds in ${\rm End}_r\ {\cal C}[P(\Sigma)]$.   Let ${\cal
C}[Q_{\Sigma}]{\cal A}$ denote the subring of ${\rm End}_r\ {\cal C}[P(\Sigma)]$ generated by ${\cal C}[Q_{\Sigma}]$ and
${\cal A}$.  Note that 
$$z^{\lambda}\tau(\mu)=q^{(\lambda,\mu)}\tau(\mu)z^{\lambda}\leqno{(3.1)}$$
for all $\lambda\in Q_{\Sigma}$ and $\tau(\mu)\in {\cal A}$. Furthermore (3.1) implies
that the nonzero elements of ${\cal C}[Q_{\Sigma}]$ form an Ore set in ${\cal
C}[Q_{\Sigma}]{\cal A}$.   Write ${\cal C}(Q_{\Sigma})$ for the quotient ring of ${\cal
C}[Q_{\Sigma}]$ and set ${\cal C}(Q_{\Sigma}){\cal A}$ equal to the localization of ${\cal
C}[Q_{\Sigma}]{\cal A}$ at the Ore set ${\cal C}[Q_{\Sigma}]\setminus \{0\}.$
 In
this section, we obtain a homomorphism of the center $Z(U)$ of 
 $U$ into ${\cal C}(Q_{\Sigma})[{\cal A}]$ which corresponds to the 
 action of $Z(U)$
on the zonal spherical functions.

We give ${\cal C}[{\cal A}]$ the structure of a left  ${\cal
C}[Q_{\Sigma}]{\cal A}$ module as follows.   Elements of ${\cal
A}$ act by left multiplication while
 $$z^{\lambda}\cdot
\tau(\mu)=q^{(\lambda,\mu)}\tau(\mu)$$ for each $\lambda\in
Q_{\Sigma}$ and $\tau(\mu)\in {\cal A}$.  In particular, we may also view ${\cal
C}[Q_{\Sigma}]{\cal A}$ as the subring of the (left) endomorphism ring ${\rm End}_l\ {\cal C}[{\cal A}]$ of ${\cal
C}[{\cal A}]$ generated by ${\cal C}[Q_{\Sigma}]$ and ${\cal A}$.

The action of ${\cal C}[Q_{\Sigma}]{\cal A}$ on ${\cal C}[{\cal A}]$ can
 be extended to an action of elements in  ${\cal C}(Q_{\Sigma}){\cal 
A}$ on certain elements of ${\cal A}$ as long as we avoid   denominator 
problems. In particular, consider $f\in {\cal
C}[Q_{\Sigma}]{\cal A}$ and $g\in {\cal
C}[Q_{\Sigma}]$.  Suppose that $\tau(\mu)\in {\cal A}$ such that 
$g\cdot \tau(\mu)\neq 0$.  Note that $(g\cdot 
\tau(\mu))\tau(\mu)^{-1}$ is just an element of ${\cal C}$.  
We denote $(f\cdot
\tau(\mu))(g\cdot 
\tau(\mu)\tau(\mu)^{-1})^{-1}$ by  $(fg^{-1})\cdot \tau(\mu)$.

Note that   
 the algebra ${\cal C}[P(2{\Sigma})]$ can be identified with  a
subspace of  the dual of ${\cal C}[{\cal A}]$ where 
$$z^{\lambda}(\tau(\mu))=q^{(\lambda,\mu)}$$ for all 
$z^{\lambda}\in {\cal C}[P(2{\Sigma})]$ and $\tau(\mu)\in {\cal 
C}[{\cal A}]$.
Moreover, the above two actions ${\cal C}[Q_{\Sigma}]{\cal A}$ are compatible with the
pairing between
${\cal C}[P(2{\Sigma})]$ and ${\cal C}[{\cal A}]$.   In particular, given $a'\in {\cal
C}[P(2{\Sigma})],
a\in {\cal C}[{\cal A}]$, and $b\in {\cal C}(Q_{\Sigma}){\cal A}$, we obtain
$$a'(b\cdot a)=(a'*b)(a).$$

Recall that $B'_{\theta}\in {\cal B}$ has been chosen so that
the image of ${ }_{B'_{\theta}}{\cal H}_{B_{\theta}}$ in ${\cal 
C}[P(2{\Sigma})]$ is $W_{\Theta}$ invariant (see the end of Section 1). For the remainder of 
the paper, we will drop the $\theta$ subscript, setting
$B=B_{\theta}$ and $B'=B'_{\theta}$.

Let ${\cal C}(Q_{\Sigma}){\cal A}_{\geq}$ denote the subalgebra of 
${\cal C}(Q_{\Sigma}){\cal A}$ generated by ${\cal C}(Q_{\Sigma})$
and ${\cal A}_{\geq}.$ Set $U_{\geq}=U^+G^-{\cal A}_{\geq}$.

\begin{lemma}  For each $\tau(\gamma)\in {\cal A}$ and $X\in U^+G^-\tau(\gamma)$, there exists
$p_X\in {\cal C}(Q_{\Sigma}){\cal A}_{\geq}\tau(\gamma)$ such that
$$X\tau(\lambda)-(p_X\cdot\tau(\lambda)) \in
 B_+(U_{\geq}\tau(\gamma+\lambda))+(U_{\geq}\tau(\gamma+\lambda)){B_+'}\leqno{(3.2)}$$
for all 
$\tau(\lambda)\in {\cal A}$
 such that $p_X\cdot\tau(\lambda)$ is defined.
 \end{lemma}

\noindent{\bf Proof:}
Given
$\beta\in Q^+(\pi)$, set
$P_{\tilde\beta,B'}=P_{\tilde\beta}$ and $R_{\tilde\beta,B}=R_{\tilde\beta}$.
By Lemma 2.1 and (2.10), we may reduce to the case when
$X\in N^+{\cal A}_{\geq}\tau(\gamma)$. Note that if $X$ and 
$X'$
both satisfy (3.2), then so does $X+X'$.   Hence we may assume
that there exists a $\beta\in Q^+(\pi)$ such that 
$X\in N^+_{\beta,r}\tau(\gamma')$ for some 
 $\tau(\gamma')\in {\cal A}_{\geq}\tau(\gamma)$.

Given $\alpha\in Q^+(\pi)$, set ${\rm 
ht}_r(\alpha)=\sum_{\alpha_i\in \pi^*}m_{\alpha}^i$ where $\tilde\alpha=\sum_{\alpha_i\in \pi^*}
m_{\alpha}^i\tilde\alpha_i$. We prove the 
lemma by induction on ${\rm ht}_r(\beta)$.  In particular, assume first
that
${\rm ht}_r(\beta)=0$. Hence $\tilde\beta=0$.  Since $N^+T\cap {\cal 
M}T=T$, it further follows that $X\in 
T$ and (3.2) holds with $p_X=1$.  Now assume that ${\rm
ht}_r(\tilde\beta)>0$ and
(3.2) is true for  all elements in $N^+_{\gamma,r}T$ with ${\rm
ht}_r(\gamma)<{\rm  ht}_r(\beta)$.

By Theorem 2.3, $P_{\tilde\beta}\circ R_{\tilde\beta}$ is an
isomorphism of $N^+_{\beta,r} $ onto itself.   Let $X_i,
1\leq i\leq m,$ be a basis for $N^+_{\beta,r} $ considered as a
vector space over ${\cal C}$ so that $P_{\tilde\beta}\circ
R_{\tilde\beta}$ is an upper triangular $m\times m$ matrix.   The
fact that $P_{\tilde\beta}\circ R_{\tilde\beta}$ is an isomorphism
ensures that the diagonal entries of this matrix, say $c_{ii}$, are nonzero.

 By
Theorem 2.3, $X_i\tau(\gamma')$ is an element of
$$R_{\tilde\beta}(X_i)\tau(\gamma'
) +\sum_{\tilde\xi<\tilde\beta}N^-_{-\xi,r}{\cal 
A}_{\geq}\tau(\gamma)
+B_+U_{\geq}\tau(\gamma).$$ Note that $\tilde\xi<\tilde\beta$ 
implies that ${\rm ht}_r(\xi')< {\rm ht}_r(\beta)$ for all $\xi'\in 
Q^+(\pi)$ satisfying $\tilde\xi'=\tilde\xi$. Thus by the inductive 
hypothesis, there
exists
$p_1\in {\cal C}(Q_{\Sigma}){\cal A}_{\geq}\tau(\gamma)$ such that
$X_i\tau(\gamma')\tau(\lambda)$ is an element of $$ R_{\tilde\beta}(X_i)\tau(\gamma' )\tau(\lambda)
+p_1\cdot
\tau(\lambda)+B_+U_{\geq}\tau(\gamma+\lambda)+U_{\geq}\tau(\gamma+\lambda)B_+'\leqno{(3.3)}$$
for all $\lambda$ such that $p_1\cdot \tau(\lambda)$ is defined.
Since
$R_{\tilde\beta}(X_i)\in N^-_{-\beta,r}$, we further have that
$$\tau(\gamma')R_{\tilde\beta}(X_i)
\in\ \tau(\gamma')P_{\tilde\beta}(R_{\tilde\beta}(X_i))
+\sum_{\tilde\xi<\tilde\beta}{\cal
A}_{\geq}N^+_{\xi,r}\tau(\gamma)+U_{\geq}\tau(\gamma)B'_+.$$ 
Applying induction to elements of $N^+_{\xi,r}$,
we can find $p_2\in {\cal C}(Q_{\Sigma}){\cal
A}_{\geq}\tau(\gamma)$ such that
$$\eqalign{ \tau(\lambda)\tau(\gamma')&R_{\tilde\beta}(X_i)-
\tau(\lambda)\tau(\gamma')P_{\tilde\beta}(R_{\tilde\beta}(X_i))-p_2\cdot
\tau(\lambda)\cr&\in
B_+U_{\geq}\tau(\gamma+\lambda)+U_{\geq}\tau(\gamma+\lambda)B_+'\cr}\leqno{(3.4)}$$
for all $\lambda$ such that $p_2\cdot
\tau(\lambda)$ is defined.
Set $p_2'=q^{(\gamma',\tilde\beta)}z^{\tilde\beta}p_2$.  Combining (3.3) and (3.4)  yields
$$\eqalign{ X_i\tau(\gamma')\tau(\lambda)-&
q^{(2\gamma'+2\lambda,\tilde\beta)}
P_{\tilde\beta}\circ
R_{\tilde\beta}(X_i)\tau(\lambda)\tau(\gamma')- (p_1+p_2')\cdot \tau(\lambda)\cr&\in
B_+U_{\geq}\tau(\gamma+\lambda)+U_{\geq}\tau(\gamma+\lambda)B_+'\cr}$$
for all $\lambda$ such that $(p_1+p_2')\cdot
\tau(\lambda)$ is defined.

By the choice of the $\{X_i\}$, it follows that
$$P_{\tilde\beta}(R_{\tilde\beta}(X_i))\in c_{ii}(X_i)
+\sum_{j\leq i}{\cal C}X_j.$$ Hence
$$\eqalign{&X_i\tau(\gamma')(1-c_{ii}
z^{2\tilde\beta}q^{(2\gamma',\tilde\beta)})\cdot \tau(\lambda)
-(p_1+p_2')\cdot\tau(\lambda)\cr
&\in \sum_{j\leq i}{\cal
C}X_j\tau(\gamma'+\lambda)
(z^{2\tilde\beta})\cdot\tau(\lambda)  +
B_+U_{\geq}\tau(\gamma+\lambda)+U_{\geq}\tau(\gamma+\lambda)B_+'.\cr}$$ By induction on
$i$, there exist
$s_j$ in ${\cal C}$ and $p_3\in {\cal C}(Q_{\Sigma}){\cal
A}_{\geq}\tau(\gamma)$ such that
$$ X_i\tau(\gamma')\prod_{1\leq j\leq i}(1-s_jz^{2\tilde\beta})\cdot
\tau(\lambda) -  p_3\cdot\tau(\lambda) \in 
B_+U_{\geq}\tau(\gamma+\lambda)+U_{\geq}\tau(\gamma+\lambda)B_+'$$
for all $\lambda$ such that $p_3\cdot \tau(\lambda)$ is defined. The lemma follows by
setting
$p_{X_i}=p_3\prod_{1\leq j\leq i}(1-s_jz^{2\tilde\beta})^{-1}$ and 
noting that $p_X$ is a linear combinaton of the $p_{X_i}$. $\Box$

\medskip
Let $\check U$ denote the simply connected quantized enveloping algebra 
([J1, Section   3.2.10].)  Recall that $\check U$ is generated by $U$ and the torus
$\check T=\{\tau(\lambda)|\lambda\in P(\pi)\}$ 
corresponding to the weight lattice.  Set $$\check {\cal A}=\{\tau(\tilde\mu)|\mu\in
P(\pi)\}.$$  The ring ${\cal C}(Q_{\Sigma})\check {\cal A}$ is defined
in an analogous way to ${\cal C}(Q_{\Sigma}){\cal
A}$ using the fact    that the right action of ${\cal A}$ on ${\cal C}[P(\Sigma)]$ extends to 
$\check {\cal A}$. 

Let $\check U^0$ denote the group algebra ${\cal C}[\check T]$ and 
set  $\check T_{\Theta}=\{\tau(\mu)|\tau(\mu)\in \check T$
and 
$\Theta(\mu)=\mu\}$.  The definition of $\check{\cal A}$ and $\check T_{\Theta}$ yields the following inclusion:
$$\check U^0\subset {\cal C}[\check{\cal A}]\oplus \check U^0{\cal C}[\check T_{\Theta}]_+.\leqno{(3.5)}$$
Hence $$\tau(\gamma)=\tau(\tilde\gamma)+ \tau(\tilde\gamma)(\tau({{1}\over{2}}(\gamma+\Theta(\gamma))-1)\in
\tau(\tilde\gamma)+U^0{\cal C}[\check T_{\Theta}]_+$$ for all $\tau(\gamma)\in \check T$. 
 It follows that Lemma 3.1
extends to elements 
$X\in U^+G^-\tau(\gamma)$  for any $\tau(\gamma)\in \check T$, where 
${\cal C}(Q_{\Sigma}){\cal A}_{\geq}\tau(\gamma)$ is
replaced by 
${\cal C}(Q_{\Sigma}){\cal A}_{\geq}\tau(\tilde\gamma)$, $B_+$ is replaced by
$(B
\check T_{\Theta})_+$, and
$B'_+$ is replaced by $(B'\check T_{\Theta})_+$.

 Set $T_{\geq}$ equal to the submonoid of $T$ generated by $t_i^2$, 
for $i=1,\dots, n$. 
Consider  $X$ and $p_X$ defined  as in the previous lemma.   Let
$g_{\lambda}$  be the zonal spherical function in $_{B'}{\cal
H}_{B}(\lambda)$
 with image $\varphi_{\lambda}\in {\cal C}[P(2{\Sigma})]$ where 
 $\varphi_{\lambda}$ is chosen as in the end of Section 1.  Assume further that
 $p_X\cdot \tau(\beta) $ is defined. Then
$$g_{\lambda}(X\tau(\beta))=g_{\lambda}(p_X\cdot\tau(\beta))=\varphi_{\lambda}
(p_X\cdot\tau(\beta)).$$ Hence
$(\varphi_{\lambda}*p_X)(\tau(\beta))=g_{\lambda}(X\tau(\beta))$
for all $\beta$ such that $p_X\cdot \tau(\beta)$ is defined.   We have established the following result.

\begin{theorem}  There is a linear map ${\cal X} :\check U\rightarrow {\cal
C}(Q_{\Sigma})\check {\cal A}$ such that
$$g_{\lambda}(u\tau(\beta))=(\varphi_{\lambda}*
{\cal X}(u))(\tau(\beta))$$ for all $u\in \check U$, $\lambda\in
P^+(2\Sigma)$ and $\tau(\beta)\in {\cal A}$ such that ${\cal
X}(u)\cdot \tau(\beta)$ is defined. Furthermore, if $u\in
U^+G^-T_{\geq}\tau(\gamma)$, then ${\cal X}(u)\in {\cal
C}(Q_{\Sigma}){\cal A}_{\geq}\tau(\tilde\gamma)$.
\end{theorem}

Let $Z(\check U)$ denote the center of $\check U$. The restriction of ${\cal X}$
to  
$Z(\check U)$  is particularly nice.  
Recall that $\check U$ admits a   direct sum decomposition
$$\check U=\check U^0\oplus(G^-_+\check U+\check UU^+_+).$$ 
Let ${\cal P}$ denote the
quantum Harish-Chandra projection  of $\check U$ onto $\check U^0$ using
this decomposition. A central
element $c$ acts on elements of $L(\lambda)$, and hence on 
the zonal spherical function $g_{\lambda}$, as multiplication by  
the scalar $z^{\lambda}({\cal P} (c))$.  In particular
$$g_{\lambda}(c\tau(\gamma))=z^{\lambda}({\cal
P}(c))(\varphi_{\lambda}(\tau(\gamma)))$$ for all $c\in Z(\check U)$ and
$\tau(\gamma)\in {\cal A}$.   It follows that $$(\varphi_{\lambda}*{\cal
X}(c)-z^{\lambda}({\cal P}(c))\varphi_{\lambda})\cdot 
\tau(\gamma)=0\leqno{(3.6)}$$
for all $\tau(\gamma)$ such that ${\cal X}(c)\cdot\tau(\gamma)$ is defined.
The next result shows that (3.6) holds for all 
$\tau(\gamma)$.

\begin{corollary}    The restriction of  ${\cal X}$ to $Z(\check U)$ 
	is an algebra homomorphism from $Z(\check U)$ to $ {\cal
C}(Q_{\Sigma})\check{\cal A}$ such that
$$\varphi_{\lambda}*{\cal
X}(c)=z^{\lambda}({\cal P}(c))\varphi_{\lambda}\leqno{(3.7)}$$ and thus
$$g_{\lambda}(c\tau(\beta))=(\varphi_{\lambda}*
{\cal X}(c))(\tau(\beta))\leqno{(3.8)}$$ for all $c\in Z(\check U)$, $\lambda\in
P^+(2\Sigma)$, and $\tau(\beta)\in {\cal A}$. Furthermore, if $z\in
Z(\check U)$  and $z\in U^+G^-T_{\geq}\tau(\gamma)$, 
then ${\cal X}(z)\in
{\cal C}(Q_{\Sigma}){\cal A}_{\geq}\tau(\tilde\gamma)$.

\end{corollary}

\noindent
{\bf Proof:} Let $c\in Z(\check U)$. Note that we can find a nonzero
element 
$p=\sum p_{\beta}z^{\beta}$ in ${\cal C}[Q_{\Sigma}]$
such that
$p{\cal X}(c)$ is in the subring of ${\cal
C}(Q_{\Sigma})\check {\cal A}$ generated by
${\cal C}[Q_{\Sigma}]$ and $\check {\cal A}$.  It follows that there 
exists $f=\sum f_{\beta}z^{\beta}$ in ${\cal 
C}[Q_{\Sigma}]$ such that $p^{-1}f=\varphi_{\lambda}*{\cal
X}(c)-z^{\lambda}({\cal P}(c))\varphi_{\lambda}.$  Hence 
$\sum f_{\beta}q^{(\beta,\gamma)}=0$ for all $\gamma$ such that
$\sum_{\beta} p_{\beta}q^{(\beta,\gamma)} \neq 0$.   

Assume that $f\neq 0$.  Choose
$\gamma$ such that 
$(\gamma,\beta)\neq 0$ for  at least one $\beta$ with $f_{\beta}\neq 0$
and $(\gamma,\beta')\neq 0$ for at least one $\beta'$ with $p_{\beta'}\neq
0$.    A standard 
Vandermonde determinant argument shows that there exists $N\geq 0$ 
such that 
$\sum_{\beta}f_{\beta}q^{(\beta,m\gamma)}\neq 0$ for all $m\geq N$.
By the previous paragraph,
$\sum p_{\beta}q^{(\beta,m\gamma)}=0$ for all  $m\geq N$.   Another
application of the Vandermonde determinant  argument yields $p=0$, a
contradiction.   Hence $f=0$.  This proves   (3.7) and (3.8) immediately
follows.  The last assertion of the  corollary is a direct consequence of
Theorem 3.2. $\Box$
 
\medskip

The restricted Weyl group
$W_{\Theta}$ acts 
on ${\cal C}(Q_{\Sigma})
\check{\cal A}$   by
$$s\cdot \tau(\mu)=\tau(s\mu)$$
and
$$s\cdot z^{\mu}=z^{s\mu}$$
for any $s\in W_{\Theta}$ and $\tau(\mu)\in \check{\cal A}$.
In the classical case, elements of the center of the classical 
enveloping algebra of $U(\gg)$  can be realized
as $W_{\Theta}$ invariant elements of ${\cal C}(Q_{\Sigma})\check{\cal A}$
with respect to their action on spherical functions.
The next theorem is a    quantum version of this result.

The group algebra ${\cal C}[Q_{\Sigma}]$ is just the Laurent polynomial ring 
corresponding to the polynomial ring ${\cal C} [z^{-\tilde\alpha_i}|\alpha_i\in \pi^*]$.
 Let
$ {\cal C}((Q_{\Sigma}))$ denote the formal Laurent series ring ${\cal C}((z^{-\tilde\alpha_i}|\alpha_i\in \pi^*)).$
In particular, the ring
$ {\cal C}((Q_{\Sigma}))$ consists of finite linear combinations of possibly infinite sums of the form
$\sum_{\gamma\leq\beta} a_{\gamma}z^{\gamma}$  where $\gamma$ and $\beta$ are elements
of 
$
Q_{\Sigma}$ and each
$a_{\gamma}\in {\cal C}$.  Note that ${\check {\cal A}}$ embeds inside ${\rm End}_r {\cal C}((Q_{\Sigma}))$
where $$\sum_{\gamma\leq \beta}a_{\gamma}z^{\gamma}*\tau(\nu)=\sum_{\gamma\leq
\beta}a_{\gamma}q^{(\gamma,\nu)}z^{\gamma}\leqno{(3.9)}$$ for all $\sum_{\gamma\leq 
\beta}a_{\gamma}z^{\gamma}\in {\cal
C}((Q_{\Sigma}))$ and $\tau(\nu)\in \check {\cal A}$. Let ${\cal C}((Q_{\Sigma}))\check {\cal A}$ denote the subring of
${\rm End}_r {\cal C}((Q_{\Sigma}))$ generated by $  {\cal C}((Q_{\Sigma}))$  and $\check {\cal A}$. 
 The quotient ring 
$ {\cal C}(Q_{\Sigma})$ embeds inside of $ {\cal C}((Q_{\Sigma}))$ in a standard way. It
follows that that ${\cal C}(Q_{\Sigma})\check{\cal A}$ is a subring
of ${\cal C}((Q_{\Sigma}))\check{\cal A}$.  

The relations in (3.9) ensures 
that the multiplication map yields  vector 
space isomorphisms  $${\cal C}((Q_{\Sigma}))\check{\cal A}\cong {\cal 
C}((Q_{\Sigma}))\otimes{\cal C}[\check{\cal A}]\cong{\cal 
C}[\check{\cal A}]\otimes {\cal 
C}((Q_{\Sigma})).\leqno{(3.10)}$$  In particular, elements of ${\cal C}((Q_{\Sigma}))\check{\cal A}$
are finite sums of the form $\sum_{i}a_ib_i$ where $a_i\in {\cal C}[\check{\cal 
A}]$ and $b_i\in {\cal C}((Q_{\Sigma}))$.   Grouping together the 
coefficients of each $z^{\beta}$, we can write any element in ${\cal 
C}((Q_{\Sigma}))\check{\cal A}$ as a finite linear combination of 
possibly infinite sums  of the form 
$\sum_{\gamma<\beta}a_{\gamma}z^{\gamma}$ where $\gamma$ and $\beta$ 
are in $Q_{\Sigma}$ and each $a_{\gamma}\in {\cal C}[\check{\cal 
A}]$.  However, the reader should be aware that not all such sums are elements of ${\cal 
C}((Q_{\Sigma}))\check{\cal A}$.

Let $\omega_i'$ be the fundamental weight in $P^+(\Sigma)$ corresponding to the restricted 
root $\tilde\alpha_i$.  Since
$Q(\Sigma)$ is a subset of
$P(\Sigma)$, it follows that ${\cal C}((Q_{\Sigma}))$ is a subring of the Laurent power series ring
${\cal
C}((z^{-\omega_i'}|\alpha_i\in \pi^*))$.  Thus we may view
${\cal C}((Q_{\Sigma}))\check {\cal A}$ as a subring of ${\rm End}_r {\cal
C}((z^{-\omega_i'}|\alpha_i\in \pi^*)).$  This interpretation will
be needed in the proof below where elements of ${\cal
C}((Q_{\Sigma}))\check {\cal A}$ act on zonal spherical functions.

\begin{theorem}   The image of $Z(\check U)$ under ${\cal X}$ is contained in
$({\cal C}(Q_{\Sigma})\check{\cal A})^{W_{\Theta}}$.
\end{theorem}

\noindent{Proof:}  Let $a\in Z(\check U)$ and $s=s_{\alpha}$ be the
reflection in 
$W_{\Theta}$ corresponding to a simple
root
$\alpha\in \Sigma$. Assume that ${\cal X}(a)\neq s\cdot {\cal X}(a)$.
We can think of ${\cal X}(a)-s\cdot {\cal X}(a)$ as an
element of ${\cal C}((Q_{\Sigma}))\check{\cal A}$.  It follows from 
(3.10) that
 there exists a finite set $\{\gamma_1,\dots,\gamma_r\}$ of 
 noncomparable elements in $Q_{\Sigma}$ and  
 elements $a_{\beta_i}\in{\cal
C}[\check{\cal A}]$ for $\beta_i\leq \gamma_i$ and $1\leq i\leq r$ such that
 $${\cal X}(a)-s\cdot
{\cal X}(a)=\sum_i\sum_{\beta_i\leq\gamma_i}a_{\beta_i}z^{\beta_i}.$$  
We may further assume that $a_{\gamma_i}\neq 0$ for each $1\leq i\leq r$.

Recall  that the zonal spherical function $\varphi_{\lambda}$ can be written
as a sum of the form
$$z^{\lambda}+\sum_{\beta<\lambda}
y_{\beta}z^{\beta}$$ where each $\beta\in Q_{\Sigma}$ and each
$y_{\beta}\in {\cal C}$ (see (1.2)).  Hence
$$0=\varphi_{\lambda}*({\cal X}(a)-s\cdot {\cal
X}(a))\in\sum_i(z^{\lambda}(a_{\gamma_i}))z^{\lambda+\gamma_i}+\sum_{i}\sum_{\beta<\lambda+
\gamma_i}{\cal
C}z^{\beta}.$$
 Hence  $z^{\lambda_i}(a_{\gamma_i})=0$  for all $1\leq i\leq r$ and all $\lambda\in P^+(2\Sigma)$.
 This forces each $a_{\gamma_i}=0$, a contradiction. $\Box$
 
 \medskip

We wish to extend Corollary 3.3 and Theorem 3.4 to $B$ invariant 
elements of $\check U$.  In particular, 
let $\check U^B$ denote the subalgebra of $\check U$ consisting of  $\adr B$ invariant
elements.      Since $B$ is not a Hopf subalgebra of $U$, it is not
obvious that the centralizer $C_{\check U}(B)=\{c\in \check U|bc=cb$  for all $b\in B\}$ of 
$B$ in ${\check U}$ is 
equal to the set ${\check U}^B$.  
Nevertheless, the next lemma shows that this is indeed true.

\begin{lemma}  $\check U^B=C_{\check U}(B)$. In
particular, $Z(\check U)$ is a subset of $\check U^B$.
\end{lemma}

\noindent
{\bf Proof:} Note that the second assertion is an immediate consequence
of the first. It is straightforward to check that  $C_{\check U}(B)\subseteq \check 
U^B$.  (The argument follows as in [J1, Lemma 1.3.3] using the right 
adjoint action instead of the left adjoint action.)

Recall that 
$\sum \sigma(a_{(1)})a_{(2)}=\epsilon(a)$ for all $a\in U$ where 
 $\epsilon$ is the counit for $U$ and the coproduct is given in Sweedler notation, 
${\it\Delta}(a)=\sum a_{(1)}\otimes a_{(2)}$. Recall further that $B$ 
is a left coideal and so ${\it\Delta}(B)\subset U\otimes B$. 
Suppose that  $c\in \check U^B$.  Then 
 $$\eqalign{ac&=\sum 
  a_{(1)}\epsilon(a_{(2)})c=\sum
 a_{(1)}\sigma(a_{(2)})ca_{(3)}\cr&=\sum\epsilon(a_{(1)})ca_{(2)}=c\sum\epsilon
 (a_{(1)})a_{(2)}=ca\cr}$$ for all $a\in B$.  Hence $c\in C_{\check 
 U}(B)$ and $\check U^B\subseteq C_{\check 
 U}(B)$. $\Box$

\medskip
Theorem 2.2(ii) and (3.5) imply the following inclusion 
$$\check U\subseteq ((B\check T_{\Theta})_++N^+_+\check {\cal A})\oplus {\cal C}[\check A].\leqno{(3.11)}$$   
Let ${\cal P}_{\cal A}$ denote the 
projection of $\check U$ onto ${\cal C}[\check{\cal A}]$ using this 
decomposition.

Let $\lambda\in P^+(2\Sigma)$.   Recall [L4, Theorem 3.2] that
$(L(\lambda)^*)^B$ is  one dimensional. Choose a nonzero generating
vector $v_{\lambda}^*$ of weight $\lambda$ for $L(\lambda)^*$. By [L4, Lemma 3.3], we
can choose a nonzero 
 vector $\xi^*_{\lambda}$
  in $(L(\lambda)^*)^B$ such that 
  $\xi^*_{\lambda}=v^*_{\lambda}+v^*_{\lambda}N^+_+$. Suppose that 
$c\in \check U^B$. By the previous lemma, 
$\xi^*_{\lambda}cu= \xi_{\lambda}^*uc=0$ for all $u\in B_+$.  Hence 
 $\xi^*_{\lambda}c$ is a scalar multiple of $\xi^*_{\lambda}$. It 
 follows from (3.11) and the definition of ${\cal P}_{\cal A}$ that 
 $$\xi^*_{\lambda}c\in \xi^*_{\lambda}({\cal P}_{\cal A}(c) 
 +N^+_+)\subseteq z^{\lambda}({\cal P}_{\cal 
 A}(c))v^*_{\lambda}+v^*_{\lambda}N^+_+.$$  Hence $\xi^*_{\lambda}c=   
 z^{\lambda}({\cal P}_{\cal 
 A}(c))\xi^*_{\lambda}$.
 
  Consider the special case when $c$ is an element in $Z(\check U)$.
Now $\xi^*_{\lambda}\in v^*_{\lambda}+v^*_{\lambda}B_+$ (see the proof
of [L3, Theorem 7.7]).    In particular, we can find a $b$ such that
$b-1\in B_+$ and 
$\xi^*_{\lambda}=v^*_{\lambda}b$.   It follows that 
$v^*_{\lambda}bc=v^*_{\lambda}cb$.   Therefore, 
$\xi^*_{\lambda}c=z^{\lambda}({\cal P}(c))\xi^*_{\lambda}$.  By the 
previous paragraph, we see that  $z^{\lambda}({\cal P}_{\cal A}(c))=z^{\lambda}({\cal P}(c))$ for all $c\in 
Z(\check U)$.  In particular, ${\cal P}_{\cal A}(c)$ agrees with the 
image of ${\cal P}(c)$ under  projection onto ${\cal C}[\check{\cal A}]$ 
using (3.5). Thus
arguing as in Corollary 3.3 and Theorem 3.4, we have the 
following generalization of Theorem 3.4.

\begin{theorem} The restriction of  ${\cal X}$ to $\check U^B$ 
	is an algebra homomorphism from $\check U^B$ to $( {\cal
C}(Q_{\Sigma})\check{\cal A})^{W_{\Theta}}$ such that
$$\varphi_{\lambda}*{\cal 
X}(c)=z^{\lambda}({\cal P}_{\cal A}(c))\varphi_{\lambda}$$ and 
$$g_{\lambda}(c\tau(\beta))=(\varphi_{\lambda}*
{\cal X}(c))(\tau(\beta))$$ for all $c\in \check U^B$, $\lambda\in
P^+(2\Sigma)$, and $\tau(\beta)\in {\cal A}$.
\end{theorem}

Consider $c\in \check U^B$. As in the classical case,
we refer to the image ${\cal 
X}(c)$ in  $({\cal
C}(Q_{\Sigma})\check{\cal A})^{W_{\Theta}}$ as the {\it radial component} 
of $c$.

\medskip
\section{Central Elements: the Rank One Case}

In this section, we study certain central elements of $\check U$ and compute 
their radial components when the restricted
root system has rank one. 
In particular, we assume that 
$\pi^*$ contains a single root $\alpha_i$ and so 
$\Sigma^+=\{\tilde\alpha_i\}$.  Assume for the moment that 
$\pi=\pi^*$ and so $\Theta(\alpha_i)=-\alpha_i$.
Then $U$ is just  $U_q({\gs\gl\ 2})$ and is  generated by $x_i,y_i,t_i,t_i^{-1}$.
Set
$q_i=q^{(\alpha_i,\alpha_i)/2}$.
Note that
$$(q_it_i+q_i^{-1}t_i^{-1}) +(q_i-q_i^{-1})^{2}y_ix_i\leqno{(4.1)}$$
is central in $U$.
 We show that the other rank one cases  contain a similar central element. 
 
Note that $\tilde\alpha_i-\widetilde{\Theta(\alpha_i)}=2\tilde\alpha_i$. Since 
we are assuming that $\Sigma$ is reduced, it follows that $\alpha_i+\Theta(-\alpha_i)$ is 
not a root in $\Delta$. Hence 
$(\alpha_i,\Theta(-\alpha_i))=0$.
Using Araki's classification of irreducible symmetric pairs, we
have the following possibilities for $\Theta(\alpha_i)$.

\begin{enumerate}
\item[(4.2)] $\gg$ is of type $A_1$ with $\pi=\{\alpha_i\}$ and 
$\Theta(\alpha_i)=-\alpha_i$.
\item[(4.3)] $\gg$ is of type $A_1\times A_1$ with 
$\pi=\{\alpha_i,\alpha_{p(i)}\}$ and 
$\Theta(\alpha_i)=-\alpha_{p(i)}$.
\item[(4.4)]
 $\gg$ is of type $A_3$ with
$\pi=\{\alpha_1,\alpha_2,\alpha_3\}$,
$\alpha_i=\alpha_2$, and
$\Theta(\alpha_2)=-\alpha_1-\alpha_3-\alpha_2$. 
\item[(4.5)]
$\gg$ is of type $B_r$ with
$\pi=\{\alpha_1,\dots,\alpha_r\}$, $\alpha_i=\alpha_1$, and
$\Theta(\alpha_1)=- \alpha_1-2\alpha_2-\cdots -2\alpha_r$.
\item[(4.6)] $\gg$ is of type $D_r$ with
$\pi=\{\alpha_1,\dots,\alpha_r\}$, $\alpha_i=\alpha_1$, and
$\Theta(\alpha_1)=-\alpha_1-2\alpha_2-\cdots
-2\alpha_{r-2}-\alpha_{r-1}-\alpha_r$. \end{enumerate}

Recall that the $(\ad U)$ module $(\ad U)\tau(-2\mu)$ for $\mu\in
P^+(\pi)$  contains a one-dimensional subspace of $Z(\check U)$.
 (See [J1, 7.1.16-7.19 and 7.1.25] or [JL1]
for more information about $Z(\check U)$.) Moreover, there exists 
a (unique) nonzero vector $c_{\mu}$ in $(\ad U)\tau(-2\mu)\cap Z(\check U)$ 
such that
$$c_{\mu}\in \tau(-2\mu)+(\ad U_+)\tau(-2\mu).$$  
We find a ``small'' 
weight $\mu\in 
P^+(\pi)$ such that $c_{\mu}$ looks like (4.1) modulo $({\cal M}\check 
 T_{\Theta})_+\check U+\check U({\cal M}\check 
T_{\Theta})_+$. In particular, $\mu$ will satisfy the conditions of 
the following lemma. Let $w_0$ denote the longest element of the Weyl group $W$
associated to the root system of $\gg$.

\begin{lemma}  There exists $\mu\in P^+(\pi)$ such that 
$$\Theta(\mu-\alpha_i/2)=\mu-\alpha_i/2{\rm \  and\ } 
\Theta(-w_0\mu-\alpha_i/2)=-w_0\mu-\alpha_i/2.\leqno{(4.7)}$$
\end{lemma}

\noindent
{\bf Proof:}  Suppose first that $\gg$ satisfies the conditions of 
(4.2) or (4.3) above.   Then $\alpha_i/2$ is in $P^+(\pi)$ and 
$w_0\alpha_i/2=-\alpha_i/2$.   Hence we set $\mu=\alpha_i/2$ in these 
cases.  

Let $\omega_j$ be the fundamental
weight corresponding to $\alpha_j,$ for each $j$.
The remaining cases are handled below.

\medskip
\noindent {\bf Case (4.4):}    Set $\mu=\omega_1$.   We have
$\omega_1=1/4(\alpha_1+2\alpha_2+3\alpha_3)$ while 
$-w_o\omega_1=\omega_3=1/4(3\alpha_1+2\alpha_2+\alpha_3) $.   Thus (4.7)
follows since $i=2$.

\medskip
\noindent {\bf Case (4.5):} Set  $\mu=\omega_r.$  In this case
$\omega_r=-w_0\omega_r=1/2(r\alpha_r+(r-1)\alpha_{r-1}+\cdots
+3\alpha_3+2\alpha_2+\alpha_1)$.  Thus (4.7) follows since $i=1$.

\medskip
\noindent {\bf Case (4.6):}   In this case, $\mu$ can be either
$\omega_r$ or $\omega_{r-1}$. Note that
$\omega_r=1/2({{r}\over{2}}\alpha_r+{{(r-2)}\over{2}}\alpha_{r-1}+(r-2)\alpha_{r-2}+\cdots
2\alpha_2+\alpha_1) $ and
$\omega_{r-1}=1/2({{r}\over{2}}\alpha_{r-1}+{{(r-2)}\over{2}}\alpha_{r}+(r-2)\alpha_{r-2}+\cdots
2\alpha_2+\alpha_1)$.  Furthermore,
$-w_o(\omega_r)=\omega_{r-1}$ and $i=1$. $\Box$

\medskip

Recall the definition of the Harish-Chandra projection ${\cal P}$ 
given in Section 3.   We have the following
description of the Harish-Chandra  projection of central elements
$c_{\mu}$ of $\check U$ ([J1, 7.1.19 and 7.1.25]): up to a nonzero scalar,
$${\cal P}(c_{\mu})=\sum_{\nu\in P^+(\pi)}\hat\tau(-2\nu)\dim 
L(\mu)_{\nu}\leqno{(4.8)}$$ where
$$\hat\tau( \beta)=\sum_{w\in W}\tau(w\beta)q^{(\rho,w\beta)}.$$
The next result is the first step in understanding the central element 
$c_{\mu}$ when $\mu$ satisfies the conditions of Lemma 4.1.

\begin{proposition} Suppose that $\mu\in P^+(\pi)$ is chosen as in Lemma 4.1 to
satisfy condition (4.7). Then
	$$c_{\mu}\in y_i{\cal M}^-{\cal M}^+x_i+{\cal P}(c_{\mu})+  ({\cal M}\check 
	T_{\Theta})_+\check U+\check U({\cal M}\check T_{\Theta})_+\ .$$
	Moreover, ${\cal P}(c_{\mu})\in {\cal C}[T_{\geq}]\tau(-2\mu)$.
\end{proposition}

\noindent{\bf Proof:} The last statement is an immediate
consequence of (4.8). Note that 
$c_{\mu}$ is a sum of zero weight vectors in $({\rm ad}\ U)\tau(-2\mu).$ By [J2, Theorem
3.3], we can construct a basis for the zero weight space of $({\rm ad}\ U)\tau(-2\mu)$
consisting of vectors in sets of the form
$$a_{-\beta}b_{\beta}\tau(-2\mu)+\sum_{\gamma\in Q^+(\pi)}\sum_{0<\gamma'<\beta-\gamma}
G^-_{-\gamma'}U^+_{\gamma'}\tau(-2\mu+2\gamma)$$
where   $a_{-\beta}$ is a weight vector of weight $-\beta$ in $(\ad U^-)\tau(-2\mu)$ and
$b_{\beta}$ is a weight vector of weight $\beta$ in $(\ad U^+)\tau(-2\mu)$. It further
follows  from [J1, 7.1.20], that the weights $\beta$ appearing in the above expression
satisfy $\beta\leq \mu-w_0\mu$. Hence
$$\eqalign{c_{\mu}&\in 
\sum_{\gamma\in
Q^+(\pi_{\Theta})}G^-_{-\mu+w_0\mu+\gamma}U^+_{\mu-w_0\mu-\gamma}T_{\Theta}\tau(-2\mu)
\cr&+\sum_{\{\beta|0\leq\tilde\beta<\tilde\mu-\widetilde{w_0\mu}\}}G^-_{-\beta}U^+_{\beta}T_{\geq}\tau(-2\mu).\cr}\leqno{(4.9)}$$

 By choice of $\mu$, we have that $\mu-w_0\mu-\alpha_i\in Q(\pi_{\Theta})$. Hence
$0\leq\tilde\beta<\tilde\mu-\widetilde{w_0\mu}$ forces $\beta\in Q^+(\pi_{\Theta})$. Thus (4.9) implies that 
$$c_{\mu}\in y_it_i{\cal M}^-{\cal M}^+x_iT_{\Theta}\tau(-2\mu)+{\cal
P}(c_{\mu})+  {\cal M}_+\check 
	 U+\check U{\cal M}_+.\leqno{(4.10)}$$ The assumption on $\mu$ further implies that 
$t_i\tau(-2\mu)=\tau(-2\mu+\alpha_i)\in \check T_{\Theta}$.  Hence 
$$t_i\tau(-2\mu)=1+(t_i\tau(-2\mu)-1)\in  {\cal C}+U{\cal C}[\check 
T_{\Theta}]_+.$$
The lemma now follows from this expression and (4.10).  $\Box$

\medskip

The next  lemma simplifies the component of $c_{\mu}$ coming from 
$y_i{\cal M}^-{\cal M}^+x_i$.

\begin{lemma} 
$ y_i{\cal M}^-{\cal M}^+x_i\subseteq {\cal C}  y_ix_i +
{\cal M}_+U+U{\cal M}_+.$
\end{lemma}

\noindent
{\bf Proof:}  Consider an element $y_ibcx_i$ of $U$ where $b\in {\cal M}^-$ and $c\in {\cal M}^+$.
Using the relations of $U$, we can rewrite this element as  
a sum of terms of the form $y_ic'b'd'x_i$ where
$c'\in {\cal M}^+$, $b'\in {\cal M}\cap U^0$, and $d'\in {\cal M}^-.$   Since $\alpha_i\notin \pi_{\Theta}$, we have
$y_ic'b'd'x_i=c'y_ib'x_id'$.   So if either $c'$    is in  ${\cal
M}_+^+$ or $d'$ is in
${\cal M}^-_+$, then $y_ic'b'd'x_i\in {\cal M}^+_+U+U{\cal
M}^-_+\subseteq {\cal M}_+U+U{\cal M}_+$. If neither of these 
conditions hold, we may assume that $c'b'd'=b'$ which is an element 
of ${\cal M}\cap U^0$.   Then $y_ib'x_i\in {\cal C}y_ix_i+y_ix_i({\cal 
M}\cap U^0)_+\subseteq {\cal C}  y_ix_i +
{\cal M}_+U+U{\cal M}_+.$
$\Box$
\medskip

 The tensor product 
decomposition (2.1) implies the following direct sum decomposition:
$$U^0={\cal C}[T']\oplus U^0{\cal C}[T_{\Theta}]_+.\leqno{(4.11)}$$
Recall the direct sum decomposition (3.5) using $\check {\cal A}$ instead of $T'$.
Let ${\cal P}'_{\cal A}$ denote the map from $\check U$ to
${\cal C}[\check {\cal A}]$ defined by composing the Harish-Chandra map ${\cal P}$ with the projection of $\check U^0$ into 
${\cal C}[\check{\cal A}]$ using (3.5).    Similarly, let ${\cal P}'$ 
denote the map from $U$ onto ${\cal C}[T']$
which is  the composition
of
${\cal P}$ with   the 
projection onto ${\cal C}[T']$ using (4.11). Note that  $${\cal 
P}'_{\cal A}(a)=\sum_ma_m\tau(\tilde\alpha_i)^m{\rm \  if
\ and\ only\ if\ } 
{\cal P}'(a)=\sum_ma_mt_i^m\leqno{(4.12)}$$ for all $a\in U$.

Recall the dotted Weyl group action on $\check T$ defined by
$$w. \tau(\mu)q^{(\rho,\mu)}=\tau(w\mu)q^{(\rho,w\mu)}\leqno{(4.13)}$$
for all $\tau(\mu)\in \check T$ and $w\in W$.  By [J1,
7.1.17 and 7.1.25], the image of $Z(\check U)$ under ${\cal P}$ is
contained in    
${\cal C}[\check T]^{W.}$.  
Define the dotted action of $W_{\Theta}$
on
$\check {\cal A}$ using the same formula as in (4.13) where now $w$ is an element of $W_{\Theta}$ and $\tau(\mu)\in
\check {\cal A}$.

\begin{lemma} For all $\mu\in P^+(2\Sigma)$,  the image 
of $c_{\mu}$ under ${\cal P}'_{\cal A}$ is
 invariant under the dotted action of $W_{\Theta}.$  Moreover,
if
$\mu$  satisfies the conditions of (4.7) then
${\cal P}'(c_{\mu})$ is a  scalar multiple of 
$q_i^{(\rho,\tilde\alpha_i)}t_i+q_i^{-(\rho,\tilde\alpha_i)}t_i^{-1}$.

\end{lemma}
 
\noindent
{\bf Proof:}  Recall that
$w_0$ is the longest  element of $W$ and let $w_0'$ denote the longest element of the Weyl 
group $<s_{\alpha}|\alpha\in \pi_{\Theta}>$.   Set $w=w_0'w_0$. 
 By [L3, Section 7] or checking (4.2)-(4.6) directly, we see that  
$w_0'\alpha_i=\Theta(-\alpha_{i})$ and $w_0'\Theta(\alpha_{i})=-\alpha_i$.
Further checking the possibilities for $w_0$ in (4.2)-(4.6) yields that
$w\alpha_i=\Theta(\alpha_{i})$ and $w\Theta(\alpha_{i})=\alpha_i$.  
 Hence $w\tilde\alpha_i=-\tilde\alpha_i$ and   we may identify $W_{\Theta}$ with the subgroup 
$\langle w\rangle$ of $W$.

It is straightforward to check using (4.2) 
through (4.6) that $w_0$ sends $\alpha_i$ to $-\alpha_i$.  Furthermore,  
$w_0$ sends a simple root in $\pi_{\Theta}$ to the negative of a 
simple root in $\pi_{\Theta}$.  It follows  that $w$ permutes the
elements  of $\pi_{\Theta}$. Thus
$q^{(\rho,\gamma)}=q^{(\rho,w\gamma)}$  for all $\gamma$ such that
$\Theta(\gamma)=\gamma$. Hence
$w. \tau(\mu)=\tau(w\mu)$ for all $\mu$ such that 
$\Theta(\mu)=\mu$. It follows that ${\cal C}[\check T_{\Theta}]_+$
is invariant under the dotted action of $W_{\Theta}$.
Thus $${\cal P}'_{\cal A}({\cal C}[\check T]^{W.})\subseteq {\cal C}[\check{\cal A}]^{W_{\Theta}.}$$
which proves the  first assertion of the lemma.

Now assume that $\mu$ satisfies
(4.7).   It follows that 
${\cal P}'_{\cal A}(\tau(-2\mu))=\tau(\tilde\alpha_i)$.  By Proposition 4.2, 
${\cal P}(c_{\mu})\in {\cal C}[T_{\geq}]\tau(-2\mu)$.   Hence (2.1) and (2.10) imply that
${\cal P}'_{\cal A}(c_{\mu})\in {\cal C}[{\cal A}_{\geq}]\tau(\tilde\alpha_i)$.
The second assertion now follows from (4.12) and   the 
fact that the only elements   of ${\cal C}[{\cal A}_{\geq}]\tau(\tilde\alpha_i)$ invariant under the dotted action of
$W_{\Theta}$ are scalar multiples of
$q_i^{(\rho,\tilde\alpha_i)}\tau(\tilde\alpha_i)+ q_i^{-(\rho,\tilde\alpha_i)}\tau(\tilde\alpha_i)^{-1}$. 
$\Box$

\medskip

 An immediate consequence of (4.8) is that
$z^{\lambda}({\cal P}(c_{\mu}))\neq 0$, and more importantly, ${\cal P}'(c_{\mu})\neq 0$ for
any choice of $\lambda$ and
$\mu$.   This fact is used in the  next result.  In particular,  we show that when $\mu$
is chosen as in Lemma 4.1,  then $c_{\mu}$ looks like the central element of $U_q(\gs\gl\
2)$ described in (4.1).

 \begin{theorem} Assume that $\mu$   satisfies the 
conditions of (4.7). Let $a$ be the nonzero scalar guaranteed by Lemma 
4.4 (and the above comments) such that 
${\cal
P}'(c_{\mu})=a(q_i^{(\rho,\tilde\alpha_i)}t_i+q_i^{-(\rho,\tilde\alpha_i)}t_i^{-1}).
$ Then 
$$\eqalign{c_{\mu}\in
a[q_i^{(\rho,\tilde\alpha_i)}t_i+&q_i^{-(\rho,\tilde\alpha_i)}t_i^{-1} +
(q_i-q_i^{-1})(q_i^{(\rho,\tilde\alpha_i)}-q_i^{-(\rho,\tilde\alpha_i)})y_ix_i]\cr +&\ ({\cal
M}\check T_{\Theta})_+\check U+ \check U({\cal M}\check T_{\Theta})_+.\cr}$$
\end{theorem}

\noindent 
{\bf Proof:}   By Proposition 4.2 and Lemma 4.3, there 
exists  $c\in {\cal C}$ such that   
$$c_{\mu}-( cy_ix_i+b)\in ({\cal 
M}\check T_{\Theta})_+\check U+\check U({\cal 
M}\check T_{\Theta})_+\leqno{(4.14)}$$
where   $$b=
a(q_i^{(\rho,\tilde\alpha_i)}t_i+q_i^{-(\rho,\tilde\alpha_i
)}t_i^{-1}).\leqno{(4.15)}$$ 

Consider the
${\bf C}$ algebra automorphism
$\psi$ of
$U$ defined by
$\psi(x_j)=y_jt_j$, $\psi(y_j)=t_j^{-1}x_j$, $\psi(t)=t$ and $\psi(q)=q^{-1}$
for all $1\leq j\leq n$ and $t\in \check T$.  A straightforward check shows that
$\psi((\ad x_i)b)=(\ad y_it_i)\psi(b)$ and $\psi((\ad y_i)b)=(\ad 
t_i^{-1}x_i)\psi(b)$ for all $b\in \check U$.
In particular, if $x\in (\ad U)\tau(\mu)$, then so is $\psi(x)$.

Recall that $c_{\mu}$ has been scaled  so that $c_{\mu}\in
\tau(-2\mu)+(\ad U_+)\tau(-2\mu)$.   Hence $\psi(c_{\mu})\in \tau(-2\mu)+(\ad
U)_+\tau(-2\mu)$.   Since $(\ad U_+)\check U\cap Z(\check U)=0$, it follows that
$\psi(c_{\mu})=c_{\mu}$.  Therefore, applying $\psi$ to (4.14) using 
the form of $b$ given in (4.15) yields $c_{\mu}$ is an element of 
  $$a
(q_i^{-(\rho,\tilde\alpha_i)}t_i+q_i^{(\rho,\tilde\alpha_i
)}t_i^{-1})+ct_i^{-1}x_iy_it_i  + ({\cal M}\check T_{\Theta})_+\check 
U+\check U ({\cal
M}\check T_{\Theta})_+.$$  
  Thus
$$\eqalign{c_{\mu}&\in
a(q_i^{-(\rho,\tilde\alpha_i)}t_i+q_i^{(\rho,\tilde\alpha_i
)}t_i^{-1})+c(q_i-q_i)^{-1}(t_i-t_i^{-1})+ cy_ix_i \cr&+  ({\cal 
M}\check T_{\Theta})_+\check U+\check U ({\cal M}\check T_{\Theta})_+.\cr}$$
Note that $q_i^{-(\rho,\tilde\alpha_i)}t_i+q_i^{(\rho,\tilde\alpha_i
)}t_i^{-1}$ is not invariant with respect to the dotted ${W_{\Theta}}$ 
action.  Hence $c$ must be 
nonzero.  Moreover, in order for $a(q_i^{-(\rho,\tilde\alpha_i)}t_i+q_i^{(\rho,\tilde\alpha_i
)}t_i^{-1})+c(q_i-q_i)^{-1}(t_i-t_i^{-1})$ to be 
invariant under the dotted $W_{\Theta}$ action, we must have
  $c =(q_i-q_i^{-1})(q_i^{(\rho,\tilde\alpha_i)}-q_i^{-(\rho,\tilde\alpha_i)})a$.
$\Box$

\medskip

  Let $\chi$ be
the Hopf algebra automorphism (which restricts to the identity on
${\cal M}T$) defined by
$$\chi(x_i)=q^{-1/2(\rho,\Theta(\alpha_i)-\alpha_i)}x_i$$
and
$$\chi(y_i)=q^{1/2(\rho,\Theta(\alpha_i)-\alpha_i)}
y_i.$$
By [L4, Section 5], we may assume that $B'=\chi(B)$.

Note that $B_+$ contains $$B_i=y_it_i+\tilde\theta(y_i)t_i$$
and $B'_+$ contains
$$B_i'=y_it_i+q^{-(\rho,\Theta(\alpha_i)-\alpha_i)}\tilde\theta(y_i)t_i .$$
By [L4, Lemma 5.1], we have
$$q^{(\rho,\Theta(\alpha_i)+\alpha_i)}
\tilde\theta(y_{p(i)})t_{p(i)}^{-1}x_{p(i)}\in t_i^{-1}x_i\tilde\theta(y_i)
+ {\cal
M}^+_+ U+U{\cal M}^+_+.$$  When $i=p(i)$, it follows that
$$q^{(\rho,\Theta(\alpha_i)+\alpha_i)}
\tilde\theta(y_{i})t_{i}^{-1}x_{i}\in t_i^{-1}x_i\tilde\theta(y_i)
+ {\cal
M}^+_+ U+ U{\cal M}^+_+.\leqno{(4.16)}$$

Assume for the moment that $i\neq p(i)$. The assumption that $\Sigma$ is reduced ensures that 
 $\Theta(\alpha_i)=-\alpha_{p(i)}$ and 
$(\alpha_i,\alpha_{p(i)})=0$.  Moreover, checking cases (4.2)-(4.6) yields that 
$\tilde\theta(y_{i})=t_{p(i)}^{-1}x_{p(i)}$.   Hence
$t_i^{-1}x_i\tilde\theta(y_i)=t_i^{-1}x_it_{p(i)}^{-1}x_{p(i)}=t_{p(i)}^{-1}x_{p(i)}t_i^{-1}x_i=\tilde\theta(y_i)t_i^{-1}x_i$.
It follows that (4.16) holds when $i\neq p(i)$ as well.

The next lemma will allow us to compute ${\cal X}(y_ix_i)$. This, in
turn, will be used to   compute the image of $c_{\mu}$ under ${\cal X}$
where $\mu$ satisfies the conditions of (4.7).

\begin{lemma} Let  $\tau(\lambda)\in T$  such that $s(\alpha_i,\alpha_i)=
 (\lambda,\alpha_i)=(\lambda,\Theta(\alpha_i))\neq 0$.
 Then
$$y_ix_i\tau(\lambda)+q_i^{-4s}{{(t_i-t_i^{-1})}\over{(q_i^{-4s}-1)(q_i-q_i^{-1})}} )\tau(\lambda)$$ is an element in $
B_+U+U{B'}_+$.
\end{lemma}

\noindent
{\bf Proof:}  Set
$a_i=2(\rho,\Theta(\alpha_i)-\alpha_i)/(\alpha_i,\alpha_i)$.
Note that
$$B_it_i^{-1}x_i\tau(\lambda)=y_ix_i\tau(\lambda)+ 
\tilde\theta(y_{i})t_it_i^{-1}x_{i}\tau(\lambda).$$
Now $2(\rho,\Theta(\alpha_i)+\alpha_i)/(\alpha_i,\alpha_i)=a_i+2$.
Hence by (4.16),
$$y_ix_i\tau(\lambda)+q_i^{-a_i}t_i^{-1}x_i\tilde\theta(y_i)t_i\tau(\lambda)
-B_it_i^{-1}x_i\tau(\lambda)\in {\cal M}_+U+U{\cal 
M}_+.$$
On the other hand $$\eqalign{t_i^{-1}x_i\tau(\lambda)B_i'=&q_i^{-2s}x_iy_i\tau(\lambda)+
q_i^{2s{-a_i}}t_i^{-1}x_i\tilde\theta(y_i)t_i\tau(\lambda)\cr
=&q_i^{-2s}y_ix_i\tau(\lambda)+q_i^{-2s}{{(t_i-t_i^{-1})}\over{(q_i-q_i^{-1})}}\tau(\lambda)+
q_i^{2s{-a_i}}t_i^{-1}x_i\tilde\theta(y_i)t_i\tau(\lambda).\cr}$$

Thus 
$$(q_i^{-4s}-1)y_ix_i\tau(\lambda)+q_i^{-4s}{{(t_i-t_i)}\over{(q_i-q_i^{-1})}}\tau(\lambda)\in B_+U+UB'_+.$$
$\Box$
\medskip

Set $\tilde t_i=\tau(\tilde\alpha_i)$.  We are  now ready to compute the radial components of the central 
elements described in Theorem 4.5.

\begin{theorem} Let $\mu$ satisfy the conditions of (4.7).
Let $a$ be the nonzero scalar such that ${\cal
P}'(c_{\mu})=a(q_i^{(\rho,\tilde\alpha_i)}t_i+q_i^{-(\rho,\tilde\alpha_i)}t_i^{-1})$.	Then
 $${\cal X}(c_{\mu})=a[q_i^{- (\rho,\tilde\alpha_i)}\tilde
 t_i(q_i^{2(\rho,\tilde\alpha_i)}z^{2\tilde\alpha_i}-1)
+q_i^{(\rho,\tilde\alpha_i)}\tilde
t_i^{-1}(q_i^{-2(\rho,\tilde\alpha_i)}z^{2\tilde\alpha_i}-1)]
(z^{2\tilde\alpha_i}-1)^{-1}.$$
\end{theorem}

\noindent {\bf Proof:} By Theorem 4.5 
and Lemma 4.6, $a^{-1}c_{\mu}\tau(\lambda)$ is  an element in 
$$(q_i^{(\rho,\tilde\alpha_i)}t_i+
q_i^{-(\rho,\tilde\alpha_i)}t_i^{-1})\tau(\lambda)+
{{(q_i^{(\rho,\tilde\alpha_i)}-q_i^{-(\rho,\tilde\alpha_i)})(t_i-t_i^{-1})\tau(\lambda)}
\over {(q_i^{4s}-1)}}+B_+\check U+\check UB_+'$$ for all $\tau(\lambda)\in {\cal 
A}$.  This set simplifies to 
$$
{{(q_i^{4s+2(\rho,\tilde\alpha_i)}-1)q_i^{-(\rho,\tilde\alpha_i)}\tilde t_i+
(q_i^{4s-2(\rho,\tilde\alpha_i)}-1)q_i^{(\rho,\tilde\alpha_i)}\tilde t_i^{-1}}
\over{(q_i^{4s}-1)}}\tau(\lambda)
+B_+\check U+\check UB_+'.$$

Set $Y$ to be the right-hand expression in the statement of the theorem. 
The desired formula now follows from the fact that $Y$ is the unique
element in ${\cal C}(Q_{\Sigma})\check {\cal A}$ such that  
$$a^{-1}Y\cdot\tau(\lambda) =\left( {{(q_i^{4s+2(\rho,\tilde\alpha_i)}-1)q_i^{-(\rho,\tilde\alpha_i)}\tilde t_i+
(q_i^{4s-2(\rho,\tilde\alpha_i)}-1)q_i^{(\rho,\tilde\alpha_i)}
\tilde 
t_i^{-1}}\over{(q_i^{4s}-1)}}\right)\tau(\lambda) 
$$
for all $\tau(\lambda)\in {\cal A}$. $\Box$

\section{ Graded Zonal Spherical Functions}

 Recall the $\ad U$ filtration on $U$ defined in    [J1, 7.1.1].  In this
section, we use a modified version of this filtration that is chosen so
that the associated graded ring of $U$ contains $B$ as a
subalgebra.

Define a degree function on $U$ by
\begin{enumerate}
\item[(5.1)] $\deg x_i=\deg y_it_i=0$ for all $i$, $1\leq i\leq n$.
\item[(5.2)] $\deg t_i^{-1}=1$ for all $i$ such that
$\alpha_i\in\pi\setminus\pi_{\Theta}$ \item[(5.3)] $\deg t=0$ for all $t\in
T_{\Theta}$.
\end{enumerate}
Let ${\cal F}$ denote the filtration on $U$ defined by the above
degree function.  Write $\gra U$ for the associated graded algebra
with respect to this filtration.  Note that elements of $U^+$, $G^-$, 
and  ${\cal M}$ are all in degree zero.  Moreover the relations 
satisfied by the elements of $U^+$ (resp. $G^-$, ${\cal M}$) are 
homogeneous of degree $0$. Therefore, the map 
$a\mapsto \gra a$ defines an isomorphism between $U^+$ and $\gra 
U^+$, between $G^-$ and $\gra G^-$, and between ${\cal M}$ and $\gra 
{\cal M}$.     Using this isomorphism, we 
write $U^+$ for $\gra U^+$,  $G^-$ for $\gra G^-$, and ${\cal M}$ for 
$\gra {\cal M}$. Furthermore, if $S$ is a subset and $a$ is an element  of $U^+$, $G^-$, 
or ${\cal M}$ then we simply write $a$ for $\gra a$ and $S$ for $\gra S$.

The next lemma shows that a similar 
identification holds for $B$.  In particular, we may  identify $B$ with 
$\gra B$ as a subalgebra of $\gra U$.

\begin{lemma} For all $b\in B$,   $\deg b=0$.   Thus the filtration
${\cal F}$ restricts to the trivial filtration on $B$ and the map 
$a\mapsto \gra a$ defines an isomorphism between $B$ and $\gra 
B$.  
\end{lemma}

\noindent {\bf Proof:} Let $\tilde B$ denote the algebra generated 
freely over ${\cal M}^+T_{\Theta}$ by elements $\tilde B_i$, $1\leq 
i\leq n$. By [L3, Theorem 7.4], there is a homomorphism from $\tilde B$ onto $B$ 
which is the identity on ${\cal M}T_{\Theta}$, sends $\tilde B_i$ to 
$B_i$ for $\alpha_i\notin\pi_{\Theta}$, and sends $\tilde B_i$ to 
$y_it_i$ for $\alpha_i\in \pi_{\Theta}$.

By (5.1) and (5.3), we have that $\deg a=0$ for
all $a\in {\cal M}T_{\Theta}$.  Consider
$B_i=y_it_i+\tilde\theta(y_i)t_i$ for some $\alpha_i\notin\pi_{\Theta}$.
We have $\deg y_it_i=0$. Recall that $\tilde\theta$ is a particular lift of the involution $\theta$ to a ${\bf C}$
algebra automorphism of $U$. It follows  from the   explicit  description of $\tilde\theta$ given in   
of [L3, Theorem 7.1]  that $\tilde\theta(y_i)\in
U^+\tau(\Theta(\alpha_i))$.  Since
$\tau(\Theta(\alpha_i))t_i=\tau(\Theta(\alpha_i)+\alpha_i)\in
T_{\Theta}$, we have that $\deg \tilde\theta(y_i)t_i=0$ as well.   So
the generators of $B$ are all in degree $0$. By the previous 
paragraph, the relations satisfied by these generators are all 
homogeneous of degree $0$.   Hence all elements in $B$ have degree 
$0$. $\Box$

\medskip

Note that not all elements of ${\cal C}[T]$ are in degree $0$.  So 
${\cal C}[T]$ does not naturally identify with its graded image in the 
same way as  the subalgebras discussed above.  However,
 the algebra map 
induced by $t_i \mapsto \gra t_i$ for   $1\leq i\leq n$  does define an
isomorphism from 
${\cal C}[T]$ to $\gra {\cal C}[T]$.  Thus any $\gra T$ module 
inherits the structure of  a $T$ module via this isomorphism.

Set $\gh^*_{\Theta}=\{\lambda\in \gh^*|\Theta(\lambda)=-\lambda\}$.
Consider $\lambda\in \gh^*_{\Theta}$ and let $v_{\lambda}$ be a
(left)
$T$ weight vector of  weight $\lambda$.  Note that $tv_{\lambda}=v_{\lambda}$ for all $t\in 
T_{\Theta}$.
 Make ${\cal
C}v_{\lambda}$ into a $\gra {\cal M}TU^+$ module by insisting that
${\cal M}_+v_{\lambda}=U^+_+v_{\lambda}=0$.  Define the left $\gra
U$ module $\bar M({\lambda})$ by
$$\bar M({\lambda})=\gra U\otimes_{(\gra {\cal M}TU^+)}v_{\lambda}.$$
By Lemma 2.1 and Theorem 2.2, $\bar M({\lambda})=\gra N^-\otimes v_{\lambda}$ as $\gra 
N^-$ modules. Since $N^-$ is a subalgebra of $G^-$, 
$ N^-$ can be identified with $\gra N^-$ via the obvious map.

The algebra $G^-$ can be given the structure of a $U^+$ module as 
in    [J1, Sections 5.3 and 7.1].  In particular,  let $x_i'$ and $x_i''$
be functions on $G^-$ such 
that
$$(\ad x_i) m=x_i'(m)+x_i''(m)t_i^2$$ for all $m\in G^-$ and for 
all $i$ such that $1\leq i\leq n$.   Given $i$ 
such that $1\leq i\leq n$, the action of $x_i$  on the
element $m\in G^-$ is defined by
$$x_i*m=x_i'(m).$$

\begin{lemma} $N^-$ is a $U^+$ submodule of $G^-$.   Moreover, 
$\bar M({\lambda})\cong N^-$ as $U^+$ modules for all  $\lambda\in 
\gh^*_{\Theta}$.
\end{lemma}

\noindent
{\bf Proof:}  First, note that by  [L3, Section 6], $N^-$ is an $\ad {\cal
M}^+$ module.   In particular, $(\ad x)n\in N^-$ for all $x\in {\cal 
M}^+$ and $n\in N^-$.   It follows that  $x_i''(n)=0$  and
$x_i*n=(\ad x_i)n$  is an element of $N^-$ for all $\alpha_i\in
\pi_{\Theta}$ and $n\in N^-$.  Thus ${\cal M}*N^-\subseteq N^-$.

Now $N^-$ is generated by elements of the form $(\ad y)y_jt_j$
([L3, Section 6]) where $y\in {\cal M}^-$ and $\alpha_j\notin\pi_{\Theta}$.
Suppose that $y\in {\cal M}^-$.
By the defining relations of $U$, $(\ad x_k)((\ad y)y_jt_j)=(\ad y)(\ad x_k)y_jt_j)$ for
all $k$ such that $\alpha_k\notin\pi_{\Theta}$.    It follows that
$$x_k'((\ad y)y_jt_j)=(\ad y)\left({{-\delta_{kj}}\over
{q_j-q_j^{-1}}}\right).$$  Thus $x_k'((\ad y)y_jt_j)=0$ for $y\in {\cal
M}^-_+$ while $x_k'( y_jt_j)$ is a scalar.  Therefore $x_k'(n)\in N^-$
for all $n\in N^-$ and $\alpha_k\notin \pi_{\Theta}$.  The fact that 
$G^-$ is generated by $N^-$ and ${\cal M}^-$ ([L3, Section 6])
yields that $N^-$ is a $U^+$ submodule
of $G^-$.

Fix $n\in N^-$. By (2.4), $(\gra x_i) (n\otimes v_{\lambda})=\gra((\ad
x_i)n)\otimes v_{\lambda}$ since $x_iv_{\lambda}=0$.   Furthermore,
$\gra (\ad x_i)n= x_i'(n)$ for all $\alpha_i\notin\pi_{\Theta}$ by
the definition of the filtration ${\cal F}$. Now consider 
$\alpha_i\in \pi_{\Theta}$.  It follows that 
$(\ad x_i)n= x_i'(n)$  and so $(\gra x_i)( n\otimes 
v_{\lambda})= x_i'(n)\otimes v_{\lambda}$ in this case as well.  Thus,  
the map $n\mapsto n\otimes v_{\lambda}$ is an isomorphism 
of $U^+$ modules under the identification of  $U^+$ with $\gra U^+$.
$\Box$
\medskip

 Recall that  $N^-$ is an $\ad T$ 
module.  Since $v_{\lambda}$ is a $T$ weight vector,  it follows that  
$\bar M({\lambda})$ is a direct sum of
its $T$ weight spaces with highest weight equal to $\lambda$. 
Moreover, suppose $b$ is a weight vector in either $G^-$ or $U^+$, and 
$m$ is  a weight vector in $\bar M({\lambda})$.  Note that  the weight 
of $(\gra b)m$ 
is equal to the sum of the weights of $\gra b$ and $m$.  

\begin{lemma}  $\bar M({\lambda})$ is a simple $\gra U$ module for each 
$\lambda\in {\gh}_{\Theta}^*$.
\end{lemma}

\noindent {\bf Proof:}   Observe that $\bar M({\lambda})$ is a cyclic $\gra
U$ module generated by $v_{\lambda}$.   So it is sufficient to
show that $1\otimes v_{\lambda}$ is in the $\gra U$ module
generated by $n\otimes v_{\lambda}$ for any weight vector
$n\in N^-$.  

By the discussion preceding the lemma,
$(\gra b)m$ has higher weight than $m$ for any nonzero 
weight vector  $\gra b$ in  $\gra U^+_+$ and nonzero 
weight vector $ m$ in $\bar M({\lambda})$.  Hence, it is
enough to show that the only vectors in $\bar M({\lambda})$ annihilated
by $ U^+_+$ are scalar multiples of $1\otimes v_{\lambda}$.
This follows from Lemma 5.2 and the fact that  the only $U^+$ invariant vectors of $G^-$ are scalar
multiples of $1$ ([JL2, Lemma 4.7(i)]).   
$\Box$

\medskip

 Now let $v_{\lambda}^r$ be a right $T$ weight vector of weight 
$\lambda$.
We can make ${\cal C}v_{\lambda}^r$ into a right $\gra U$ module as follows.
Set $v_{\lambda}^r{\cal M}_+=v_{\lambda}^rG^-_+=0$ and set 
$$\bar M({\lambda})^r={\cal C} v_{\lambda}^r\otimes_{\gra {\cal
M}TG^-}\gra U.$$ Replacing $x_i$ with $y_it_i$ and $N^-$ with $N^+$,  we
can give
$N^+$ a $G^-$ module structure analogous to the $U^+$ module structure of
$N^-$.  Furthermore, as in  Lemma 5.2, 
$\bar M({\lambda})^r$ is isomorphic to $N^+$ as a $G^-$ module. 
Moreover,   as in Lemma 5.3, $\bar M({\lambda})^r$ is a simple right $\gra
U$ module.

\begin{lemma} Let $\lambda$ and $ \lambda'$ be elements in $\gh^*_{\Theta}.$
The map which sends $$m\otimes v_{\lambda}\mapsto m\otimes v_{\lambda'}$$
for all $m\in N^-$ defines an isomorphism from  $\bar M({\lambda})$ onto $
\bar M({\lambda'})$ as $ U^+$ modules, $  G^-$ modules,   and $B$ modules. Similarly, the map which sends
$$v_{\lambda}^r\otimes m\mapsto  v_{\lambda'}^r\otimes m$$
for all $m\in N^+$ defines an isomorphism from  $\bar M({\lambda})^r$
onto
$\bar  M({\lambda'})^r$
as $ U^+$ modules, $  G^-$ modules,   and $B$
modules.
\end{lemma} 

\noindent
{\bf Proof:} We prove the first assertion.   The second follows 
in a similar fashion.  
The proof of Lemma 5.2 shows that the map $n\mapsto 
n\otimes v_{\lambda}$ is a $U^+$ module isomorphism from $N^-$ 
onto $\bar M({\lambda})$.   This isomorphism is independent of
$\lambda$.    Hence $\bar M({\lambda})$ is isomorphic to
$\bar M(\lambda')$ as
$U^+$ modules  for all $\lambda,\lambda'\in \gh^*_{\Theta}$. 

Now $N^-$ is a ${\cal C}[T_{\Theta}]$ module via the adjoint action.   
Recall that
$tv_{\lambda}=v_{\lambda}$ for all $t\in 
T_{\Theta}$.
Thus, $(\gra t)( n\otimes v_{\lambda})=(\gra 
tnt^{-1})\otimes v_{\lambda}$ for all $t\in T_{\Theta}$ and $n\in N^-$. Hence  
$n\mapsto n\otimes v_{\lambda}$ is an isomorphism of ${\cal 
C}[T_{\Theta}]$ modules. This isomorphism is independent of the 
choice of  
$\lambda\in \gh^*_{\Theta}$.  Therefore
$\bar M({\lambda})\cong \bar M({\lambda'})$ as ${\cal C}[T_{\Theta}]$
modules for all 
$\lambda,\lambda'\in \gh^*_{\Theta}$.

 Recall that $\bar M({\lambda})= \gra U\otimes v_{\lambda}=
N^-\otimes v_{\lambda}$ as left $ N^-$ modules. In particular, the action of an element in 
$ N^-$ on  $\bar M({\lambda})$ just corresponds to left multiplication by
that element. Thus the action of $ N^-$ on $\bar M({\lambda})$ is 
independent of $\lambda$.

Now consider $y_it_i$ where $\alpha_i\in
\pi_{\Theta}$.  In particular, $y_it_i$ is an element of ${\cal M}$.
Recall that $N^-$ is $\ad y_it_i$ invariant while $
y_it_i\otimes v_{\lambda}=0$.  Hence, given
$n\in N^-$, we have $$( y_it_in)\otimes v_{\lambda}= ((\ad y_it_i)n)\otimes v_{\lambda}.$$  Thus the action
of $ y_it_i$ on $\bar M({\lambda})$ corresponds to the action of $(\ad
y_it_i)$ on
$N^-$ for all $\alpha_i\in \pi_{\Theta}$. As mentioned   earlier ([L3, Section 6]) $G^-$ is generated by 
${\cal M}^-$ and
$N^-$.  It follows that the action 
of $ G^-$ on $\bar M({\lambda})$ is independent of $\lambda$. 

Recall the identification of $B$ with $\gra B$.  It follows from the
proof of Lemma 5.1 that
$\gra B$ is contained in the subalgebra of $\gra U$ generated by $\gra U^+$,
$\gra G^-$, and $\gra {\cal C}[T_{\Theta}]$. Thus the isomorphism of
$\bar M({\lambda})$ to $\bar M({\lambda'})$ as $B$ modules follows from
their isomorphism as $ U^+$, $ G^-$, and $ {\cal
C}   [T_{\Theta}]$ modules. 
$\Box$

\medskip

Let $\bar M({\lambda})^*$ denote the dual of $\bar M({\lambda})$ given its
natural right $\gra U$ module structure.  The locally finite $T$ part,
$F_T(\bar M({\lambda})^*)$, of $\bar M({\lambda})^*$ is the direct sum of
its
$T$ weight spaces.  Note further that the $\beta$ weight space of
$F_T(\bar M({\lambda})^*)$ is the dual of the $\beta$ weight space of
$\bar M({\lambda})$.  In particular, $$\dim
F_T(\bar M({\lambda})^*)=\bar M({\lambda})_{\beta}$$ for all $\beta$.
Let $v_{\lambda}^*$ be a nonzero vector in 
$F_T(\bar M({\lambda})^*)_{\lambda}$.
Then $v_{\lambda}^*$ generates a simple $\gra U$ module isomorphic to
$\bar M(\lambda)^r$. A comparison of the dimension of the weight spaces
yields $F_T(\bar M(\lambda)^*)=v_{\lambda}^*\, \gra U$.  

Let $\hat
 M(\lambda)^r$ denote the completion of $\bar M(\lambda)^r$ consisting of
possibly infinite sums of distinct weight vectors
$\sum_{\gamma}a_{\gamma}$ for $a_{\gamma}\in (\bar M(\lambda)^r)_{\gamma}$. 
We can identify
$\bar M(\lambda)^*$
with $\hat M({\lambda})^r$.   Similar considerations allow us to identify 
$\bar M(\lambda)^{r*}$ with the completion $\hat M({\lambda})$ consisting of 
possibly infinite sums of distinct weight vectors in $\bar M(\lambda)$.

For the remainder of the paper, given $u\in \gra U$  we  write $uv_{\lambda}$
for the element
$u(1\otimes  v_{\lambda})$ of $\bar M(\lambda)$ and   
$v_{\lambda}^ru$ for  the element $(v_{\lambda}^r\otimes 1)u$ of $\bar
M(\lambda)^r$.  In light of  the  isomorphisms of Lemma 5.4, we often abbreviate   
$(\gra a)w $ as $aw$ and $w'(\gra a)$ as $w'a$ 
for $a\in G^-\cup B\cup U^+$,
$w\in \bar M(\lambda)$, and $w'\in \bar M(\lambda)^r$.

\begin{lemma} Let $V$ be a finite dimensional simple right 
$B$-module and $W$
be a finite dimensional simple 
left $B$ module. There are  vector space isomorphisms
$${\rm Hom}_B(V,\bar M(\lambda)^*)\cong {\rm Hom}_{\cal
M}({\cal C}v_{\lambda}, V^*)$$ 
and $${\rm Hom}_B(W, \bar M(\lambda)^{r*})\cong 
{\rm Hom}_{\cal M}({\cal C}v_{\lambda}^r, W^*).$$
\end{lemma}
 
\noindent {\bf Proof:} 
Given  $\psi\in {\rm Hom}_B(V,\bar M(\lambda)^*)$, define a
linear map $\tilde\psi$ from ${\cal C}v_{\lambda}$ to $ V^*$
by $\tilde\psi(v_{\lambda})(w)=\psi(w)(v_{\lambda})$. Note that 
$\psi(w)(mv_{\lambda})=(\psi(w)m)v_{\lambda}=\psi(wm)v_{\lambda}$
for all $m\in {\cal M}$.  
Hence
$\tilde\psi(mv_{\lambda})(w)=\tilde\psi(v_{\lambda})(wm)$ for $m\in {\cal M}$.  Thus 
$\tilde\psi\in {\rm Hom}_{\cal
M}({\cal C}v_{\lambda}, V^*)$.  

Using Lemma 2.1, we obtain a graded version of 
Theorem 2.2.    In particular, there is an isomorphism of vector 
spaces via the (graded) multiplication map:
$$\gra U\cong B\otimes \gra {\cal C}[T']\otimes N^-.$$
It follows that $\bar M(\lambda)=B v_{\lambda}.$  Hence $\psi(w)$ is
completely determined by its action on $v_{\lambda}$.  Thus the
map from $\psi$ to $\tilde\psi$ is one-to-one. The first isomorphism 
now follows from the fact that this map is 
clearly invertible.  A similar argument verifies the second isomorphism.
$\Box$

\medskip
Note that Lemmas 5.1, 5.4, and 5.5 hold  when we replace $B$ by
any subalgebra in ${\cal B}$.      In particular, these lemmas apply to $B'$.
 Let
$V_1$ denote the trivial one dimensional left
$B'$ module.  It follows that $V_1$ is annihilated by $B'_+$.  Let 
$V_1^r$ denote the trivial one-dimensional right $B$ module.   Then by
Lemma 5.5, $$\dim {\rm Hom}_{B'}(V_1,\hat M({\lambda}))=\dim {\rm
Hom}_{\cal M}({\cal C}v_{\lambda}, V_1)=1.$$ Similarly,
$\dim {\rm Hom}_B(V_1^r,\hat M({\lambda})^r)=1.$  In particular, 
the space of $B'$ invariants in $\hat M({\lambda})$ is one 
dimensional and the space of $B$ invariants in   $\hat
M({\lambda})^r$ is one dimensional. Let
$\zeta^r_{\lambda}$ be a nonzero vector in $(\hat M({\lambda})^r)^B$
and $\zeta_{\lambda}$ be a nonzero vector in $\hat M({\lambda})^{B'}$.

Let  $\hat N^-$
be the space consisting of
possibly infinite sums  of the form $\sum_{\gamma\leq 0}a_{\gamma}$
where $a_{\gamma}$ is a weight vector of weight $\gamma$ in $N^-$. 
Similarly, let  $\hat N^+$
be the space consisting of
possibly infinite sums  of the form $\sum_{\gamma\geq 0}a_{\gamma}$
where $a_{\gamma}$ is a weight vector of weight $\gamma$ in $N^+$. Note
that $\hat M({\lambda})^r=v_{\lambda}^r\hat N^+$ and $\hat
M({\lambda})=\hat N^- v_{\lambda}$. 

\begin{lemma}  There exists $b\in \hat N^-$ and $b^r\in \hat N^+$  such that
$\zeta_{\lambda}=b v_{\lambda}$ and 
$\zeta^r_{\lambda}=v_{\lambda}^rb^r$  for all  $\lambda\in 
\gh_{\Theta}^*$.  Moreover, both $b$ and $b^r$ have nonzero constant 
terms.
\end{lemma}

\noindent {\bf Proof:} The fact that there is a universal element $b$  which
satisfies  
$\zeta_{\lambda}=bv_{\lambda}$    for any choice of  $\lambda\in 
\gh_{\Theta}^*$   follows immediately from Lemma
5.4. Similarly, Lemma 5.4 ensures the existence of unique element $b^r$ satisfying 
$\zeta^r_{\lambda}=v_{\lambda}^rb^r$  for all  $\lambda\in 
\gh_{\Theta}^*$. We prove the last statement of the lemma.  Fix
$\lambda\in
\gh^*_{\Theta}$ and  write
$b=\sum_{\gamma\leq \beta}b_{\gamma}$ where $b_{\beta}\neq 0$.   If 
$x_jb_{\gamma}v_{\lambda}\neq 0$,
then $x_jb_{\gamma}v_{\lambda}$ has weight $\alpha_j+\gamma+\lambda$.  Hence 
$$x_jbv_{\lambda}=x_jb_{\beta}v_{\lambda}+{\rm \ terms\ of\ weight\ 
lower\ than\ 
}\alpha_j+\beta+\lambda.$$  Now if $\alpha_j\in \pi_{\Theta}$, then 
$x_jbv_{\lambda}=0$.  It follows that $x_jb_{\beta}v_{\lambda}=0$ for all $\alpha_j\in 
\pi_{\Theta}$.  On the other hand,  as in the proof of Lemma 2.1, if
$\alpha_i\notin\pi_{\Theta}$,  then  $B_+$ contains an element of the form $x_i+Y_i$, 
 where  $Y_i$ is a weight vector in 
$G^-T_{\Theta}$ of weight $\Theta(\alpha_i)$. In particular, $Y_ib_{\gamma}v_{\lambda}$ 
has weight strictly lower than $\gamma$.   So for 
$\alpha_i\notin\pi_{\Theta}$, we have 
$$0=(x_i+Y_i)bv_{\lambda}=x_ib_{\beta}v_{\lambda}+{\rm \ terms \ of \ weight \ 
lower \ 
than \ }\alpha_i+\beta+\lambda.$$  Hence $x_ib_{\beta}v_{\lambda}=0$ for all $i$ such 
that $1\leq i\leq n$.   Since the only highest weight vectors in 
$\bar M(\lambda)$ are scalar multiples of $v_{\lambda}$, it follows that 
$b_{\beta}$ is a nonzero scalar.   The same argument works for $b^r$ 
using  $y_it_i+\tilde\Theta(y_i)t_i$ instead of $x_i+Y_i$ and $y_jt_j$ 
instead of $x_j$.  $\Box$
\medskip

Note that we can consider elements of $\bar M(\lambda)^*\otimes 
\bar M(\lambda)$ as functions on $T$ where $(m^*\otimes m)(t)=m^*(tm)$ 
for all $m\in \bar M(\lambda)$, $m^*\in \bar M^*({\lambda})$, and $t\in T$.   
This gives rise to a linear map $\bar\Upsilon$ from $\bar M(\lambda)^*\otimes 
\bar M(\lambda)$ to ${\cal C}[P(\pi)]$ such that 
$(\bar\Upsilon(m^*\otimes m))(t)=(m^*\otimes m)(t)$. Recall the 
identification of $\hat M(\lambda)^r$ with $\bar M(\lambda)^*$. Consider  weight vectors
$m_{\gamma}^r\in  \bar M({\lambda})^r_{\gamma}$ and $m_{\gamma'}\in 
\bar M(\lambda)_{\gamma'}$.
Note that $m_{\gamma}^r(tm_{\gamma'})$ is zero if
$\gamma\neq\gamma'$.   It follows that the map
$\bar\Upsilon$ can be extended to a  linear map, which we also refer to as $\bar\Upsilon$,  from $\hat
M({\lambda})^r\otimes 
\hat M({\lambda})$ to the formal Laurent series ring
${\cal C}((z^{-\alpha}|\alpha\in \pi))$ such that (with the obvious interpretations)
$$\bar\Upsilon(\sum_{\gamma\leq \lambda}m_{\gamma}^r\otimes 
\sum_{\gamma\leq \lambda}m_{\gamma})=\sum_{\gamma\leq \lambda}\bar\Upsilon
(m_{\gamma}^r\otimes 
 m_{\gamma}).$$
Let ${\cal C}[[z^{-\tilde\alpha_i}|\alpha_i\in \pi^*]]$ denote  the subring of  ${\cal C}((z^{-\alpha}|\alpha\in \pi))$
consisting of elements of the form $\sum_{\tilde\gamma\geq 0}a_{\tilde\gamma}z^{-\tilde\gamma}$ where each
$\tilde\gamma\in Q^+(\Sigma)$ and  $a_{\tilde\gamma}\in
{\cal C}$.
 
 Note that $\gra U$ inherits a triangular decomposition from $U$.   
Let $\bar {\cal P}$ be the projection of $\gra U$ onto $\gra U^0$ 
using the direct sum decomposition $$\gra U=\gra U^0\oplus\gra 
(G^-_+U+ UU^+_+).$$  Now $\gra 
x_iy_jt_j=\delta_{ij}(q_i-q_i^{-1})^{-1}+\gra q_i^{-1}y_jt_jx_i$ for 
all $\alpha_i\notin\pi_{\Theta}$.   Hence $$\bar
{\cal P}(\gra U^+G^-)\subseteq {\cal C}[T_{\Theta}].$$
 
 \begin{lemma}  There exists  $p\in {\cal
C}[[z^{-\tilde\alpha_i}|\alpha_i\in \pi^*]]$
	 such that $\bar\Upsilon(\zeta_{\lambda}^r\otimes 
	 \zeta_{\lambda})=z^{\lambda}p$ 
	 for all  $\lambda\in \gh^*_{\Theta}$.  Moreover, $p$ has a nonzero 
	 constant term.
\end{lemma}

\noindent
{\bf Proof:}   Let
$b$ and
$b^r$ be as in Lemma 5.6.   We can write 
$b=\sum_{\gamma}b_{-\gamma}$ where each $b_{-\gamma}$ is a weight vector of 
weight $-\gamma$ in $N^-$.  Similarly, we can write 
$b^r=\sum_{\gamma}b^r_{\gamma}$ where each $b_{\gamma}^r$ is a weight 
vector of weight $\gamma$ in $N^+$. Let $b''_{\gamma}$ be the scalar such that 
$\bar {\cal P}(\gra b_{\gamma}^rb_{-\gamma})\in b''_{\gamma}+{\cal
C}[T_{\Theta}]_+$.  Note  that $s v_{\lambda}=v_{\lambda}^rs=0$ for all
$s\in {\cal  C}[T_{\Theta}]_+$ and $\lambda\in \gh^*_{\Theta}$.  
It follows that 
$$\eqalign{(v^r_{\lambda}b^r_{\gamma})(\tau(\beta)b_{-\gamma}v_{\lambda})=&(v^r_{\lambda})((\gra 
b^r_{\gamma}\tau(\beta)b_{-\gamma})v_{\lambda})\cr=&q^{(\beta,-\gamma+\lambda)}
v^r_{\lambda}(
(\gra b_{\gamma}^rb_{-\gamma})v_{\lambda})\cr=&q^{(\beta,-\gamma+\lambda)}
v^r_{\lambda}
(b''_{\gamma}v_{\lambda})\cr}$$ for all $\tau(\beta)\in T$ and $\lambda\in 
\gh^*_{\Theta}$.    Thus $$z^{-\lambda}\bar\Upsilon(b_{-\gamma}v_{\lambda}\otimes 
v^r_{\lambda}b^r_{\gamma})=b''_{\gamma}z^{-\gamma}$$ for each 
$\gamma$ and for all $\lambda\in 
\gh^*_{\Theta}$.
Therefore, by Lemma 5.6,  $z^{-\lambda}\bar\Upsilon(\zeta_{\lambda}^r\otimes 
	 \zeta_{\lambda})=z^{-\lambda'}\bar\Upsilon(\zeta_{\lambda'}^r\otimes 
	 \zeta_{\lambda'})$ for all $\lambda $ and $\lambda'$ in $\gh^*_{\Theta}$. 
	 This proves the first assertion. The second assertion    follows from 
	 the fact  that both $b$ and $b^r$ have nonzero constant terms (Lemma 
	 5.6). $\Box$
  	 \medskip

We can give ${\cal C}(Q_{\Sigma}){\cal A}$ a filtration by setting 
$\deg f=0$ for all $f\in {\cal C}(Q_{\Sigma})$ and setting the degree of
an element $t\in {\cal A}$ equal to its degree in the ${\cal F}$
filtration of $U$. Note that ${\cal C}(Q_{\Sigma}){\cal A}$ is 
isomorphic to its associated graded ring under this filtration.   
Given a homogeneous element $g\tau(\beta)$ where $g\in {\cal 
C}(Q_{\Sigma})$ and $\tau(\beta)\in {\cal A}$, we write $\gra 
g\tau(\beta)$ as just $g\tau(\beta)$.

Note  that ${\cal F}$ extends in an obvious way to a filtration 
on $\check U$. Similarly, the above filtration extends to ${\cal 
C}(Q_{\Sigma})\check {\cal A}$.
 Now suppose
$z\in Z(\check U)$.   Then $\gra  z$ is in the center of $\gra \check U$.   
As explained in Section 4, given $\mu\in P^+(\pi)$, there exists a  central element 
$c_{\mu}$ in $(\ad U)\tau(-2\mu)$ such that $c_{\mu}\in 
\tau(-2\mu)+(\ad U_+)\tau(-2\mu)$. Furthermore, by (4.8) (which 
applies in general and not just to the rank one cases), there exists a
nonzero scalar multiple $c_{\mu}'$ of $c_{\mu}$ such that 
$$c_{\mu}'\in \tau(-2\mu)+G^-_+U_{\geq}U^+_+\tau(-2\mu).\leqno{(5.4)}$$
  Hence by Lemma 3.1,
there exists
$p_{\mu}\in {\cal C}(Q_{\Sigma})$ such that 
$\gra c_{\mu}'\tau(\beta)$ is an element of $$ B_+\gra 
(G^-U^+\tau(-2\mu+\beta))+
\gra (G^-U^+\tau(-2\mu+\beta))B'_++\tau(-2\tilde\mu)(p_{\mu}\cdot\tau(\beta)).$$ 
 It follows that  $$\gra({\cal
X}( c_{\mu}'))= \tau(-2\tilde\mu)p_{\mu}$$ and
$$\eqalign{( \zeta^r_{\lambda}\otimes \zeta_{\lambda})(\gra 
c_{\mu}\tau(\beta))&=\bar\Upsilon(\zeta^r_{\lambda}\otimes 
\zeta_{\lambda})(
(\gra {\cal X}( c_{\mu}'))\cdot\tau(\beta))\cr&=\bar\Upsilon(\zeta^r_{\lambda}\otimes
\zeta_{\lambda}) ((\tau(-2\tilde\mu)p_{\mu})\cdot\tau(\beta))\cr}$$ for all
$\tau(\beta)\in {\cal A}$.

\begin{theorem} There exists   $p\in {\cal
C}[[z^{-\tilde\alpha_i}|\alpha_i\in \pi^*]]$ such that
$\bar\Upsilon(\zeta^r_{\lambda}\otimes
\zeta_{\lambda})=z^{\lambda}p$ for all $\lambda\in \gh^*_{\Theta}$ and 
$\gra( {\cal X} (  c_{\mu}'))=p^{-1}\tau(-2\tilde\mu)p$ for all $\mu\in
P^+(\pi)$.    
\end{theorem}  
 
\noindent{\bf Proof:} 
The first assertion is simply Lemma 5.7. 
By (5.4), it follows that $(\gra
c_{\mu}')v_{\lambda}=q^{(\lambda,\mu)}v_{\lambda}$.  Furthermore, 
$\lambda\in \gh^*_{\Theta}$ ensures that 
$(\lambda,\mu)=(\lambda,\tilde\mu)$.  Since $(\gra c_{\mu}')$ is
central in $\gra U$, it follows that $\gra c_{\mu}'$ acts on $\bar
M(\lambda)$ as multiplication by the scalar $q^{(\lambda,\tilde\mu)}$.  
A similar argument yields that $\gra c_{\mu}'$ acts on $\bar M(\lambda)^r$
as multiplication by the same scalar $q^{(\lambda,\tilde\mu)}$. Hence
$\gra c_{\mu}'$ acts on elements of $\hat M(\lambda)^r\otimes \hat
M(\lambda)$ as multiplication by the scalar $q^{(\lambda,\tilde\mu)}$.  
It follows that
$$\eqalign{\bar\Upsilon(\zeta^r_{\lambda}\otimes \zeta_{\lambda})*
p^{-1}\tau(-2\tilde\mu)p=& z^{\lambda}p* 
p^{-1}\tau(-2\tilde\mu)p\cr=& 
q^{(\lambda,\tilde\mu)}z^{\lambda}p\cr=&z^{\lambda} p*(\gra  
{\cal X}(c_{\mu}')).\cr}$$
 In particular, $$\gra ({\cal X}(
c_{\mu}'))-p^{-1}\tau(-2\tilde\mu)p\leqno(5.5)$$  acts as zero on
$z^{\lambda}p$ for all $\lambda\in \gh^*_{\Theta}$. 
But $\gra ({\cal X}(
c_{\mu}'))-p^{-1}\tau(-2\tilde\mu)p$ is an element of ${\cal
C}[[z^{-\tilde\alpha_i}|\alpha_i\in \pi^*]]\tau(-2\tilde\mu)$.  The only element of this set which acts as
zero on $z^{\lambda}p$ is zero.     This forces   
the expression in (5.5) to be identically equal to zero.
$\Box$

\section{Computing Graded Radial Components}
	 
In this section, we compute the graded image of  radial components using 
information about rank one quantum symmetric pairs from Section 4. 
In particular, for each  $i$ such that $\alpha_i\in \pi^*$, we 
associate a semisimple Lie subalgebra $\gg_i$ of $\gg$ such that 
$\gg_i,\gg_i^{\theta}$ is an irreducible symmetric pair with rank one 
restricted root system as follows.  Recall that $\omega_j$ denotes 
the fundamental weight corresponding to the root $\alpha_j\in \pi$.
For all $i$ such that $\alpha_i\in \pi^*$, set
$\pi_i=\{\alpha_j|(\omega_j,\Theta(-\alpha_i))\neq 0 {\rm\ or \ } 
(\omega_j,\Theta(-\alpha_{p(i)}))\neq 0 \}$. Let 
 $\gg_{i}\subseteq \gg$ be the semisimple  Lie subalgebra 
generated by the root vectors $e_j$ and $f_j$ with $\alpha_j\in 
\pi_{i}$.  It follows that $\pi_i$ is the set of simple roots 
associated to $\gg_i$.  Moreover, the choice of $\pi_i$ ensures that $\theta$ 
restricts to an involution of $\gg_i$ which we also refer to as 
$\theta$.   Set $\Sigma_i=\{\pm\tilde\alpha_i\}$ and note that 
$\Sigma_i$ is precisely the set of restricted roots associated to the 
symmetric pair 
$\gg_i,\gg_i^{\theta}$. 

Let  $\Delta_i$ denote the root system associated to $\gg_i$ and set  
  $U_{i}=U_q(\gg_{i})$. Set $\Delta_i^+=\Delta^+\cap \Delta_i$. Note that $B\cap U_i$ can be though of as a 
  (standard) quantum analog of $U(\gg_i^{\theta_i})$ inside of $U_i$.   In 
  particular, results in the previous sections of this paper apply to 
  the quantum symmetric pair $U_i, B\cap U_i$.  A similar comment can
  be made in reference to the subalgebra $B'$ of $U$.

For most standard subsets of $U$, we use the subscript $i$ to 
denote the intersection of this subset with $U_i$. For example, we write 
$U_i^+$ for $U^+\cap U_i$.  The  
exception to this rule is    $B\cap U_i$ since $B_i$ has already been 
defined as something different in (1.1).

Set 
$\gh^*_{\Theta i}=\{\lambda\in 
\gh^*|(\lambda,\eta)=0$ for 
all $\eta\in Q(\pi_{i})$ such that $\Theta(\eta)=\eta\}$.
  Note that $\gh^*_{\Theta}$ is 
a subset of $\gh^*_{\Theta i}$.  Given $\lambda\in \gh^*_{\Theta i}$, 
let $v_{\lambda}$ be a 
$T$ weight vector of weight $\lambda$ and give ${\cal C}v_{\lambda}$ 
the structure of a trivial ${\cal M}_iU_i^+$ module.  Write 
$\bar M_i(\lambda)$ for the (left) $\gra U_iT$ module 
induced from the $\gra T{\cal M}_iU_i^+$ module ${\cal C}v_{\lambda}$.

\begin{lemma}
	Fix $\lambda\in \gh^*_{\Theta}$ and let $w$ be a weight 
	vector in $\bar M(\lambda)$ of weight $\gamma$. Assume further that
	$sw=0$ for all $s\in \gra ({\cal M}_iT_{\Theta i}U_i^+)_+$.  
	Then $\gamma\in \gh^*_{\Theta i}$ 
	and the map $ uv_{\gamma}\mapsto uw$ is a $\gra U_{i}T$ 
	module isomorphism from 
	$\gra U_{i}Tv_{\gamma}$ onto $\gra U_{i}Tw$.
\end{lemma}

\noindent
{\bf Proof:}   Since $sw=0$ for all $s\in {\cal C}[T_{\Theta i}]_+$, it 
follows that $\tau(\eta)w=q^{(\eta,\gamma)}w=w$ for all $\tau(\eta)\in T_{\Theta i}$. 
Now $\tau(\eta)\in T_{\Theta i}$ if and only if 
$\eta\in Q(\pi_{i})$ and 
$\Theta(\eta)=\eta$.
Hence $\gamma\in \gh^*_{\Theta i}$.

Set $I=UU^+_++\sum_{\tau(\beta)\in T}
U(\tau(\beta)-q^{(\gamma,\beta)})$.  Note that 
$\bar M_i({\gamma})$ is isomorphic to the left 
$\gra U_{i}T$ module $(\gra U_{i}T)/(\gra (I\cap  
U_{i}T))$.
  Since $w$ is annihilated 
by $I$, it follows that the map $uv_{\gamma}\mapsto 
	uw$ is a $\gra U_{i}T$ module map. 
	
	By our assumptions on $w$, we have that $(\gra U_{i}T) 
	w=N^-_iw$.  Now $\bar M(\lambda)$ is 
	a free $N^-$ module.  Hence  the subspace
	$(\gra U_{i}T)w$ is a cyclic free  
	$N^-_i$ module.   The lemma now follows from the 
	fact that 
	$\bar M_i({\gamma})$ 
	is also a cyclic free  $N^-_i$ module. $\Box$

\medskip
Let $\hat  M_i(\lambda)$ be the subspace of $\hat M(\lambda)$ 
consisting of  possibly infinite sums of weight vectors 
 in $\bar M_i(\lambda)$. Similarly, let $\hat N^-_i$ be the subspace
of 
 $\hat N^-$ consisting of  possibly infinite sums of weight vectors 
 in $N^-_i$. Given $\lambda\in \gh^*_{\Theta i}$, let 
$\zeta_{\lambda i}$ denote the $B'\cap U_i$ 
invariant vector of $\hat  M_i(\lambda)$.
By Lemma 5.6, there is an element  $b_{i}$  in $\hat N^-_i$ 
such that $\zeta_{\lambda i}=b_{i}v_{\lambda}$.
	 
\begin{lemma} 
	There exists $w=\sum_{\gamma}w_{\gamma}\in \hat N^-$ such that  each 
	$w_{\gamma}$ is a weight vector of weight $\gamma$, 
	$sw=0$ for all $s\in \gra ({\cal M}_iT_{\Theta i}U_i^+)_+$, and 
	$b_{i}wv_{\lambda}=\zeta_{\lambda}$
	for all $\lambda\in \gh^*_{\Theta}$.
	\end{lemma}
	
\noindent
{\bf Proof:}  Fix $\lambda\in \gh^*_{\Theta}$.  Let 
$m=\sum_{\gamma}m_{\gamma}$ be a vector in $\hat  M(\lambda)$ such 
that $(B'\cap U_i)_+m=0$. It follows that $sm_{\gamma}=0$ for all 
$s\in {\cal C}[T_{\Theta i}]_+$ and all 
$\gamma$.   Thus $m_{\gamma}$ nonzero implies that $\gamma\in 
\gh^*_{\Theta i}$.

Let $\beta$ be the highest weight 
such that $m_{\beta}\neq 0$. The same argument as in Lemma 5.6 shows 
that $x_jm_{\beta}=0$ for all $\alpha_j\in \pi_{i}$. By 
the previous paragraph, 
$\beta$ is a weight in $\gh^*_{\Theta i}$.  

Note that each $y_jt_j$ is a highest weight vector for the action of 
$\ad {\cal M}^+$.   By  [JL1,Section 4], 
$\ad y_k$ acts ad nilpotently on  $y_jt_j$ whenever $k\neq j$.
 Hence   [JL1, Theorem 5.9] ensures that  each $y_jt_j$ such that 
$\alpha_j\notin\pi_{\Theta}$ generates a locally  finite
$\ad {\cal M}$ module.  By  the definition of $N^-$, it 
follows that $N^-$ is a 
locally finite module with respect to the action of $\ad {\cal M}$. 
Hence $m_{\beta}$ generates a finite 
dimensional $\ad {\cal M}_i$ module. The fact that
$(\beta,\alpha_j)=0$ for all $\alpha_j\in\pi_{i}\cap 
\pi_{\Theta}$ further implies that $m_{\beta}$ generates a 
one-dimensional trivial $\ad {\cal M}_i$ module.
Thus by Lemma 6.1, $b_{i}m_{\beta}$ is a 
$B'\cap U_i$ invariant vector.  Moreover, rescaling if 
necessary, we may assume by Lemma 5.6 that 
$b_{i}m_{\beta}=m_{\beta}+$ terms of weight lower than 
$\beta$.  Now $m'=m-b_{i}m_{\beta}$ is  a $B'\cap U_i$ 
invariant vector. Moreover, when $m'$ is written as a sum of weight 
vectors, the highest weight of a nonzero summand is strictly less 
than $\beta$.   By induction, we can find a sequence of weight vectors 
$\{w_{\gamma}\}$ such that  each $w_{\gamma}$ generates a trivial 
$U^+_i{\cal 
M}_iT_{\Theta i}$ module and 
$$m-b_{i}\sum_{0\geq\gamma\geq \gamma'}w_{\gamma}\in 
\sum_{\beta<\gamma'}N^-_{\beta}.$$  Thus by the definition of $\hat 
N^-$, we obtain $m=b_{i}\sum_{\gamma}w_{\gamma}$.$\Box$

\medskip

 Let $G_{\pi\setminus\pi_i}^-$ denote the subalgebra of $G^-$ generated 
by $(\ad U_{i}^-){\cal 
C}[y_it_i|\alpha_i\notin\pi_{i}]$.  By [L3, Section 6],
we see that multiplication induces an isomorphism of vector spaces
$$G^-\cong G_{i}^-\otimes G_{\pi\setminus\pi_i}^-.\leqno{(6.1)}$$
Furthermore, $G_{\pi\setminus\pi_i}^-$ is generated by elements of 
the form $(\ad y)y_jt_j$ where $\alpha_j\notin\pi_{i}$,
$y\in G_i^-$, and the weight of $(\ad y)y_jt_j$
is a root  in $\Delta$. It follows that the weight of $(\ad 
y)y_jt_j$ cannot be an element of $\Delta_i$.  
Hence the weights of vectors in $G^-_{\pi\setminus\pi_i}$ are elements of
 $$\sum_{\gamma\in \Delta^+\setminus\Delta^+_i}{\bf N}(-\gamma).$$

  Suppose that $\beta\in \Delta$
and
$\tilde\beta=\tilde\alpha_i$.    Note that
$2\tilde\alpha_i\in
\sum_{\alpha\in
\pi_i}{\bf N}\alpha_i$. Hence $\beta$ must be a positive root. So both
$\beta$ and
$-\Theta(\beta)$ are elements of $\sum_{\alpha\in \pi}{\bf N}\alpha$. 
Hence
$2\tilde\beta\in
\sum_{\alpha\in
\pi_i}{\bf N}\alpha$ forces both $\beta$ and $-\Theta(\beta)$ to be
elements of $\sum_{\alpha\in
\pi_i}{\bf N}\alpha$.   In particular  
$$\{\beta\in{\Delta}^+|\tilde\beta=\tilde\alpha_i\}
\subset \Delta_i^+.\leqno{(6.2)}
$$

By (6.2), if $\gamma\in  \Delta^+\setminus\Delta^+_i$, then 
 $\tilde\gamma\notin\Sigma_{i}$.    In particular, if $\beta$ 
 is a weight of an element of $G^-_{\pi\setminus\pi_i}$, then 
 $$\tilde\beta\in \sum_{\tilde\gamma\in 
	  \Sigma^+\setminus\{\tilde\alpha_i\}}{\bf N}(-\tilde\gamma).$$	
	 
\begin{lemma} Let $\lambda\in \gh^*_{\Theta}$. Suppose $w$ is a weight vector of weight $\beta$ in $N^-$
	  such that $wv_{\lambda}$ generates a trivial 
	  $\gra(U_i^+{\cal M}_iT_{\Theta i})$ module
	  and $\beta\in Q(\Sigma)$.
	  Then the weight of $w$ is contained in $\sum_{\tilde\gamma\in 
	  \Sigma^+\setminus\{\tilde\alpha_i\}}{\bf N}(-\tilde\gamma).$
	  \end{lemma}
	  
\noindent
{\bf Proof:}  Write $w=\sum_j w_{1j}  w_{2j}$ where $w_{1j}\in 
 G_i$ and $w_{2j}\in G_{\pi\setminus \pi_i}^-$.  For each $j$,  set $\gamma_{1j}$ equal to 
 the weight of     $w_{1j}$ and $\gamma_{2j}$ equal to the weight of  $w_{2j}$. 
 We may further assume that  $\{w_{2j}\}_j$ is a 
linearly independent set. Choose $\beta'$ minimal in the set 
$\{\gamma_{2j}\}_j$ using the standard partial ordering on $Q(\pi)$.  
By the discussion preceding the lemma, we have that $$\tilde\beta'\in 
\sum_{\tilde\gamma\in 
	  \Sigma^+\setminus\{\tilde\alpha_i\}}{\bf N}(-\tilde\gamma).$$  Thus it is 
	  sufficient to show that $\beta=\beta'$.
	  
Consider $\alpha_k\in \pi_i$. In particular, $x_kwv_{\lambda}=0$.  By 
Lemma 5.2, $\bar M(\lambda)$ is isomorphic to the submodule $N^-$ of 
$G^-$ as $U^+$ 
modules.   Hence $x_k*w=\gra ((\ad x_k)w)=0$.	On 
the other hand,
$$\gra ((\ad x_k)w)\in \sum_{\gamma_{2j}=\beta'}\gra ((\ad 
x_k)w_{1j})w_{2j}+ 
G_i^-\sum_{\gamma\not<\beta'}(G_{\pi\setminus\pi_i}^-)_{\gamma}.$$ 
Hence, by (6.1),  $\gra (\ad x_k)w_{1j}=0$ for all $j$ such that 
$\gamma_{2j}=\beta'$. Fix $j$ such that $\gamma_{2j}=\beta'$. It follows that
$w_{1j}$ is a highest weight vector with 
respect to the action of  $U_i^+$  on
$G_i^-$.   By [JL2, Lemma 4.7(i)], $w_{1j}$ is a scalar.  Thus the weight of $w$ agrees 
with the weight of $w_{2j}$ which is just $\beta'$.   $\Box$

\medskip  
Given $\lambda\in \gh^*_{\Theta i}$, 
let $v_{\lambda}^r$ be a right $T$ weight vector of 
weight $\lambda$ and give ${\cal C}v_{\lambda}^r$ the structure of a 
one dimensional trivial $\gra ({\cal M}_iG_i^-)$ module. 
Define $\bar M_i(\lambda)^r$ to be the right $\gra U_iT$ 
module induced from 
the one dimensional $\gra ({\cal M}_iG_i^-T)$ module ${\cal C}v_{\lambda}^r$. 
Note that versions of Lemmas 6.1-6.3 hold for  the right $\gra U_i$ 
modules
$\bar M_i(\lambda)^r$.   In particular, let $b^r$ be chosen as in Lemma 
5.6.  Let $b_i^r$ be also chosen as in Lemma 5.6 for the modules 
$\bar M_i(\lambda)^r$. 
As in Lemma 6.2, there exists $w^r\in \hat N^+$ such that $w^rb_{i}^r=b^r$.
Moreover, $w^r=\sum_{\gamma}w_{\gamma}^r$ where each $w_{\gamma}$ is
annihilated by elements in $\gra (G^-_i{\cal
M}_iT_{\Theta i})_+$.   Furthermore, $w^r_{\gamma}\neq 0$ implies that $\gamma$
is an element of $\sum_{\tilde\beta\in
\Sigma^+\setminus\{\tilde\alpha_i\}}{\bf N}\tilde\beta.$

Recall the graded version $\bar {\cal P}$ of the  Harish-Chandra map
defined in the last section.  The vector space decomposition (3.5) extends  in the obvious way to the 
corresponding graded algebras.   Let $\bar{\cal P_{\cal A}}$ denote 
composition of $\bar{\cal P}$ with projection onto $\gra {\cal 
C}[\check{\cal A}]$ using a graded version of (3.5).

Given
$\lambda\in \gh^*_{\Theta}$, assume  that
$v_{\lambda}^r$ has been chosen so that 
$v_{\lambda}^r( v_{\lambda})=1.$  It follows that 
$v_{\lambda}^r(\tau(\beta)v_{\lambda})=q^{(\lambda,\beta)}$ for all 
$\tau(\beta)\in T$. In particular,
$v_{\lambda}^r(tv_{\lambda})=1$ for all $t\in T_{\Theta}$.   Hence 
$v_{\lambda}^rc\tau(\beta)dv_{\lambda}=q^{(\lambda-\gamma,\beta)}(\bar{\cal P}_{\cal A}
(\gra cd))$
for all $c\in U^+_{\gamma}$, $\tau(\beta)\in T$, and $d\in G^-_{-\gamma}$.
Thus $$\bar\Upsilon (v_{\lambda}^rc\otimes 
dv_{\lambda})=z^{\lambda-\gamma}(\bar{\cal P}_{\cal A}(\gra cd))$$ for 
all $c\in U^+_{\gamma}$  and $d\in G^-_{-\gamma}$.

Choose $p_i\in {\cal 
C}[[z^{-\tilde\alpha_i}|\alpha_i\in \pi^*]]$ as in Lemma 5.7  such that 
$\Upsilon(\zeta_{\lambda i}^r\otimes \zeta_{\lambda i})=z^{\lambda}p_i$.
Let ${\cal C}[[z^{-\tilde\beta}|\tilde\beta\in 
\Sigma^+\setminus{\tilde\alpha_i}]]$ denote the subring of ${\cal 
C}[[z^{-\tilde\alpha_i}|\alpha_i\in \pi^*]]$  consisting of possibly infinite linear combinations of 
the $z^{-\nu}$ for $\nu\in \sum_{\beta\in 
\Sigma^+\setminus\{\tilde\alpha_i\}} {\bf
N}\beta$.

\begin{lemma} There exists $k_i\in {\cal C}[[z^{-\tilde\beta}|\tilde\beta\in 
\Sigma^+\setminus\{\tilde\alpha_i\}]]$ such that   $p=p_{i}k_i$. Furthermore, 
$k_i$ has a nonzero constant term.
\end{lemma}

\noindent
{\bf Proof:} By Lemma 5.7, $p$ has a nonzero constant term.   Hence 
if we can write $p=p_ik_i$, then both $p_i$ and $k_i$ have nonzero 
constant terms.  Thus the second assertion follows from the first.

Fix $\lambda\in \gh^*_{\Theta}$. Using Lemma 5.6, choose $b\in \hat N^-$ such 
that $\zeta_{\lambda}=bv_{\lambda}$  and 
 $b_i\in \hat N_i^-$ such that  $\zeta_{\lambda i}=b_iv_{\lambda}$ 
Write $$b_i=\sum_{\delta}b_{i\delta}{\rm
\  and \ }b_i^r=\sum_{\delta}b_{i\delta}^r$$ where each $b_{i\delta}\in 
(N^-_i)_{-\delta}$ and $b_{i\delta}^r\in (N_i^+)_{\delta}$.
Note that
$$z^{\lambda}p_i=\sum_{\delta}\bar\Upsilon(v_{\lambda}^rb_{i\delta}^r\otimes 
b_{i\delta}v_{\lambda})=\sum_{\delta}z^{\lambda-\delta}(\bar{\cal 
P}_{\cal A}
(\gra b^r_{i\delta}b_{i\delta})).$$

 Let   
 $w=\sum_{\gamma}w_{\gamma}$ satisfy the conditions of   Lemma 6.2. 
 By Lemma 6.3, each $\gamma$ with $w_{\gamma}$ nonzero satisfies 
 $-\gamma\in \sum_{\beta\in 
\Sigma^+\setminus\{\tilde\alpha_i\}} {\bf
N}\beta$. Choose $\gamma$ so 
that $w_{\gamma}\neq 0$.
Note  that $ 
U^+_{i+}w_{\gamma}v_{\lambda}=0$.
Using the identification  of $\bar M(\lambda)$ with $N^-$ as  $U^+$ 
modules (Lemma 5.2),
we obtain
$\gra (\ad  x_j) w_{\gamma}=0$  for all  $\alpha_j\in
\pi_{i}$. Hence $\gra x_jw_{\gamma}=\gra
t_jw_{\gamma}t_j^{-1}x_j$ for all $x_j\in U_i$. It
follows that
$\bar{\cal P}_{\cal A}(\gra (UU^+_{i+})w_{\gamma})=0$.   Similarly $\bar{\cal 
P}_{\cal A}
(w_{\gamma}^r(\gra G^-_{i+}U))=0$. To make the next computation 
easier to read,
 we shorten 
$\bar{\cal P}_{\cal A}(\gra u)$ to $\bar{\cal P}_{\cal A}( u)$
for $u\in U$. Since $b_{i\delta}\in  (N^-_i)_{-\delta}$
and 
$b_{i\delta}^r\in  (N^+_i)_{\delta}$, it follows that 
$$\eqalign{\sum_{\delta}\bar{\cal P}_{\cal A}
(w_{\gamma}^rb^r_{i\delta}\tau(\beta)
b_{i\delta}w_{\gamma})=&\sum_{\delta}\bar{\cal 
P}_{\cal A}(w_{\gamma}^r\bar{\cal P}_{\cal A}(b^r_{i\delta}\tau(\beta)
b_{i\delta})w_{\gamma})\cr
=&\sum_{\delta}q^{(-\delta,\beta)}\bar{\cal P}_{\cal A}
(b^r_{i\delta}
b_{i\delta})\bar{\cal P}_{\cal 
A}(w_{\gamma}^r\tau(\beta)w_{\gamma})\cr
=&\bar{\cal P}_{\cal A}(w_{\gamma}^r\tau(\beta)w_{\gamma})p_i\cdot 
\tau(\beta).\cr}$$ 
Now $\bar{\cal P}_{\cal A} (w_{\gamma}^r\tau(\beta)w_{\gamma}
)$ is equal to 
$a_{\gamma}z^{\gamma}\cdot \tau(\beta)$ for some scalar $a_{\gamma}$ 
independent of $\beta$.
Thus $$\bar\Upsilon(v_{\lambda}^rb^r\otimes bv_{\lambda})=
z^{\lambda}(\sum_{\gamma}a_{\gamma}z^{\gamma})p_{i}.$$   
  The lemma now follows from the fact that  $-\gamma$ is in the ${\bf N}$
span of the set $\Sigma^+\setminus\{\tilde\alpha_i\}$ whenever 
$a_{\gamma}\neq 0$.
$\Box$
\medskip

Let $\tilde s_i$ be the reflection in $W_{\Theta}$ corresponding to the restricted 
root $\tilde\alpha_i$.  Recall that $\tilde s_i$ restricts to a permutation on the set 
$\Sigma^+\setminus\{\tilde\alpha_i\}$.   Hence $\tilde s_i$ induces a linear 
map on ${\cal
C}[[z^{-\tilde\beta}|\tilde\beta\in 
\Sigma^+\setminus\{\tilde\alpha_i\}]]$ defined by  
$$\tilde s_i
\sum_{\tilde\gamma}z^{-\tilde\gamma}=\sum_{\tilde\gamma}z^{\tilde 
s_i(-\tilde\gamma)}.$$

\begin{lemma} Choose $k_i\in {\cal C}[[z^{-\tilde\beta}|\tilde\beta\in 
\Sigma^+\setminus\{\tilde\alpha_i\}]]$ such that $p=p_ik_i$.  Then $\tilde s_ik_i=k_i$. 
\end{lemma}

\noindent{\bf Proof:} Recall that $\omega_j$ denotes the fundamental 
weight corresponding to $\alpha_j\in \pi$. Given $\alpha_j$ and 
$\alpha_k$ in $ \pi^*$, 
we know
that $(\omega_j,\tilde\alpha_k)$ is a 
nonzero scalar multiple of $\delta_{jk}$. Set 
$\nu_i=(\sum_{\{r|\ \alpha_r\in \pi^*\}}\omega_r)-\omega_i$.  It follows that
$(\tilde\alpha_j,\nu_i)\neq 0$ for all $j\neq i$ and 
$(\tilde\alpha_i,\nu_i)=0$.  Hence $p_i$ commutes with 
$\tau(-2\nu_i)$. More generally,   $$z^{-\tilde\beta}\tau(-2\nu_i)
=q^{(-2\nu_i,-\tilde\beta)}\tau(-2\nu_i)z^{-\tilde\beta}$$ for all 
$\tilde\beta\in Q(\Sigma)$. Note that if
$\tilde\beta\in \Sigma^+$ is not a  scalar multiple of $\tilde\alpha_i$, then 
$(\tilde\beta,\nu_i)\neq 0$. Hence if $k\in {\cal 
C}[[z^{-\tilde\beta}|\tilde\beta\in \Sigma^+\setminus\{\tilde\alpha_i\}]]$ and $k$ 
commutes with $\tau(-2\nu_i)$ then $k$ is a scalar.  
Now     
by Corollary 3.3, Theorem 3.4 and Theorem 5.8,
$${\cal X}(c_{\nu_i})=\sum_{w\in 
W_{\Theta}}w(k_i^{-1}\tau(-2\tilde\nu_i)k_i)+\sum_{\{\mu\in 
P^+(2\Sigma)|\mu<\tilde\nu_i\}}\sum_{w\in 
W_{\Theta}}w(f_{\mu})\tau(-2w\mu)$$
up to a nonzero scalar for some $f_{\mu}\in {\cal C}(Q_{\Sigma})$.
Since ${\cal X}(c_{\nu_i})$ is $\tilde s_i$ invariant
and $\tilde s_i\tau(-2\nu_i)=\tau(-2\nu_i)$, it follows that $\tilde 
s_i(k_i\tau(-2\nu_i)k_i^{-1})=k_i\tau(-2\nu_i)k_i^{-1}$.  Hence 
$\tilde s_ik_i=k_ik$ where $k$ commutes with $\tau(-2\nu_i)$. 
Thus 
$k$ must be a scalar. Since $k_i$ has 
a nonzero constant term, and $\tilde s_i$ fixes constants, it follows 
that $k=1$. 
$\Box$
\medskip

Let $a$ and $x$ be an indeterminates
and define $$(x;a)_{\infty}=\prod_{i=0}^{\infty}(1-xa^i).\leqno(6.3)$$ Set
$z_i=z^{2\tilde\alpha_i}$.  Let $\rho_i$ denote the half sum of the
positive roots in $\Delta_i$.  Note that $\rho$ is just the sum of
the fundamental weights corresponding to the simple roots in $\pi$. 
A similar statement can be made concerning $\rho_i$ with respect to
$\pi_i$.  Hence $(\rho,\beta)=(\rho_i,\beta)$ for all $\beta\in
Q(\pi_i)$.  

For each $\tilde\alpha\in\Sigma^+$, we let ${\cal C}[[z^{-\tilde\alpha}]]$ denote the subring of ${\cal
C}[[z^{-\tilde\alpha_i}|\alpha_i\in \pi^*]]$ consisting of elements of the form $\sum_{m\geq 0}a_{m}z^{-m\tilde\alpha}$ 
where the $a_m$ are scalars. 
Recall that $z_i=z^{2\tilde\alpha_i}$. Using
the rank one computations found in Section 4, we  determine
$p_i$.

\begin{lemma}Given
$\alpha_i\in
\pi^*$,
we
have
$$p_{i}={{(g_iz_{i}^{-1};a_{i})_{\infty}}\over{(
z_{i}^{-1}; a_{i})_{\infty}} }$$  where $a_i
=q^{(2\tilde\alpha_i,\tilde\alpha_i)}$ and
$g_i=q_i^{2(\rho,\tilde\alpha_i)}$. 
\end{lemma}

\noindent 
{\bf Proof:} 
Let $\mu_i$ be chosen to satisfy the conditions of (4.7)
with respect to $\gg_i$.  By Theorem 4.7 and Theorem 5.8 
$$p_i^{-1}\tilde t_i^{-1}p_i=\gra ({\cal
X}(c_{\mu_i}'))=\tilde
t_i^{-1}(1-g_iz_i^{-1} )
(1-z_i^{-1})^{-1}. $$  
Note that $z_i^{-1}\tilde t_i^{-1}=a_i\tilde t_i^{-1}z^{-1}_i$.
Hence, $$\eqalign{{{(z_{i}^{-1};a_{i})_{\infty}}\over{(g_{i}
z_{i}^{-1}; a_{i})_{\infty}} }\tilde
t_i^{-1}{{(g_{i}
z_{i}^{-1}; a_{i})_{\infty}}\over {(z_{i}^{-1};a_{i})_{\infty}}}
=&\tilde t_i^{-1}\prod_{j=0}^{\infty}{{(1-z_i^{-1}a_i^{j+1})}\over
{(1-g_iz_i^{-1}a_i^{j+1})}}
{{(1-g_iz_i^{-1}a_i^j)}\over{(1-z_i^{-1}a_i^j)}}
\cr=&\tilde t_i^{-1}(1-g_iz_i^{-1} )
(1-z_i^{-1})^{-1}.\cr}$$
The lemma now follows from the fact that  the only elements in ${\cal
C}[[z^{-\tilde\alpha_i}]]$ which commute with $\tilde t_i$ are the
scalars. $\Box$
\medskip

Consider $\alpha\in \Delta$.  Since $\Sigma$ is reduced, $(\alpha,-\Theta(\alpha))$ equals $0$ if 
$-\alpha\neq \Theta(\alpha)$ and equals $(\alpha,\alpha)$ otherwise.
Hence it is straightforward to check that $\tilde\alpha=\tilde\beta$ 
implies $(\alpha,\alpha)=(\beta,\beta)$ for all  $\beta$ 
in $\Delta$. In particular, the length of a root whose restriction equals 
$\tilde\alpha$ is just a function of $\tilde\alpha$ and does not 
depend on the choice of root.
Set
$$mult({\tilde\alpha})=|\{\beta\in{\Delta}|\tilde\beta=\tilde\alpha\}|.$$
We recall two well known facts about the $mult$ function.
First, $mult(\tilde\alpha_i)=
2(\rho,\tilde\alpha_i)/(\tilde\alpha_i,\tilde\alpha_i)$ for all 
$\alpha_i\in \pi^*$.  Moreover,
$mult(\tilde\alpha)=mult(w\tilde\alpha)$  for all  $ w\in W_{\Theta}$
and $\alpha\in\Delta$.

Given $\tilde\alpha\in \Sigma^+$, set $a_{\tilde\alpha}=q^{(2\tilde\alpha,\tilde\alpha)}$ and
$g_{\tilde\alpha}=q^{mult(\tilde\alpha)(\tilde\alpha,\tilde\alpha)(\alpha,\alpha)/2}.$ 
Set 
$$p_{\tilde\alpha}={{(g_{\tilde\alpha}z^{-2\tilde\alpha};a_{\tilde\alpha})_{\infty}}
\over{(
z^{-2\tilde\alpha};
a_{\tilde\alpha})_{\infty}}
}$$ for $\tilde\alpha\in \Sigma$.  Note that
$p_i=p_{\tilde\alpha_i}$.   Furthermore, by the previous paragraph,
$wp_{\tilde\alpha}=p_{\tilde\beta}$ whenever $w\in W_{\Theta}$ 
satisfies $w\tilde\alpha=\tilde\beta$.
Now $p_{\tilde\alpha}$ is clearly an element of ${\bf
C}((q))[[z^{-\tilde\alpha}]]$ where ${\bf C}((q))$ is the Laurent polynomial
ring in the one variable $q$.  However, by definition,
$p_i$ is actually an element of ${\cal C}[[z^{-2\tilde\alpha}]]$.  Hence $p_{\tilde\alpha}$
is an element of
${\cal C}[[z^{-2\tilde\alpha}]]$ for each $\tilde\alpha\in\Sigma^+$.

\begin{theorem} Set $$p=\prod_{\tilde\alpha\in \Sigma^+}p_{\tilde\alpha}.$$
Then for each $\mu\in P^+(\pi)$,
$$\gra ({\cal X}(  c_{\mu}'))=p^{-1}\tau(-2\tilde\mu)p.\leqno{(6.4)}$$
\end{theorem}

\noindent{\bf Proof:} By Theorem 5.8, there exists an element  $p$ in $ {\cal
C}[[z^{-\tilde\alpha_i}|\alpha_i\in 
\pi^*]]$ such that $\gra ({\cal X}(c_{\mu}'))=p^{-1}\tau(-2\tilde\mu)p$
for all $\mu\in P^+(\pi)$.  
Since each $p_{\tilde\alpha}$ is an element of 
${\cal C}[[z^{-\tilde\alpha}]]$, it follows that $\prod_{\tilde\alpha\in 
\Sigma^+}p_{\tilde\alpha
}$ is an element of ${\cal C}[[z^{-\tilde\alpha_i}|\alpha_i\in 
\pi^*]]$.   Choose $p'$ in  ${\cal C}[[z^{-\tilde\alpha_i}|\alpha_i\in 
\pi^*]]$ so that $pp'=\prod_{\tilde\alpha\in 
\Sigma^+}p_{\tilde\alpha
}$. 

Set $k_i'=\prod_{\tilde\alpha\in 
\Sigma^+\setminus\{\tilde\alpha_i\}}p_{\tilde\alpha}.$  Note that 
$p_ik_i'=pp'$.  Choose $k_i$ as in Lemma 6.4.
It follows that 
$k_ip'=k_i'$ for each $i$ such that $\alpha_i\in \pi^*$.  
Write $k_i'=\sum_{\tilde\gamma\geq 0}d_{\tilde\gamma}'z^{-\tilde\gamma}$,
$k_i=\sum_{\tilde\gamma}d_{\tilde\gamma\geq 0}z^{-\tilde\gamma}$,
and  $p'=\sum_{\tilde\gamma\geq 0}a_{\tilde\gamma}z^{-\tilde\gamma}$.

Note that  by the definition of $k_i'$ and Lemma 6.5,
$$ \sum_{\tilde\gamma\geq 0}d_{\tilde\gamma}'z^{-\tilde 
s_i\tilde\gamma}=k_i'{\rm\ and \ }\sum_{\tilde\gamma\geq 0}d_{\tilde\gamma}z^{-\tilde 
s_i\tilde\gamma}=k_i$$ where equality holds in ${\cal C}[[z^{-\tilde\alpha_i}|\alpha_i\in 
\pi^*]]$.  It follows that 
$$\sum_{\tilde\gamma\geq 0}a_{\tilde\gamma}z^{-\tilde s_i\tilde\gamma}=p'$$
in ${\cal C}[[z^{-\tilde\alpha_i}|\alpha_i\in 
\pi^*]]$ for all simple reflections $\tilde s_i\in 
W_{\Theta}$.   Hence 
$$\sum_{\tilde\gamma\geq 
0}a_{\tilde\gamma}z^{-w\tilde\gamma}=p'\leqno{(6.5)}$$
in ${\cal C}[[z^{-\tilde\alpha_i}|\alpha_i\in 
\pi^*]]$
for all $w\in W_{\Theta}$. Suppose $\tilde\gamma>0$ is such that 
$a_{\tilde\gamma}$ is nonzero.  We can find $w\in W_{\Theta}$ such 
that $w\tilde\gamma<0$.   But then $z^{-w\tilde\gamma}\notin {\cal C}[[z^{-\tilde\alpha_i}|\alpha_i\in 
\pi^*]]$.   This contradicts 
(6.5).   Hence $p'$ must be a scalar and so 
$k_i$ is a scalar multiple of $k_i'$.  $\Box$

\medskip
Let $p$ be defined as in Theorem 6.7.   It follows that $p$ is the unique element of ${\cal
C}[[z^{-\tilde\alpha_i}|\alpha_i\in\pi^*]]$ such that (6.4) holds for all $\mu\in P^+(\pi)$. 
Recall the definition of the projection map ${\cal P}_{\cal A}$ given 
at the end of Section 3.  The next result extends Theorem 6.7 to
other elements in
$Z(\check  U)$ and, more generally, to elements of $\check U^B$.

\begin{corollary} For each $c\in \check U^B$,  
$$\gra ({\cal X} (  c))=p^{-1}\gra ({\cal P}_{\cal A}(c))p$$ 
where
$$p=\prod_{\tilde\alpha\in \Sigma^+}p_{\tilde\alpha}.$$ Moreover,  $c\notin 
B_+\check U$ if and only if
   ${\cal P}_{\cal A}(c)\neq 0$. 
\end{corollary}

\noindent
{\bf Proof:}  Let $c\in \check U^B$. The definition of the map ${\cal 
P}_{\cal A}$ ensures that  ${\cal P}_{\cal A}(\check U)$ 
is contained in ${\cal C}[\check {\cal A}]$. Thus $z^{\lambda}({\cal P}_{\cal A}(c))=0$ for 
all $\lambda\in P^+(2\Sigma)$ if and only if ${\cal P}_{\cal 
A}(c)=0$.  Hence, by Theorem 3.6, ${\cal X}(c)=0$ if and only if 
${\cal P}_{\cal A}(c)=0$.

Suppose  $c\in B_+\check U$. Then ${\cal X}(c)=0$. By the previous paragraph, ${\cal
P}_{\cal A}(c)=0$.  Thus the lemma follows trivially in this case.   Thus we may 
assume that $c\notin B_+\check U$. 

Let $c'\in \check 
{\cal A}N^+$ such that $c-c'\in B_+\check U$. In particular, ${\cal
X}(c)={\cal X}(c')$.  
We can write $\gra c'=(\gra an)+(\gra n')$
where
$a$ is in ${\cal C}[\check {\cal A}]$, $n$ is a weight vector in $N^+$ of weight 
$\gamma$, and $n'\in \sum_{\beta>\gamma}\check{\cal A}N^+_{\beta}$.
Note that $a\neq 0$ since $c\notin B_+\check U$.  
Choose $\lambda\in P^+(2\Sigma)$   such that  $z^{\lambda}(\gra
a)\neq 0$. Since $v_{\lambda}^r$ generates the $\gra U$ module $\bar M(\lambda)^r$, it follows that 
$v_{\lambda}^rn\neq 0$. Hence $$\zeta_{\lambda}^r(\gra c')\in 
v_{\lambda}^r(\gra c')+\sum_{\beta>0}v_{\lambda}^rN^+_{\beta}(\gra 
c')=z^{\lambda}(\gra 
a)v_{\lambda}^rn+\sum_{\beta>0}v_{\lambda}^rN^+_{\beta}.\leqno{(6.6)}$$
In particular, $\zeta_{\lambda}^r(\gra c')$ is nonzero.

Given $b\in B$, we have $c'b=cb+(c'-c)b\in bc+B_+\check U$.  Thus $c'b\in B_+\check U$ for all $b\in B_+$. Hence 
$$\zeta_{\lambda}^r(\gra c')B_+=0.$$    Therefore, by the
discussion following Lemma 5.5, 
$\zeta_{\lambda}^r(\gra c')$ is a scalar multiple of 
$\zeta_{\lambda}^r$.   
This forces $n$ to be an element of ${\cal C}$. Without loss of generality, we may 
assume that $n=1$. Thus  (6.6) implies that   $\zeta_{\lambda}^r(\gra c')=z^{\lambda}(\gra a)\zeta_{\lambda}^r$,
which is a nonzero multiple of $\zeta_{\lambda}^r$. 
Using the
definition of the map ${\cal P}_{\cal A}$ (end of Section 3), we see that $\gra  a=\gra ({\cal P}_{\cal A}(c'))=\gra ({\cal P_A}(c))$. Hence 
${\cal P_A}(c)$ is nonzero, which completes
the proof of the second assertion of the lemma.

Arguing as in Theorem 5.8, yields
$$\gra({\cal X}(c'))=p^{-1}\gra ({\cal P_A}(c'))p.$$ Hence $\gra ({\cal 
X}(c))=p^{-1}\gra ({\cal P_A}(c))p$. The first
assertion now follows from Theorem 6.7.
$\Box$

\medskip
An immediate consequence of Corollary 6.8 is that 
$$\check U^B\cap( \check {\cal A}N^+_++B_+\check U)\subset B_+\check U.$$ 
A similar argument switching the right and left actions yields that 
$$ \check U^B\cap( \check {\cal A}N^+_++\check UB_+)\subset \check UB_+.$$

 \section{Minuscule and Pseudominuscule Weights}
 
In this section, we find  ``small'' elements in $\check U^B$ which 
correspond to minuscule or pseudominuscule weights of $\Sigma$ and 
determine their radial components.
For all but three types of irreducible symmetric pairs $\gg,\gg^{\theta}$, 
this small element is $c_{\mu}'$ plus a constant term 
where $\tilde\mu$ is   either a minuscule or a pseudominuscule weight 
in $\Sigma$.  Most of this section is devoted to finding a suitable 
element in $\check U^B$ in the   remaining  three problematic cases.  
This involves a  
separate construction using fine information about finite 
dimensional $\adr U$ submodules  of $\check U$. 

Recall that $\gg,\gg^{\theta}$ is an irreducible symmetric pair
and $\Sigma$ is reduced.
It follows from the classification of irreducible symmetric pairs 
that $\Sigma$ is an irreducible root system corresponding to a simple 
Lie algebra as classified in [Hu, Chapter III]. 
 Before discussing elements of $\check U^B$,
we briefly review facts concerning minuscule and 
pseudomninuscule weights  associated to root systems 
of simple Lie algebras. Since we will be applying this information 
to the restricted root system, 
all the results will be stated with respect to $\Sigma$. 

A fundamental weight
$\beta$ is called minuscule with respect to the root system $\Sigma$ if
$$0\leq {{2(\beta,\alpha)}\over{(\alpha,\alpha)}}\leq 1\leqno{(7.1)}$$
for all $\alpha\in {\Sigma}^+$.  Since $\Sigma$ is simple, it
admits a minuscule weight if and only if it is not of type $E8,
F4, $ or $G2$.  It is straightforward to check that the
minuscule weights for $\Sigma$ are exactly the smallest fundamental weights
not contained in $Q^+(\Sigma)$.  In particular, a minuscule weight
$\beta$ satisfies the following condition:
\begin{enumerate}
\item[(7.2)] There does not exist $\gamma\in P^+(\Sigma)$ such that
$\beta-\gamma\in Q^+(\Sigma)$. \end{enumerate}

The longest root of $\Sigma$, when $\Sigma$ is simply laced, and
the longest short root of $\Sigma$, when $\Sigma$ has two root
lengths, is called a pseudominuscule weight.   A
pseudominuscule weight $\beta$ satisfies (7.1) for all $\alpha\in
\Sigma^+\setminus{\beta}$.  Moreover, one checks easily that a
pseudominuscule weight $\beta$ satisfies the following condition
similar to (7.2):
\begin{enumerate}
\item[(7.3)] The weight  $\gamma\in P^+(\Sigma)$ satisfies
$\beta-\gamma\in Q^+(\Sigma)$ if and only if $\gamma=0$.
\end{enumerate}

Recall the definition of the element $p$ given in  Theorem 6.7. Using the
previous sections, we can compute the image of
$c_{\mu}$ under ${\cal X}$ when $\mu$ satisfies one of the
above conditions. In particular, we have the following.

\begin{lemma} 
 (i) If $\mu\in P^+(\pi)$ is such that $\tilde\mu$
is minuscule in $\Sigma$ then      there exists  
	    a central element $c$ in $Z(\check U)$  with
 $${\cal X}(c)=\sum_{w\in 
	W_{\Theta}}w(p^{-1}\tau(-2\tilde\mu)p).$$

\noindent (ii)  If $\mu\in P^+(\pi)$ is such that $\tilde\mu$ is 
pseudominuscule in $\Sigma$  then   there exists   a central element $c$ 
in $Z(\check U)$ with
$${\cal X}(c)=\sum_{w\in 
	W_{\Theta}}w(p^{-1}\tau(-2\tilde\mu)p(1-\tau(2\tilde\mu))).$$
	\end{lemma}
	
	\noindent
	{\bf Proof:}   
Suppose first that  
 there exists $\mu\in P^+(\pi)$ such that $\tilde\mu$
is minuscule in $\Sigma$. Let $({\cal C}(Q_{\Sigma}){\cal A}_{\geq})_+$ denote the subalgebra of 
${\cal C}(Q_{\Sigma}){\cal A}_{\geq}$ generated by ${\cal C}(Q_{\Sigma})$ and ${\cal C}[{\cal A}_{\geq}]_+.$ By
Theorem 6.7 and Corollary 3.3, 
$${\cal X}(c_{\mu}')\in p^{-1}\tau(-2\tilde\mu)p+({\cal C}(Q_{\Sigma}){\cal 
A}_{\geq})_+ \tau(-2\tilde\mu).$$ Note that since $\tilde\mu$ satisfies (7.2), 
there does not exist $\tau(-2\tilde\beta)\in {\cal C}[{\cal A}_{\geq}]\tau(-2\tilde\mu)$ such
that $\tilde\beta$ is dominant.  Hence by Theorem 3.4, 
$${\cal X}(c_{\mu}')=\sum_{w\in W_{\Theta}}w(p^{-1}\tau(-2\tilde\mu)p).$$  This
proves (i).   

Now assume that there exists $\mu\in P^+(\pi)$ such that $\tilde\mu$
is pseudominuscule in $\Sigma$.  The same reasoning as in the previous paragraph, yields
that $${\cal X}(c_{\mu}')=\sum_{w\in W_{\Theta}}w(p^{-1}\tau(-2\tilde\mu)p)
+g$$ for some $g\in {\cal C}(Q_{\Sigma})$.  Now the zonal spherical
functions $\varphi_{\lambda}$ are eigenvectors for the action of ${\cal
X}(c_{\mu})$. When $\lambda=0$, $\varphi_{\lambda}$ is just $1$.   Hence
the action of ${\cal X}(c_{\mu}')$ on $1$ must be a scalar.   Now 
the action of $$\sum_{w\in W_{\Theta}}w(p^{-1}\tau(-2\tilde\mu)p)
-\sum_{w\in 
	W_{\Theta}}w(p^{-1}\tau(-2\tilde\mu)p\tau(2\tilde\mu))$$
on $1$ is zero.  Furthermore, $w(p^{-1}\tau(-2\tilde\mu)p\tau(2\tilde\mu))$ is in 
${\cal C}(Q_{
\Sigma})$ for each $w\in W_{\Theta}$. Hence $g+
\sum_{w\in 
	W_{\Theta}}w(p^{-1}\tau(-2\tilde\mu)p\tau(2\tilde\mu))$ is a scalar,
say $g_0$.   It follows that ${\cal X}(c_{\mu}'-g_0)$ has the required
form. $\Box$
\medskip

Assume for the moment that $\gg,\gg^{\theta}$ is not of type EIV, 
EVII, or EIX. Suppse that $\beta$ is a minuscule weight or 
pseudominuscule weight in $P^+(\Sigma)$.
A straightforward computation shows that $P^+(\pi)$ contains  a 
 fundamental weight $\mu$ such that $\beta=\tilde\mu$.
 These values of $\mu$ and $\tilde\mu$ are given in the appendix.  

Now assume that $\gg, \gg^{\theta}$ is of type EIV.   Then $\Sigma$ is of type A3
with set of simple roots $\{\tilde\alpha_1,\tilde\alpha_6\}$. Recall 
that $\omega_i$ denotes the fundamental weight corresponding to the 
simple root $\alpha_i$.
For $i=1$ and $i=6$, let $\omega'_i$ denote the fundamental weight in the
weight lattice of
$\Sigma$ corresponding to  
$\tilde\alpha_i$.  Note that both $\omega'_1$ and $\omega'_6$ are
minuscule.   It is straightforward to check that neither $\omega'_1$ nor
$\omega'_6$ is in the span of the set $\{\tilde\omega_i|1\leq i\leq 6\}$.
A similar computaton shows that if $\gg,\theta$ is of  type EVII or EIX
then the minuscule or pseudominuscule weights associated $\Sigma$ are not the 
restriction of  elements in $P^+(\pi)$. Thus, in these special 
 cases, $Z(\check U)$ does 
 not contain elements whose radial components are of the form 
 described in Lemma 7.1. The remainder of this section 
 is devoted to finding  elements in $\check U^B$  in the 
 remaining   cases which 
 play the role of $c_{\mu}$ for $\mu$ minuscule or pseudominuscule.

 We recall   basic facts   about the structure of
$\check U$ as a
$U$ module with respect to the adjoint action (see [JL1] and [JL2], or [J1, Section 7]). Here we use the
right adjoint action instead of the left and will translate the results accordingly. For each 
$\eta\in P^+(\pi)$, observe that  $\tau(
 2\eta)$ generates a finite dimensional $\adr U$ 
 module.  Let $F_r(\check U)$ denote the locally finite part of $\check 
 U$ with respect to the right adjoint action. One has that $F_r(\check U)$ is a
subalgebra of $\check U$.   As an $(\adr U)$ module,
$F_r(\check U)$ is   isomorphic  to the direct sum of the $(\adr
U)\tau(2\eta)$ as $\eta$ runs over elements in $P^+(\pi)$.  

 Let $G^+$ be the subalgebra of
$U$ generated by the $x_jt_j^{-1}$ for 
$j=1,\dots, n$.  
By [J2, Theorem 3.3], one can construct a basis of $(\adr U)\tau(2\eta)$
consisting of weight vectors  contained in sets of the form
$$a_{-\beta}b_{\beta'}\tau(2\eta)+\sum_{\xi\in
Q^+(\pi)}\sum_{\gamma\leq\beta-\xi}\sum_{\gamma'\leq\beta'-\xi}
U^-_{-\gamma}G^+_{\gamma'}\tau(2\eta-2\xi)\leqno{(7.4)}$$ where
$a_{-\beta}\in  U^-_{-\beta}$ and $b_{\beta'}\in
 G^+_{\beta'}$. Moreover, $0\leq \xi\leq \eta-w_0\eta$ and   both
$\beta$ and
$\beta'$ are less than or equal to
$\eta-w_0\eta$  ([J1,
7.1.20).

 \begin{lemma} Let $\eta\in P^+(\pi)$ and suppose that $u\in (\adr U)\tau(2\eta)$.  Then $${\cal
P}_{\cal A}(u)\in \tau(2\widetilde{w_0\eta}){\cal C}[{\cal
A}_{\geq}].\leqno{(7.5)}$$
\end{lemma}

\noindent
{\bf Proof:} Suppose that  
$u\in U^-_{-\gamma}G^+_{\gamma'}\tau(2\eta-2\xi)$ where $\xi$ and
$\gamma$ are    elements of 
$Q^+(\pi)$ such that $0\leq \xi\leq \eta-w_0\eta$ and $0\leq
\gamma\leq \eta-w_0\eta-\xi$.  We argue that
$u$ satisfies (7.5).  The lemma then follows by the  linearity of ${\cal
P}_{\cal A}$ and the description of a basis of $(\adr U)\tau(2\eta)$ using
(7.4).

  Note
that (2.7) implies that 
$ U^+_+\check U^0\subseteq {\cal M}^+_+ N^+\check U^0+N_+^+\check {\cal
A}$.  Hence
${\cal P}_{\cal A}(U^+_+\check U^0)=0$. It further follows from (3.11)
that  $$\check UU^+_+\check U^0\subseteq (B_+\check U+N^+_+\check{\cal
A}+{\cal C}[\check{\cal A}])U^+_+\check U^0\subseteq B_+\check
U+N^+_+\check{\cal A}.$$ Hence ${\cal P}_{\cal A}(U^-G^+_+\check
U^0)=0$.  
Thus, we may reduce to the case
where
$\gamma'=0$ and so 
$u\in U^-_{-\gamma}\tau(2\eta-2\xi)$.  It follows that $u\in
G^-_{-\gamma}\tau(2\eta-2\xi-\gamma)$.  

By Lemma 2.1 and (2.10), we have
that
$${\cal P}_{\cal A}(u)\in \tau(2\tilde\eta-2\tilde\xi-\tilde\gamma){\cal
C}[{\cal A}_{\geq}].$$  To complete the proof of the lemma, we argue that
$\tau(2\tilde\eta-2\tilde\xi-\tilde\gamma)\in
\tau(2\widetilde{w_0\eta}){\cal C}[{\cal A}_{\geq}]$.  Our assumptions on
$\xi$ and
$\gamma$ ensure that 
$2\eta-2\xi-\gamma\geq 2w_0\eta.$ Furthermore, $2\tilde\eta$, $2\tilde\xi$,
and $2w_0\tilde\eta$ are all elements of $Q(2\Sigma)$. Thus it is
sufficient to show that $\tilde\gamma\in Q(2\Sigma)$ whenever ${\cal
P}_{\cal A}(u)\neq 0$. 

Note that $u\tau(-2\eta+2\xi-\gamma)$ can be written as a
sum of  monomials
$y_{i_1}t_{i_1}\cdots y_{i_m}t_{i_m}$ in the $y_it_i$ of weight $\gamma$.  Assume that $y_{i_1}\notin
B_+$.
  Arguing as in the
proof of Lemma 2.1, we obtain
$$y_{i_1}t_{i_1}\cdots y_{i_m}t_{i_m}\in
-\tilde\theta(y_{i_1})t_{i_1}y_{i_2}t_{i_2}\cdots y_{i_m}t_{i_m}+
B_+y_{i_2}t_{i_2}\cdots y_{i_m}t_{i_m}.$$
Thus by (2.5), $$y_{i_1}t_{i_1}\cdots y_{i_m}t_{i_m}\in
G^-\tilde\theta(y_{i_1})t_{i_1}+G^-_{-\gamma+2\tilde\alpha_i}T'_{\geq}T_{\Theta}
+B_+U.$$   It follows that  ${\cal P}_{\cal A}(y_{i_1}t_{i_1}\cdots y_{i_m}t_{i_m})\in
{\cal P}_{\cal A}(G^-_{-\gamma+2\tilde\alpha_i}T'_{\geq}T_{\Theta}).$  Furthermore, 
$${\cal P}_{\cal A}(G^-_{-\gamma+2\tilde\alpha_i}T'_{\geq}T_{\Theta})={\cal P}_{\cal
A}(G^-_{-\gamma+2\tilde\alpha_i}){\cal P}_{\cal A}(T'_{\geq}T_{\Theta}).$$
  By induction on
$\gamma$, we obtain 
${\cal P}_{\cal A}(G^-_{-\gamma})=0$ unless $\tilde\gamma\in Q^+(2\Sigma)$.
Therefore ${\cal P}_{\cal A}(u)\neq 0$ implies $\tilde\gamma\in
Q^+(2\Sigma)$. $\Box$
\medskip

Let
$\phi$ be the Hopf algebra automorphism of
$U$ which fixes elements in
${\cal  M}T$ and $\phi(x_i)=q^{(-2\rho,\tilde\alpha_i)}x_i $ for all 
$1\leq i\leq n$.  Note
 that 
$ 
\phi(B_i)=
=q^{(\rho,-\Theta(\alpha_i)-\alpha_i)}(q^{(2\rho,\alpha_i)}y_it_i+
q^{(2\rho,\Theta(\alpha_i))}\tilde\theta(y_i)t_i)$ for all
$\alpha_i\notin\pi_{\Theta}$.
 The next lemma provides information about $(\adr 
B_+)\check U$ which will be applied later to elements of $\check U^B$.

\begin{lemma}  Suppose $u\in (\adr B_+)\check U$.   Then $u\in 
(\phi(B)\check T_{\Theta})_+\check U+\check U(B\check T_{\Theta})_+.$
\end{lemma}

\noindent
{\bf Proof:}  Suppose $u\in (\adr ({\cal M}T_{\Theta})_+)\check U$. 
Now ${\cal M}T_{\Theta}$ is  a Hopf subalgebra of $U$.   Hence 
$(\adr ({\cal M}T_{\Theta})_+)\check U\subset ({\cal
M}T_{\Theta})_+\check U+\check U({\cal
M}T_{\Theta})_+$.  The lemma now follows in this case since ${\cal M}T_{\Theta}$ is a subalgebra of both
$\phi(B)$ and $B$. Thus, it is sufficient to
show that 
$(\adr B_i)a\in (\phi(B)\check T_{\Theta})_+\check U+\check
U(B\check T_{\Theta})_+$ for all $i$ such that
$\alpha_i\notin\pi_{\Theta}$ and $a\in \check U$.

Fix $i$ with $\alpha_i\notin\pi_{\Theta}$.
By [L3, Theorem 7.1], there exists a sequence $i_1,\dots,i_s$  and a nonzero scalar $g$
such that $\tilde\theta(y_i)$ $=g(\adr x_{i_1})\dots (\adr
x_{i_s})(t_{p(i)}^{-1}x_{p(i)})$. Using the fact that $(\adr
x_j)a=-t_j^{-1}x_ja+t_j^{-1}ax_j$ for all $j$ and for all $a\in \check U$, we obtain 
$$ \tilde\theta(y_i)t_i\in {\cal
M}_+U+g\tau(\Theta(\alpha_i))x_{p(i)}x_{i_s}\cdots 
x_{i_1}t_i\leqno{(7.6)}$$
and $$\tilde\theta(y_i)t_i\in U{\cal M}_++(-1)^sgt_{i_1}^{-1}x_{i_1}\cdots
t_{i_s}^{-1}x_{i_s}t_{p(i)}^{-1}x_{p(i)}t_i.\leqno{(7.7)}$$ Note that 
$\tau(\Theta(-\alpha_i))t_i^{-1}\in T_{\Theta}$.
Applying the
antipode  to (7.7) and using (7.6) to simplify  yields
$$\eqalign{\sigma(\tilde\theta(y_i)t_i)&\in {\cal
M}_+U-
gt_i^{-1}t_{p(i)}^{-1}x_{p(i)}t_{p(i)}t_{i_s}^{-1}x_{i_s}t_{i_s}\cdots 
t_{i_1}^{-1}x_{i_1}t_{i_1}\cr
&=({\cal
M}T_{\Theta})_+U-q^{(2\rho,\Theta(\alpha_i))}g\tau(\Theta(\alpha_i))x_{p(i)}x_{i_s}\cdots 
x_{i_1}
\cr
&=({\cal 
M}T_{\Theta})_+U-q^{(2\rho,\Theta(\alpha_i))}\tilde\theta(y_i).\cr}$$
Hence 
$$\eqalign{\sigma(B_i)&=\sigma(y_it_i+\tilde\theta(y_i)t_i)\cr&=
-q^{(2\rho,\alpha_i)}y_i-q^{(2\rho,\Theta(\alpha_i))}\tilde\theta(y_i)\cr&=-q^{(\rho,\Theta(\alpha_i)+\alpha_i)}\phi(B_i)t_i^{-1}.\cr}$$

Recall that $\Delta(y_it_i)=y_it_i\otimes 1+t_i\otimes y_it_i$.   It follows 
from [L3, (7.14)] that
$${\it\Delta}(\tilde\theta(y_i)t_i)\in \tilde\theta(y_i)t_i\otimes
\tau(\Theta(\alpha_i))t_i+t_i\otimes
\tilde\theta(y_i)t_i+U\otimes ({\cal M}T_{\Theta})_+.$$
Note that $\tau(\Theta(\alpha_i))t_i\in T_{\Theta}$.
Thus $\tau(\Theta(\alpha_i)t_i\in 1+{\cal
C}[T_{\Theta}]_+$.  Combining the  formulas for
${\it\Delta}(y_it_i)$ and $\Delta(\tilde\theta(y_i))$ gives us  $${\it\Delta}(B_i)\in B_i\otimes 1+t_i\otimes
B_i+U\otimes ({\cal M}T_{\Theta})_+.$$
Hence $(\adr B_i)a\in 
-q^{(\rho,\Theta(\alpha_i)+\alpha_i)}\phi(B_i)t_i^{-1}a+t_i^{-1}aB_i+ ({\cal
M}T_{\Theta})_+\check U+\check U({\cal M}T_{\Theta})_+$ for all 
$a\in \check U$.$\Box$

\medskip

Recall the definition of the Hopf algebra automorphism $\chi$ following the proof of Theorem 4.5.
In particular,
$\chi(x_i)=q^{(\rho,\tilde\alpha_i)}x_i$
for all $1\leq i\leq n$.  
Since $(\rho,\tilde\alpha_i)=(\tilde\rho,\alpha_i)$ for all $i$, 
it follows that
$\tau(-\tilde\rho)x_i\tau(\tilde\rho)=\phi\chi(x_i)$ for $1\leq i\leq 
n$.
Similarly, $\tau(-\tilde\rho)y_i\tau(\tilde\rho)=\phi\chi(y_i)$ for all $1\leq i\leq n$.
Thus  $\tau(-\tilde\rho)u\tau(\tilde\rho)=\phi\chi(u)$ for all $u\in
 U$.

\begin{lemma} Suppose that $a'\in \check U^B$ and $a\in \check U^0$ 
such that  $$a'\in a+ ({\cal M}T_{\Theta})_+\check U+\check U({\cal M}T_{\Theta})_++(\adr
B)_+F_r(\check U).$$  Then 
$$(\varphi*{\cal X}(a'))(\tau(\tilde\rho))=\varphi(a\tau(\tilde\rho))\leqno{(7.8)}$$
for all $\varphi\in {\cal C}[P(2\Sigma)]^{W_{\Theta}}$.
\end{lemma}

\noindent
{\bf Proof:} By [L4, Corollary 5.4],  
$\{\varphi_{\lambda}|\ \lambda\in P^+(2\Sigma)\}$ is a basis for ${\cal 
C}[P(2\Sigma)]^{W_{\Theta}}$.   Hence it is sufficient to verify (7.8) 
when $\varphi=\varphi_{\lambda}$ for $\lambda\in P^+(2\Sigma)$. 
We use here facts from  Section 1 introduced before Theorem 1.1.  In 
particular,  $_{\phi\chi(B)}{\cal H}_{\phi(B)}$ is the subspace of 
$R_q[G]$ consisting of left $\phi\chi(B)$ and right $\phi(B)$ 
invariants. 
Moreover, the space $_{\phi\chi(B)}{\cal H}_{\phi(B)}$  contains a distinguished 
basis $\{g'_{\lambda}|\ \lambda\in P^+(2\Sigma)\}$ such that 
the set $\{\Upsilon(g'_{\lambda})|\ \lambda\in 
P^+(2\Sigma)\}$  satisfies (1.2). For each $\lambda\in P^+(2\Sigma)$, 
$\Upsilon(g'_{\lambda})$ is written as
$\varphi^{\lambda}_{\phi(B),{\phi\chi(B)}}$. By [L4, Theorem 6.3], 
$\varphi^{\lambda}_{\phi(B),{\phi\chi(B)}}=\varphi^{\lambda}_{B,B'}=\varphi_{\lambda}$ 
for all $\lambda\in P^+(2\Sigma)$.

 The discussion preceding this lemma implies that $$a'
\tau(\tilde\rho)\in a\tau(\tilde\rho)+\phi(B)_+\check U+\check U\phi\chi(B)_+.$$ Hence
$$g'_{\lambda}(a'\tau(\tilde\rho))=\varphi_{\lambda}(a\tau(\tilde\rho)).$$
 Recall that $g'_{\lambda}\in L(\lambda)^{\phi\chi(B)}\otimes 
 (L(\lambda)^*)^{\phi(B)}.$   Thus 
 arguing as in Section 3 (preceding Theorem 3.6), 
 $$g'_{\lambda}(a'\tau(\tilde\rho))=z^{\lambda}({\cal P}_{\cal 
  A}(a'))g'_{\lambda}(\tau(\tilde\rho)).\leqno{(7.9)}$$
  The lemma now follows from Theorem 3.6 and the fact that 
  $g'_{\lambda}(\tau(\tilde\rho))=\varphi_{\lambda}(\tau(\tilde\rho))$. $\Box$

\medskip
By [L4, Theorem 3.1], $F_r(\check U)$ can be written as a direct sum of finite
dimensional simple $(\adr B)$ modules.   Consider an element $a\in F_r(\check U)$.   Since the
action of
$\adr B$ on $F_r(\check U)$ is locally finite, $(\adr B)a$ is a finite
dimenionsal submodule of $F_r(\check U)$.  Now suppose that
$a\notin(\adr B)_+\check U$.  In particular,    $(\adr B_+)a$ has
codimension $1$ in $(\adr B)a$.   Since the action of $(\adr B)$ on
$F_r(\check U)$ is completely reducible, it follows that there exists a
nonzero element $a'\in \check U^B$ such that $$(\adr B)a={\cal
C}a'\oplus (\adr B_+)a.$$  We may further assume that $a'$ has been chosen so that $$a'\in a+(\adr B_+)a.$$
Now $(\adr B_+)\check U\cap \check U^B=0$.   Hence the choice of $a'$ is unique.  Thus we have a linear
map $L:F_r(\check U)\rightarrow \check U^B$ such that  $L(a)$ is the unique element of $\check U^B$ in the set
$a+(\adr B_+)a$.  (Note that $L$ is very similar to the so-called 
Letzter map studied in [J3].) 

Recall that $W_{\Theta}$ acts on    ${\cal C}[ \check {\cal A}]$ (see 
the discussion  following the proof of Corollary 3.3.)  Given $w\in W_{\Theta}$ and $u\in 
{\cal C}[ \check {\cal A}]$ we write $w\cdot u$ for the action of $w$ on $u$.

\begin{lemma}  Suppose that $a'$ is an element of $F_r(\check U)$ such that $a'\in a
	+({\cal M}T_{\Theta})_+\check U+\check U({\cal M}T_{\Theta})_+$ 
	where $a\in {\cal C}[\check {\cal A}]$ and 
 $$\sum_{w\in 
 W_{\Theta}}w\cdot (a\tau(\tilde\rho))\neq 0$$   Then
${\cal P}_{\cal A}(L(a'))$ is nonzero.
\end{lemma}

\noindent
{\bf Proof:}  It is straightforward to check that the bilinear
pairing between 
${\cal C}[P(2\Sigma)]^{W_{\Theta}}$ and  ${\cal C}[\check {\cal A}]^{W_{\Theta}}$ given by
$<f,k>=f(k)$ for all $f\in {\cal
C}[P(2\Sigma)]^{W_{\Theta}}$ and $k\in {\cal C}[\check {\cal A}]^{W_{\Theta}}$  is nondegenerate.   Note that 
$$\varphi(a\tau(\tilde\rho))=|W_{\Theta}|^{-1}\varphi(\sum_{w\in W_{\Theta}}
w\cdot (a\tau(\tilde\rho)))$$ for all $\varphi\in {\cal 
C}[P(2\Sigma)]^{W_{\Theta}}$.
Hence the assumptions on $a$ ensure that there exists $\varphi\in {\cal
C}[P(2\Sigma)]^{W_{\Theta}}$ such that $\varphi(a\tau(\tilde\rho))\neq 0$.

   By the previous lemma and (7.9), 
$g_{\lambda}(L(a')\tau(\tilde\rho))=z^{\lambda}({\cal P}_{\cal 
A}(L(a'))g_{\lambda}(\tau(\tilde\rho))=\varphi_{\lambda}(a\tau(\tilde\rho))$.   Thus 
it is sufficient to find $\lambda$ such that
$\varphi_{\lambda}(a\tau(\tilde\rho))$ is nonzero.   But the 
$\{\varphi_{\lambda}|\lambda\in P^+(2\Sigma)\}$ form a basis for ${\cal
C}[P(2\Sigma)]^{W_{\Theta}}$. The result now follows using the nondegenerate pairing described in the first
paragraph.
$\Box$

\medskip
We are now ready to associate an element of $\check U^B$ to a 
minuscule or pseudominuscule restricted weight.

\medskip
 \begin{lemma}  There exists  $a\in \check U^B$ such that $${\cal P}_{\cal A}(a)\in
\tau(-2\tilde\eta)+\tau(-2\tilde\eta){\cal C}[{\cal 
A}_{\geq}]_+\leqno{(7.10)}$$ where $\tilde\eta$ is a minuscule or 
pseudominuscule weight.
\end{lemma}

\noindent
{\bf Proof:}  By the discussion preceding Theorem 
3.6, ${\cal P}_{\cal A}$ restricted to $Z(\check U)$ agrees with the 
composition of ${\cal P}$ followed by the projection onto ${\cal C}[\check {\cal A}]$ using 
(3.5).  Suppose that $\eta\in P^+(\pi)$ such that 
$\tilde\eta$ is a minuscule or pseudominuscule restricted weight with 
respect to $\Sigma$.
Then the description (4.8) of the image of the central elements 
under ${\cal P}$ ensures that $c_{\eta}$ satisfies (7.10).   Thus we may 
reduce to the cases when $\gg,\gg^{\theta}$ is of type EIV, EVII, or EIX. 
We assume first that $\gg,\gg^{\theta}$ is of type EIV.

Recall that $\omega_i$ is the fundamental weight corresponding to the 
simple root $\alpha_i$ in $\pi$ and $\omega_i'$ is the fundamental weight
corresponding to the restricted root $\tilde\alpha_i$ in $\Sigma$. One 
checks that $\tilde\omega_6=2\omega_6'$, $\tilde\omega_1=2\omega_1'$, and 
$w_0\omega_6=\omega_1$.
Consider the
$(\adr U)$ submodule
$(\adr U)\tau(2\omega_6)$ of
$F_r(\check U)$.   Note that  the only restricted dominant integral 
weights less than $2\omega_1'$ is $\omega_6'$. 

By Lemma 7.2,  
$${\cal P}_{\cal A}(a)\in
\tau(-4\omega_1'){\cal C}[{\cal A}_{\geq}]_+$$ for any $a\in \check U^B\cap (\adr
U)\tau(2\omega_6)$.  On the other hand, Corollary 6.8 and Theorem 
3.6 ensure that $\gra({\cal P}_{\cal 
A}(a))$ is a linear combination of terms of the form $\tau(-2\gamma)$ 
where $\gamma\in P^+(\Sigma)$.  Hence if $a\in \check U^B\cap (\adr 
U)\tau(2\omega_6)$, then 
$${\cal P}_{\cal A}(a)\in
\tau(-2\gamma)+\tau(-2\gamma){\cal C}[{\cal A}_{\geq}]_+$$  where 
 $\gamma=\omega_6'$, or
$\gamma=2\omega_1'$. Thus it is sufficient to find two elements $a_1$ and $a_2$ in $\check U^B\cap (\adr
U)\tau(2\omega_6)$ such that ${\cal P}_{\cal A}(a_1)$ and ${\cal P}_{\cal 
A}(a_2)$ are linearly independent.

Note that $\tau(2\omega_6)$ is in $(\adr U)\tau(2\omega_6)$.   It is 
straightforward to check that 
$4\omega_6'-2\tilde\alpha_6=2\omega_1'$.  Hence 
$\tau(2\omega_1')\in \tau(2\omega_6)\tau(-2\alpha_6)+{\cal C}[\check 
T_{\Theta}]_+.$  

The action $\adr $ is a right action. In particular, $(\adr uv)w=(\adr 
v)(\adr u)w$ for all $u,v,$ and $w\in \check U$.  A direct computation 
yields that $$\eqalign{(\adr 
y_6x_6)&\tau(2\omega_6)-q^{-1}\tau(2\omega_6)\cr&=-q^{-2}(1-q^2)^2\tau(2\omega_6)t_6^{-1}y_6x_6-
q^{-1}\tau(2\omega_6)t_6^{-2}.\cr}$$
Recall that $x_5$ and $y_5$ are elements of ${\cal M}$.
Arguing as in Lemma 4.3, one has, up to a nonzero scalar, that
$y_6y_5x_5x_6$ is an element of $y_6x_6+({\cal M}T_{\Theta})_+ U+
U({\cal M}T_{\Theta})_+.$  Using this fact, it is 
straightforward to show that up to a nonzero scalar $$(\adr y_5 x_5)(\adr 
y_6x_6)\tau(2\omega_6)\in 
\tau(2\omega_6)t_6^{-1}y_6x_6+({\cal M}T_{\Theta})_+\check U+\check 
U({\cal M}T_{\Theta})_+.$$
Thus a suitable linear combination of the terms $\tau(2\omega_6)$, $(\adr 
y_6
x_6)\tau(2\omega_6)$ and $(\adr y_6x_6y_5x_5)
\tau(2\omega_6)$ yields an element $b\in (\adr 
U)\tau(2\omega_6)$ such that   $$b\in \tau(2\omega_1')+({\cal M}T_{\Theta})_+\check 
U+\check U({\cal M}T_{\Theta})_+.$$

Set $a_1=L(\tau(2\omega_6))$ and $a_2=L(b)$.   Note that
$2\omega_1'+\tilde\rho$ and 
$4\omega_6'+\tilde\rho$ are distinct dominant 
restricted weights
and hence are in different orbits with 
respect to the action of $W_{\Theta}$.   It follows from Lemma 
7.5 that ${\cal P}_{\cal A}(a_1)$, and ${\cal P}_{\cal 
A}(a_2)$ are linearly independent.   This completes the proof of the lemma when 
$\gg,\gg^{\theta}$ is of type EIV for the restricted minuscule weight $\omega_1'$.
The same argument works for the restricted minuscule weight 
$\omega_6'$ using the diagram automorphism of E$_6$ which sends 
$\alpha_1$ to $\alpha_6$.

Now consider the case when 
$\gg,\gg^{\theta}$ is of type EVII. In this case, we use the $(\adr 
U)$ module $(\adr U)\tau(2\omega_7)$. Note that 
$\omega_7'=\tilde\omega_7$ and $\omega_1'$ are the only dominant 
integral restricted weights less than or equal to $\tilde\omega_7$. Furthermore, as above, there exists a  linear  combination $b$ 
of $\tau(2\omega_7)$, $(\adr 
y_7y_6x_7
x_6)\tau(2\omega_7)$ and $(\adr y_7y_6x_7x_6y_5x_5)
\tau(2\omega_7)$  such that  $$b\in \tau(2\omega_1')+{\cal 
C}\tau(2(\tilde\omega_7-\tilde\alpha_7))+({\cal M}T_{\Theta})_+\check 
U+\check U({\cal M}T_{\Theta})_+.$$ One checks that 
$2\omega_1'+\tilde\rho,2(\tilde\omega_7-\tilde\alpha_7)+\tilde\rho, $ and 
$\tau(2\tilde\omega_7)+\tilde\rho$ are in different $W_{\Theta}$ orbits. 
Hence Lemma 7.5 ensures that ${\cal P}_{\cal A}(L(\tau(2\omega_7)))$ and ${\cal P}(L(b))$ 
are linearly independent and the proof follows in this case. Similarly, if 
$\gg,\gg^{\theta}$ is of type EIX, one can find a  linear 
combination $b$ of 
$\tau(2\omega_8)$, $(\adr 
y_8y_7y_6x_8
x_7x_6)\tau(2\lambda)$ and $(\adr y_8y_7y_6x_8
x_7x_6y_5x_5)
\tau(2\lambda)$  such that $b$ is an element of $$\tau(2\omega_1')+{\cal 
C}\tau(2(\tilde\omega_8-\tilde\alpha_8))+{\cal 
C}\tau(2(\tilde\omega_8-\tilde\alpha_8-\tilde\alpha_7))+({\cal M}T_{\Theta})_+\check 
U+\check U({\cal M}T_{\Theta})_+.$$ In this case, there are three 
dominant integral weights, $1$, $\omega_1'$, and $\omega_8'$, less than or equal to 
$\tilde\omega_8=\omega'_8$. Now the $W_{\Theta}$ orbits of $\tilde\rho$, $2\omega_1'+\tilde\rho$,
$2(\tilde\omega_8-\tilde\alpha_8)+\tilde\rho,(\tilde\omega_8-\tilde\alpha_8-\tilde\alpha_7)+\tilde\rho$,
 and 
$2\tilde\omega_8+\tilde\rho$ are  distinct. This combined with Lemma 
7.5 implies that 
the three elements ${\cal P}_{\cal A}(1), {\cal P}_{\cal 
A}(L(\tau(2\omega_8)),$ and ${\cal P}_{\cal A}(L(b))$ are linearly 
independent.  The desired element $a$ in $\check U^B$ which 
satisfies (2.10) is then a suitable linear combination of $1$, 
$L(\tau(2\omega_8)$, and $L(b)$. $\Box$

\medskip
Using the proof of  Lemma 7.1, we can compute the radial 
components of the ``small'' elements in $\check U^B$ described   in Lemma 7.6.

\begin{theorem}  Suppose that $\tilde\mu$ is a minuscule weight of $\Sigma$.  Then there
exists an element $c\in \check U^B$ such that $${\cal X}(c)=\sum_{w\in 
	W_{\Theta}}w(p^{-1}\tau(-2\tilde\mu)p).$$ If $\Sigma$ is of type F4, G2, or E8,
and $\tilde\mu$ is the pseudominuscule weight of $\Sigma$, then there 
  exists   an element $c$ in $\check U^B$ such
that 
$${\cal X}(c)=\sum_{w\in 
	W_{\Theta}}w[p^{-1}\tau(-2\tilde\mu)p(1-\tau(2\tilde\mu))].$$
\end{theorem}

\noindent
{\bf Proof:}  If $\gg,\gg^{\theta}$ is not of type EIV, EVII, or EIX, then the result
follows from Lemma 7.1 and the paragraph immediately following the lemma. More generally, let $\tilde\mu$ 
be a minuscule or pseudominuscule weight of $\Sigma$.  Then Lemma 7.6 
guarantees the existence of   an element 
$c\in \check U^B$ such that $\gra {\cal P}_{\cal 
A}(c)=\tau(-2\tilde\mu)$.   Arguing as in the proof of Lemma 
7.1, using Corollary 6.8 and Theorem 3.6 (instead of Theorem 6.7 and 
Corollary 3.3) we see that ${\cal X}(c)$ has the desired form. 
$\Box$

\section{Macdonald Polynomials}
In this section we express the zonal spherical functions as Macdonald
polynomials. Formally, Macdonald polynomials associated to a root system 
are  
Laurent polynomials which depend on an indeterminate $a$ and a system of 
parameters $g=\{g_{\alpha}|\alpha $ is a root$\}$ such that 
$g_{\alpha}=g_{w\alpha}$ for all $w$ in the corresponding Weyl group.  Often, however, the 
$g_{\alpha}$ are taken to be powers of $a$.   We will take this point 
of view here in reviewing  basic facts and notations concerning
these polynomials.  To further simplify the presentation, we assume 
that
$a$ is a power of $q$ and  that the $g_{\alpha}$ are rational powers of $a$.
In particular, each $g_{\alpha}$ is an element of ${\cal C}$.

 Let $\check \Sigma$ denote the dual root system 
to $\Sigma$. 
For each $\lambda\in P^+(2\Sigma)$, set $$m_{\lambda}=\sum_{w\in
W_{\Theta}}z^{w\lambda}.$$
The set of  Macdonald polynomials $\{P_{\lambda}(a,g)|\lambda\in
P^+(2\Sigma)\}$ associated to $2\Sigma,\check 2\Sigma$ is  the unique basis 
of ${\cal
C}[P(2\Sigma)]^{W_{\Theta}}$ 
which satisfy the following conditions:
\begin{enumerate}
\item[(8.1)] The polynomials are orthogonal with respect to the Macdonald inner product at $a,g$.
\item[(8.2)]  There exists scalars $a_{\lambda,\mu}$ in ${\cal C}$ such 
that
$P_{\lambda}(a,g)=m_{\lambda}+\sum_{\mu<\lambda} a_{\lambda,\mu}m_{\mu}$
for all $\lambda\in P^+(2\Sigma)$.
\end{enumerate}   For more details, the reader is referred to [M2] or 
[Ki].

Macdonald polynomials can also be 
characterized using certain difference operators associated to 
minuscule and pseudominuscule weights.  
In particular,
set $a_{\tilde\alpha}= a^{2(\tilde\alpha,\tilde\alpha)}$ for $\tilde\alpha\in \Sigma$.
Recall (6.3) the definition of 
$(x;a)_{\infty}.$ Let
$${\Delta}^+_{a,g}=\prod_{\tilde\alpha\in
\Sigma^+}{{(z^{2\tilde\alpha};a_{\tilde\alpha})_{\infty}}\over{(g_{\tilde\alpha}z^{2\tilde\alpha};
a_{\tilde\alpha})_{\infty}}}.$$
Given $\beta\in P(2\Sigma)$, define an operator $T_{\beta}$ on
${\cal C}[Q(2\Sigma)]$ by
$$T_{\beta}z^{\alpha}=q^{(\beta,\alpha)}z^{\alpha}.$$  Let
$D_{\beta}(a,g)$ be the operator on ${\cal C}[Q(2\Sigma)]$
defined by
$$D_{\beta}(a,g)(f)=\sum_{w\in
W_{\Theta}}(\Delta^+_{a,g})^{-1}\left({T_{\beta}(\Delta^+_{a,g}f)}
\right).$$
Similarly, let $E_{\beta}(a,g)$ be the operator on
${\bf C}[a,g][z^{\alpha}|\alpha\in \Sigma]$
defined by
$$E_{\beta}(a,g)(f)=\sum_{w\in
W_{\Theta}}({\Delta^+_{a,g}})^{-1}\left({(T_{\beta}(\Delta^+_{a,g}))((T_{\beta}-1)f)}\right).$$

Set $F_{\beta}$ equal to the difference operator $D_{\beta}$ if $\beta$ 
is a minuscule weight of $2\Sigma$ and
equal to the difference operator $E_{\beta}$ if $\beta$ is 
pseudominuscule of $2\Sigma$. 
The following result summarizes the relationship between Macdonald polynomials and these
difference operators.

\begin{theorem} ([M2])  Let $S$ be a subset of $2\Sigma$ consisting of 
\begin{enumerate} 
\item[(i)] one minuscule
weight if $2\Sigma$ is not of type $D_n$, E8, F4, or G2
\item[(ii)] both minuscule weights if $2\Sigma$ is of type $D_n$
\item[(iii)] one pseudominuscule weight   if
$\Sigma$ is   of type   E8, F4, or G2.
\end{enumerate}
  Then the
set $\{P_{\lambda}(a,g)|\lambda\in P^+(2\Sigma)\}$ is the unique basis of
${\cal C}[P(2{\Sigma})]^W$ satisfying (8.2) which consists of eigenvectors for the action of 
$F_{\beta}$ as $\beta$ ranges over $S$.  
\end{theorem}

\medskip
  In order to identify the zonal
spherical functions with Macdonald polynomials, we show that the action of certain central
elements (or more generally $B$ invariant elements) on the right correspond to the left action of 
the difference operators described above.

\begin{theorem}
Let $\gg,\gg^{\theta}$ be an irreducible symmetric pair.   Let
$\{\varphi_{\lambda}|\lambda\in P^+(2\Sigma)\}$ denote the unique
$W_{\Theta}$ invariant zonal spherical family associated to ${\cal B}$.  
Then $$\varphi_{\lambda}=P_{\lambda}(a,g) $$ for all $\lambda\in P^+(2\Sigma)$.
Here
$a
=q$ and
$g_{\tilde\alpha}=q^{mult(\tilde\alpha)(\tilde\alpha,\tilde\alpha)(\alpha,\alpha)/2}$
for each $\tilde\alpha\in \Sigma$. 
\end{theorem}

\noindent
{\bf Proof:} Let  $\tilde\mu $ be a 
  minuscule weight in $\Sigma$ if $\Sigma$ is not of type F4, EVIII, or G2 and
pseudominuscule otherwise. Note that when $\tilde\mu$ is minuscule, so is 
$w_0'\tilde\mu$ and similarly if $\tilde\mu$ is pseudominuscule then so is 
$w_0'\tilde\mu$. Consider first the case when $\tilde\mu$ is minuscule.   Applying 
Theorem  7.7 to the weight $w_0'\tilde\mu$ shows that  there exists
$c\in \check U^B$ such that
$${\cal X}(c)=\sum_{w\in 
	W_{\Theta}}w(p^{-1}\tau(-2{w_0'\tilde\mu})p).\leqno{(8.3)}$$  It is straightforward to
see that 
$w_0'(p^{-1})=\Delta^+_{a,g}$ using the definition of $p$ given in Theorem 6.7. Hence we
can rewrite (8.3) as
$${\cal X}(c)=\sum_{w\in 
	W_{\Theta}}w((\Delta_{a,g}^{+}\tau(2{\tilde\mu})(\Delta_{a,g}^+)^{-1})$$    It follows
that the  right action of
${\cal X}(c)$ agrees with
  the left action of 
$D_{2\tilde\mu}(a,g)$.  Thus the basis
$\{\varphi_{\lambda}|\lambda\in P^+(2\Sigma)\}$ is a basis of eigenvectors for the action
of
$D_{2\tilde\mu}(a,g)$ satifying (8.2) when $\mu$ is minuscule.  The theorem now follows
from Theorem 8.1 when
$\tilde\mu$ is minuscule.  A similar argument using the operators $E_{2\tilde\mu}(a,g)$
works for
$\tilde\mu$ pseudominuscule.
$\Box$

\bigskip
\section{Appendix}
We include in this appendix a complete list of all irreducible 
symmetric pairs with reduced restricted roots using Araki's 
classification [A].  For each case, we  explicitly describe the 
involution $\Theta$ on $\Delta$ and the restricted roots $\Sigma$.
(This information can also be read off the table in [A]).   We further 
give the values of 
$g_{\tilde\alpha}$ and $a_{\tilde\alpha}$ for  $\alpha\in \pi^*$
in each case.
Moreover, if $\gg,\gg^{\theta}$ is not of type EIV,EVII, or EIX, then we 
explicitly give $\mu\in P^+(\pi)$ such that $\tilde\mu$ is the 
pseudominuscule or minuscule weight of $\Sigma$.

Note that in Case 1 below, when
$\gg_1$ and
$\gg_2$ are simply laced, we have $a_{\tilde\alpha}=g_{\tilde\alpha}$ for
all restricted roots $\tilde\alpha$.   In particular,  $a=g$ in these
cases and the resulting Macdonald polynomials are just the Weyl character
formulas.

\medskip
\noindent
{\bf Case I:} $\gg= \gg_1\oplus \gg_2$ where both $\gg_1$ and $\gg_2$ 
are simple and isomorphic to each other. Let $\pi_1=\{\alpha_1,\dots, 
\alpha_n\}$ denote the simple roots corresponding to $\gg_1$ and 
$\pi_2=\{\alpha_{n+1},\dots,\alpha_{2n}\}$ denote the roots for 
$\gg_2$.  We further assume that $\alpha_i\mapsto \alpha_{i+n}$ 
defines the isomorphism between the two root systems.

\noindent
$\Theta(\alpha_i)=-\alpha_{i+n}$ for $1\leq i\leq n$;
$\Sigma=\{\tilde\alpha_1,\dots, \tilde\alpha_n\}$ and is of the same 
type as $\gg_1$;
 $\tilde\mu$ is minuscule (resp. pseudominuscule) when $\mu$ is a
minuscule (resp. pseudominuscule)
weight of $\gg_1$.

\noindent
$a_{\tilde\alpha_i}=q^{(\alpha_i,\alpha_i)}$  and
$g_{\tilde\alpha_i}=q^{(\alpha_i,\alpha_i)(\alpha_i,\alpha_i)/2} $
for each $i$.

\medskip
\noindent
{\bf Case II:} $\gg$ is simple

\noindent
{\bf Type AI}
 $\gg$ is of type $A_n$.
 
 \noindent
 $\Theta(\alpha_i)=-\alpha_i$, for $1\leq i\leq n$.
 $\Sigma=\{\alpha_1,\dots, \alpha_n\}$ is of type $A_n$.

 \noindent
 minuscule weight: $\tilde\mu=\omega_1$  where $\mu=\omega_1$.
 
 \noindent
$a_{\tilde\alpha_i}=q^4$ and
 $g_{\tilde\alpha_i}=q^2$ for all $1\leq i\leq n$.

 \medskip
\noindent
{\bf Type AII}
 $\gg$ is of type $A_n$, where $n=2m+1$ is odd and $n\geq 3$.

 \noindent
 $\Theta(\alpha_{i})=\alpha_{i}$ for $i=2j+1$, $0\leq j\leq m$,
 $\Theta(\alpha_i)=-\alpha_{i-1}-\alpha_i-\alpha_{i+1}$ for 
 $i=2j, 1\leq j\leq m$.
 $\Sigma=\{\tilde\alpha_i|\ i=2j+1, 0\leq j\leq m\}$ is of type $A_n$.
 
 \noindent
 minuscule weight: $\tilde\mu=\omega'_1$ is minuscule  where $\mu=\omega_1$.

 \noindent
 $a_{\tilde\alpha_i}=q^2$ and
 $g_{\tilde\alpha_i}=q^4$ for all $i$, $i=2j, 1\leq j\leq m$.

\medskip
\noindent
{\bf Type AIII}
{\bf{\it Case 2}}:  $\gg$ is of type $A_n$ where $n=2m+1$. 

\noindent
$\Theta(\alpha_i)=-\alpha_{n-i+1}$ for $1\leq i\leq n$.
$\Sigma=\{\tilde\alpha_i|\ 1\leq i\leq m\}$ is of type $C_m$.

\noindent
minuscule weight: $\tilde\mu=\omega'_1$  where $\mu=\omega_1$.

 \noindent
$a_{\tilde\alpha_i}=q^2$ for $1\leq i\leq m-1$, 
$a_{\tilde\alpha_m}=q^4$, and
$g_{\tilde\alpha_i}=q^2$ for $1\leq i\leq m$.

\eject
\noindent
{\bf Type BI}  
 $\gg$ is of type $B_n$, $r$ is an integer such that $2\leq r\leq n$.

 \noindent
 $\Theta(\alpha_i)=\alpha_i$ for $r+1\leq i\leq n$,
 $\Theta(\alpha_i)=-\alpha_i$  for $1\leq i\leq r-1$,
 $\Theta(\alpha_r)=-2\alpha_{r+1}-2\alpha_{r+2}-\cdots-2\alpha_n.$ 
 $\Sigma=\{\tilde\alpha_1,\dots, \tilde\alpha_r\}$ is of type $B_r$.
 
 \noindent
 minuscule weight:  $\tilde\mu=\omega'_r$ where $\mu=\omega_n$.
 
 \noindent
 $a_{\tilde\alpha_i}=q^4$ and $g_{\tilde\alpha_i}=q^2$ for $1\leq i\leq r-1$, 
 $a_{\tilde\alpha_r}=q^2$ and
 $g_{\tilde\alpha_r}=q^{2(n-r)+1}$.
 
 \medskip
 \noindent
{\bf Type BII}  
$\gg$ is of type $B_n$.

\noindent
 $\Theta(\alpha_i)=\alpha_i$ for $2\leq i\leq n$,
$\Theta(\alpha_1)=-2\alpha_{2}-2\alpha_{3}-\cdots-2\alpha_n.$
$\Sigma=\{\tilde\alpha_1\}$ is of type $A_1$.

\noindent
minuscule weight: $\tilde\mu=\omega'_1$ where $\mu=\omega_n$.

\noindent
 $a_{\tilde\alpha_1}=q^2$ and
$g_{\tilde\alpha_1}=q^{2n-1}$.

\medskip
\noindent
{\bf Type CI}
 $\gg$ is of type $C_n$.
 
 \noindent $\Theta(\alpha_i)=-\alpha_i$  for $1\leq 
i\leq n$.
$\Sigma=\{\tilde\alpha_1,\dots, \tilde\alpha_n\}$ is of type $C_n$.

\noindent
minuscule weight: $\tilde\mu=\omega'_1$ where $\mu=\omega_1$.

\noindent
$a_{\tilde\alpha_i}=q^4$ and $g_{\tilde\alpha_i}=q^2$ for $1\leq i\leq n-1$, 
$a_{\tilde\alpha_n}=q^8$, and
$g_{\tilde\alpha_n}=q^4$.

\medskip
\noindent
{\bf Type CII}
{\bf{\it Case 2}}:   $\gg$ is of type $C_n$ where $n\geq 3$ is even.

\noindent
$\Theta(\alpha_i)=\alpha_i$ for $i=2j-1, 1\leq j\leq n/2,$
$\Theta(\alpha_i)=-\alpha_{i-1}-\alpha_i-\alpha_{i+1}$
for $i=2j, 1\leq j\leq (n-2)/2$, 
$\Theta(\alpha_n)=-\alpha_n-2\alpha_{n-1}$.
$\Sigma=\{\tilde\alpha_i|\ i=2j, 1\leq j\leq n/2\}$ is of type 
$C_{n/2}$.

\noindent
minuscule weight: $\tilde\mu=\omega'_2$ where $\mu=\omega_1$.

\noindent
$a_{\tilde\alpha_i}=q^2$ and $g_{\tilde\alpha_i}=q^4$ for $i=2j, 1\leq j\leq (n-2)/2$,
$a_{\tilde\alpha_{n}}=q^4$ and
$g_{\tilde\alpha_{n}}
=q^{12}$.

\medskip
\noindent
{\bf Type DI}
{\bf{\it Case 1}}:  
  $\gg$ is of type $D_n$,  $r$ is an integer such that $2\leq r\leq 
  n-2$
  
  \noindent
$\Theta(\alpha_i)=\alpha_i$ for $r+1\leq i\leq n$, 
$\Theta(\alpha_i)=-\alpha_i$ for 
 $1\leq i\leq r-1$,
$\Theta(\alpha_r)=-\alpha_{r}-2\alpha_{r+1}-\cdots 
-2\alpha_{n-2}-\alpha_{n-1}-\alpha_n$. 
$\Sigma=\{\tilde\alpha_1,\dots, \tilde\alpha_r\}$ is of type $B_r$.

\noindent
minuscule weight: $\tilde\mu=\omega'_r$ where $\mu=\omega_n$.

\noindent
$a_{\tilde\alpha_i}=q^4$ and $g_{\tilde\alpha_i}=q^2$ for $1\leq i\leq r-1$,
$a_{\tilde\alpha_r}=q^2$ and
$g_{\tilde\alpha_r}=q^{2(n-r)}$.

\medskip
\noindent
{\bf Type DI}
{\bf{\it Case 2}}:  $\gg$ is of type $D_n$, $n\geq 4$. 

\noindent
$\Theta(\alpha_i)=-\alpha_i$ for $1\leq i\leq n-2$, 
$\Theta(\alpha_{n-1})=-\alpha_{n}, \Theta(\alpha_{n})=-\alpha_{n-1}$.
$\Sigma=\{\tilde\alpha_1,\dots, \tilde\alpha_{n-1}\}$ is of type 
$B_{n-1}$.

\noindent
minuscule weight: $\tilde\mu=\omega'_{n-1}$ where $\mu=\omega_{n-1}$.

\noindent
$a_{\tilde\alpha_i}=q^4$ for $1\leq i\leq n-2$,
$a_{\tilde\alpha_{n-1}}=q^2$, and
$g_{\tilde\alpha_i}=q^2$ for $1\leq i\leq n-1$.

\medskip
\noindent
{\bf Type DI}
{\bf{\it Case 3}}: $\gg$ is of type $D_n$, $n\geq 4$. 

\noindent
$\Theta(\alpha_i)=-\alpha_i$ for $1\leq i\leq n$.
 $\Sigma=\{\alpha_1,\dots, \alpha_n\}$ is of type $D_n$.
 
 \noindent
 minuscule weight: $\tilde\mu_i=\omega_{n-i}$ where 
 $\mu_i=\omega_{n-i}$, for $i=0$ and $i=1$.
 
 \noindent
 $a_{\tilde\alpha_i}=q^4$ and
 $g_{\tilde\alpha_i}=q^2$ for all $1\leq i\leq n$.
 
 \medskip
\noindent
{\bf Type DII} 
  $\gg$ is of type $D_n$, $n\geq 4$.
  
  \noindent
$\Theta(\alpha_i)=\alpha_i$ for $2\leq i\leq n$, 
$\Theta(\alpha_1)=-\alpha_{1}-2\alpha_{2}-\cdots 
-2\alpha_{n-2}-\alpha_{n-1}-\alpha_n$. 
$\Sigma=\{\tilde\alpha_1\}$ is of type $A_1$.

\noindent
minuscule weight: $\tilde\mu=\omega'_1$ where $\mu=\omega_n$.

\noindent
$a_{\tilde\alpha_1}=q^2$ and
$g_{\tilde\alpha_1}=q^{2(n-1)}$.

\medskip
\noindent
{\bf Type DIII}
{\bf{\it Case 1}}:
 $\gg$ is of type $D_n$ where $n$ is even, 

 \noindent
 $\Theta(\alpha_i)=\alpha_i$ for $i=2j-1, 1\leq j\leq n/2$.
$\Theta(\alpha_i)=-\alpha_{i-1}-\alpha_i-\alpha_{i+1}$ for 
for $i=2j, 1\leq j\leq (n-2)/2$,
$\Theta(\alpha_n)=-\alpha_n$.
$\Sigma=\{\tilde\alpha_2,\tilde\alpha_4,\dots, \tilde\alpha_n\}$ is 
of type $C_n$.

\noindent
minuscule weight: $\tilde\mu=\omega'_2$ where $\mu=\omega_1$.

\noindent
$a_{\tilde\alpha_i}=q^2$ and $g_{\tilde\alpha_i}=q^4$for $i=2j, 1\leq j\leq (n-2)/2$,
$a_{\tilde\alpha_n}=q^4$, and 
$g_{\tilde\alpha_n}=q^2$.

\medskip
\noindent
{\bf Type E1}
 $\gg$ is of type E6.
 
 \noindent
 $\Theta(\alpha_i)=-\alpha_i$ for all $i$.
 $\Sigma=\{\tilde\alpha_1,\dots, \tilde\alpha_6\}$.
 
 \noindent
 minuscule weight: $\tilde\mu=\omega'_6$ where $\mu=\omega_6$.
 
 \noindent
 $a_{\tilde\alpha_i}=q^4$ and
 $g_{\tilde\alpha_i}=q^2$ for all $1\leq i\leq n$.

 \medskip
\noindent
{\bf Type EII}
 $\gg$ is of type E6.
 
 \noindent
 $\Theta(\alpha_i)=-\alpha_{p(i)}$ where 
$p(1)=6,p(3)=5,p(4)=4,p(2)=2,p(5)=3,$ and $p(6)=1$.
$\Sigma=\{\tilde\alpha_2,\tilde\alpha_4, \tilde\alpha_3, 
\tilde\alpha_1\}$ is of type F4.

\noindent
pseudominuscule weight: $\tilde\mu=\omega'_1$ where $\mu=\omega_1$.

\noindent
$a_{\tilde\alpha_i}=q^2$ for $i=3,1$, $a_{\tilde\alpha_i}=q^4$ for 
$i=2,4$, and
$g_{\tilde\alpha_i}=q^2$, for $i=1,2,3,4$.

\medskip
\noindent
{\bf Type EIV}
 $\gg$ is of type E6.
 
 \noindent
 $\Theta(\alpha_i)=\alpha_i$ for $i=2,3,4,5$,
$\Theta(\alpha_1)=-\alpha_1-2\alpha_3-2\alpha_4-\alpha_5-\alpha_2$,
$\Theta(\alpha_6)=-\alpha_6-2\alpha_5-2\alpha_4-\alpha_3-\alpha_2$.
$\Sigma=\{\tilde\alpha_1,\tilde\alpha_6\}$ is of type A3.

\noindent
$a_{\tilde\alpha_i}= q^2$ and $g_{\tilde\alpha_i}= q^8$ for $i=1,6$.

\medskip
\noindent
{\bf Type EV}
 $\gg$ is of type E7.
 
 \noindent
 $\Theta(\alpha_i)=-\alpha_i$ for all $i$.
 $\Sigma=\{\tilde\alpha_1,\dots, \tilde\alpha_7\}$.
 
 \noindent
 minuscule weight: $\tilde\mu=\omega'_7$ where $\mu=\omega_7$.
 
 \noindent
 $a_{\tilde\alpha_i}=q^4$ and
 $g_{\tilde\alpha_i}=q^2$ for all $1\leq i\leq 7$.

\medskip
\noindent
{\bf Type EVI}  $\gg$ is of type E7.

\noindent
$\Theta(\alpha_i)=\alpha_i$ for $i=2,5,7$,
$\Theta(\alpha_6)=-\alpha_6-\alpha_5-\alpha_7$,
$\Theta(\alpha_4)=-\alpha_2-\alpha_5-\alpha_4$,
$\Theta(\alpha_i)=-\alpha_i$ for $i=1,3$.
$\Sigma=\{\tilde\alpha_1,\tilde\alpha_3,\tilde\alpha_4, 
\tilde\alpha_6\}$ is of type F4.

\noindent
pseudominuscule weight: $\tilde\mu=\omega'_6$ where $\mu=\omega_7$.

\noindent
$a_{\tilde\alpha_i}=q^4$ and $g_{\tilde\alpha_i}=q^2$ for $i=1,3$, $a_{\tilde\alpha_i}=q^2$
and $g_{\tilde\alpha_i}=q^4$ for 
$i=4,6$.

\eject
\noindent
{\bf Type EVII} 
 $\gg$ is of type E7.
 
 \noindent
 $\Theta(\alpha_i)=\alpha_i$ for $i=2,3,4,5$, 
 $\Theta(\alpha_1)=-\alpha_1-2\alpha_3-2\alpha_4-\alpha_2-\alpha_5$,
 $\Theta(\alpha_6)=-\alpha_6-2\alpha_4-2\alpha_5-\alpha_2-\alpha_3$,
 $\Theta(\alpha_7)=-\alpha_7$.
 $\Sigma=\{\tilde\alpha_1,\tilde\alpha_6,\tilde\alpha_7\}$ is of 
 type C3.

 \noindent
 $a_{\tilde\alpha_i}=q^2$ and $g_{\tilde\alpha_i}=q^8$ for $i=1,6$, 
 $a_{\tilde\alpha_i}=q^4$ and $g_{\tilde\alpha_i}=q^2$ for 
 $i=7$.
 
\medskip
\noindent
{\bf Type EVIII}
 $\gg$ is of type E8.
 
 \noindent
 $\Theta(\alpha_i)=-\alpha_i$ for all $i$.
 $\Sigma=\{\tilde\alpha_1,\dots, \tilde\alpha_8\}$ is of type E8.
 
 \noindent
 pseudominuscule weight: $\tilde\mu=\omega'_8$ where $\mu=\omega_8$.
 
 \noindent
 $a_{\tilde\alpha_i}=q^4$ and
 $g_{\tilde\alpha_i}=q^2$ for all $1\leq i\leq 8$.

\medskip
\noindent
{\bf Type EIX} 
 $\gg$ is of type E8.
 
 \noindent
 $\Theta(\alpha_i)=\alpha_i$ for $i=2,3,4,5$, 
 $\Theta(\alpha_1)=-\alpha_1-2\alpha_3-2\alpha_4-\alpha_2-\alpha_5$,
 $\Theta(\alpha_6)=-\alpha_6-2\alpha_4-2\alpha_5-\alpha_2-\alpha_3$,
 $\Theta(\alpha_i)=-\alpha_i$ for $i=7,8$.
 $\Sigma=\{\tilde\alpha_8,\tilde\alpha_7,\tilde\alpha_6,\tilde\alpha_1\}$ is of 
 type F4.

 \noindent
 $a_{\tilde\alpha_i}=q^2$ and $g_{\tilde\alpha_i}=q^8$ for $i=1,6$, 
 $a_{\tilde\alpha_i}=q^4$ and $g_{\tilde\alpha_i}=q^2$ for 
 $i=7,8$.
 
 \medskip
\noindent
{\bf Type FI} 
 $\gg$ is of type F4.
 
 \noindent
 $\Theta(\alpha_i)=-\alpha_i$  for $i=1,2,3,4$.
 $\Sigma=\{\tilde\alpha_1,\tilde\alpha_2,\tilde\alpha_3,\tilde\alpha_4\}$ 
 is of type F4.
 
 \noindent
 pseudominuscule weight: $\tilde\mu=\omega'_4$ where $\mu=\omega'_4$.
 
 \noindent
 $a_{\tilde\alpha_i}=q^4$ and $g_{\tilde\alpha_i}=q^2$ for $i=1,2$, 
 $a_{\tilde\alpha_i}=q^2$ and $g_{\tilde\alpha_i}=q$ for 
 $i=3,4$.
 
 \medskip
\noindent
{\bf Type G}
$\gg$ is of type G2.
 
\noindent
 $\Theta(\alpha_i)=-\alpha_i$  for $i=1,2$.
 $\Sigma=\{\tilde\alpha_1,\tilde\alpha_2\}$ 
 is of type G2.
 
 \noindent
 pseudominuscule weight: $\tilde\mu=\omega'_1$ where $\mu=\omega'_1$.
 
 \noindent
 $a_{\tilde\alpha_i}=q^4$ and $g_{\tilde\alpha_i}=q^2$  for $i=1$, 
 $a_{\tilde\alpha_i}=q^{12}$ and $g_{\tilde\alpha_i}=q^{18}$ for 
 $i=2$.

\bigskip
\centerline{REFERENCES}
\medskip
\medskip
\noindent   [A] S. Araki, On root systems and an infinitesimal 
classification of irreducible symmetric spaces, {\it Journal of 
Mathematics, Osaka City University} {\bf 13} (1962), no. 1, 1-34.

\medskip
\noindent [D] M.S. Dijkhuizen, Some remarks on the construction
of quantum symmetric spaces, In: {\it Representations of Lie
Groups, Lie Algebras and Their Quantum Analogues, } 
Acta Applicandae Mathematicae {\bf 44} (1996), no. 1-2, 59-80.

\medskip
\noindent [DS] M.S. Dikhhuizen and J.V. Stokman, Some limit transitions between BC type orthogonal polynomials interpreted on
quantum complex Grassmannians, {\it Publ. Res. Inst. Math. Sci. 35} (1999), 451-500.

\medskip
\noindent [DN] M.S. Dijkhuizen and M. Noumi, A family of quantum
projective spaces and related $q$-hypergeometric orthogonal
polynomials, {\it Transactions of the A.M.S.} {\bf 350} (1998), no.
8, 3269-3296.

\medskip\noindent [H] S. Helgason, 
Fundamental solutions of invariant differential operators on symmetric spaces, 
{\it American Journal of Mathematics} {\bf 86 } (1964) 565--601.

\medskip
\noindent [Hu] J.E. Humphreys, {\it Introduction to Lie Algebras
and Representation Theory}, Springer-Verlag, New York
(1972).

\medskip
\noindent [J1] A. Joseph, {\it Quantum Groups and Their Primitive
Ideals}, Springer-Verlag, New York (1995).

\medskip
\noindent [J2] A. Joseph, A completion of the quantized enveloping algebra of a Kac-Moody algebra,
{\it  Journal of Algebra} {\bf 214} (1999), no. 1,
235-275. 

\medskip
\noindent [J3] A. Joseph, 
On a Harish-Chandra homomorphism,
{\it C. R. Acad. Sci. Paris SŽr. I Math.} {\bf 324} (1997), no. 7, 759-764.

\medskip \noindent [JL1] A. Joseph and G. Letzter, Local
finiteness of the adjoint action for quantized enveloping
algebras, {\it Journal of Algebra} {\bf 153} (1992), 289 -318.

\medskip
\noindent [JL2] A. Joseph and G. Letzter, Separation of variables
for quantized enveloping algebras, {\it American Journal of
Mathematics} {\bf 116} (1994), 127-177.

\medskip
\noindent [JL3]  A. Joseph and G. Letzter, Evaluation of the quantum affine PRV determinants,
{\it Mathematical Research Letters} {\bf 9}  (2002), 307-322.

\medskip
\noindent [K1] T.H.  Koornwinder, Askey-Wilson polynomials as zonal 
spherical functions on the SU(2) quantum group, {\it SIAM Journal on 
Mathematical Analysis} {\bf 24} (1993), no. 3, 795-813.

\medskip
\noindent [K2] T.H. Koornwinder, Askey-Wilson polynomials for root 
systems of type BC, In: {\it Hypergeometric functions on domains of 
positivity, Jack polynomials and applications}, Contemporary 
Mathematics {\bf 138} (1992), 189-204. 

\medskip
\noindent [Ke] M.S. K\'eb\'e, ${\cal O}$-alg\`ebres quantiques, {\it C. R. Acad\'emie des Sciences,
S\'erie I, Math\'ematique } {\bf 322} (1996), no. 1, 1-4.

\medskip
\noindent
[Ki]ÊA.A. Kirillov, Jr., Lectures on affine Hecke algebras and Macdonald's conjectures, {\it Bulletin of the A.M.S., 
New Series}  {\bf 34} (1997), no. 3, 251-292. 

\medskip
\noindent [L1] G. Letzter, Subalgebras which appear in quantum
Iwasawa decompositions, {\it Canadian Journal of Mathematics}
  {\bf 49} (1997), no. 6, 1206-1223.
   
\medskip \noindent [L2] G. Letzter,  Symmetric pairs for quantized
enveloping algebras, {\it Journal of Algebra} {\bf 220} (1999), no. 2,
729-767.

\medskip \noindent [L3] G. Letzter, Coideal subalgebras and quantum 
symmetric pairs, In: {\it New Directions in Hopf Algebras, MSRI publications}
{\bf 43}, Cambridge University Press (2002), 117-166.

\medskip\noindent [L4]  G. Letzter, Quantum symmetric pairs and their 
zonal spherical functions, {\it Transformation Groups}, to appear.

\medskip
\noindent
[M1] I.G. Macdonald,  
{\it Spherical functions on a group of $p$-adic type},
  Publications of the Ramanujan Institute  {\bf 2},
 Madras (1971).

\medskip
\noindent
[M2] I.G. Macdonald,
Orthogonal polynomials associated with root systems, {\it S\'eminaire Lotharingien de Combinatoire} 
{\bf 45} (2000/01) 40 pp.

 \medskip
\noindent [N] M. Noumi, Macdonald's symmetric polynomials as zonal
spherical functions on some quantum homogeneous spaces, {\it
Advances in Mathematics} {\bf 123} (1996), no. 1, 16-77.

\medskip
\noindent [NS] M. Noumi and T. Sugitani, Quantum symmetric spaces
and related q-orthogonal polynomials, In: {\it Group Theoretical
Methods in Physics (ICGTMP)} (Toyonaka, Japan, 1994), World
Science Publishing, River Edge, New Jersey (1995), 28-40.

\medskip
\noindent [NDS] M. Noumi, M.S. Dijkhuizen, and T. Sugitani, Multivariable Askey-Wilson polynomials and quantum complex
Grassmannians, {\it Fields Institute Communications} {\bf 14} (1997), 167-177.

\medskip
\noindent [R]  M. Rosso, Groupes Quantiques, Repr\'esentations
lin\'eaires et applications, Th\`ese, Paris 7  (1990).

\medskip
\noindent [S] T. Sugitani, Zonal spherical functions on quantum Grassmann manifolds, {\it J. Math. Sci. Univ. Tokyo}
{\bf 6} (1999), no. 2, 335-369.

\medskip
\noindent 
[W] ÊN.R.ÊWallach, {\it Real Reductive Groups, I,}
 Pure and Applied Mathematics {\bf 132} Boston
(1988).

\end{document}